\documentclass[12pt]{amsart}
\usepackage{amsmath,amssymb,amsthm,mathrsfs,amsfonts}
\newtheorem{thm}{Theorem}[section]
\newtheorem{lemma}[thm]{Lemma}
\newtheorem{defn}[thm]{Definition}

\newtheorem{cor}[thm]{Corollary}
\newtheorem{prop}[thm]{Proposition}
\newtheorem{assump}[thm]{Assumptions}
\newtheorem{problem}[thm]{Problem}
\theoremstyle{remark}

\numberwithin{equation}{section}
\newcommand{\bbmatrix}[1]
{\left[ \begin{array}{cccccccccc} #1 \end{array} \right]}
\newcommand{\bU}{{\mathbf U}}

\newcommand{\CI}{\boldsymbol{\mathcal{SI}}}
\newcommand{\RCI}{\boldsymbol{\mathcal{RSI}}}
\newcommand{\SC}{\boldsymbol{\mathcal{S}}}
\newcommand{\RSC}{\boldsymbol{\mathcal{RS}}}
\newcommand{\TR}{{\operatorname{tr}}}
\title[On the class $\CI$]
{On the class $\CI$ of $J$-contractive functions intertwining solutions of linear
differential equations}
\author[D. Alpay]{Daniel Alpay}
\author[A. Melnikov]{Andrey Melnikov}
\author[V. Vinnikov]{Victor Vinnikov}
\thanks{D. Alpay thanks the
Earl Katz family for endowing the chair which
supported his research. The research of the authors was
supported in part by the
Israel Science Foundation grant 1023/07}
\address{Department of Mathematics \newline
Ben Gurion University of the Negev \newline
P.O.B. 653, \newline
Be'er Sheva 84105, \newline
ISRAEL }
\email{\newline dany@math.bgu.ac.il,\newline andreym@math.bgu.ac.il, \newline
vinnikov@math.bgu.ac.il}

\begin{document}
\begin{abstract}
In the PhD thesis of the second author under the supervision of the third author was defined
the class $\CI$ of $J$-contractive functions, depending on a parameter and arising as transfer functions of
overdetermined conservative $2D$ systems invariant in one direction.
In this paper we extend and solve in the class $\CI$, a number of problems originally set for the class
$\SC$ of functions contractive in the open right-half plane, and
unitary on the imaginary line with respect to some preassigned
signature matrix $J$. The problems we consider include the Schur algorithm, the
partial realization problem and the Nevanlinna-Pick interpolation problem. 
The arguments rely on a correspondence between elements in a given subclass of $\CI$ and elements in $\SC$.
Another important tool in the arguments is a new result pertaining to the classical tangential Schur algorithm.
\end{abstract}

\subjclass{Primary:47A48, 93C15, 93C35; Secondary: 16E45, 46E22}

\keywords{Schur algorithm, time-invariant $2D$-systems, linear
differential equations}
\maketitle
\tableofcontents

\section{Introduction}
Functions $S(\lambda)$, which are ${\mathbb C}^{p\times p}$-valued,
analytic and contractive in the open right half plane $\mathbb C_+$,
or equivalently, such that the kernel
\[ 
K(\lambda,w) = \dfrac{I_p - S(w)^* S(\lambda)}{\lambda+w^*}
\]
is positive\footnote{positive-definite in the classical terminology}
in $\mathbb C_+$, play an important role in system theory,
inverse scattering theory, network theory and related topics; see for
instance \cite{bib:livsic}, \cite{bib:bc}, \cite{bib:dd}, \cite{bib:helton},
\cite{bib:AlpSchur}. 
Here, positivity of the kernel means that for every $n\in\mathbb N$ and every
choice of points $w_1,\ldots,w_n \in \mathbb C_+$ and vectors 
$\xi_1,\ldots,\xi_n\in\mathbb C^{1\times p}$ the $n\times n$ Hermitian matrix
\[ \big[ \xi_i K(w_i,w_j) \xi_j^* )\big]_{i,j=1,\ldots n}
\]
is positive (that is, has all its eigenvalues greater or equal to $0$).\\

Far reaching generalizations of this class were introduced in
\cite{bib:MyThesis, bib:MelVin1, bib:MelVinC}, in the study
of $2D$-linear systems (say, with
respect to the variables $(t_1,t_2)$), invariant with respect to the
variable $t_1$. To introduce the classes defined in these papers we
first need a definition.

\begin{defn}
\label{def:VesPar}
Let $\sigma_1$, $\sigma_2$, $\gamma$ and $\gamma_*$ be
${\mathbb C}^{p\times p}$-valued
functions, continuous on an interval ${\mathrm I}=[a,b]$.
Suppose moreover that $\sigma_1$ and $\sigma_2$ take self-adjoint
values, and that $\sigma_1$ is differentiable and invertible on
${\mathrm I}$, and that the following relations hold:

\[ \gamma(t_2) + \gamma(t_2)^* = \gamma_*(t_2) +
\gamma_*(t_2)^* = -\dfrac{d}{dt_2} \sigma_1(t_2),\quad t_2\in \mathrm I.
\]
Then
$\sigma_1,\sigma_2, \gamma, \gamma_*$
and the interval $\mathrm I$ are called \textbf{vessel parameters}.
\end{defn}

The class of functions $\CI$ corresponding to thr vessel parameters $\sigma_1$, $\sigma_2$,
$\gamma$, $\gamma_*$  and the interval ${\mathrm I}$
was introduced in \cite{bib:MyThesis, bib:MelVinC} (see Definition \ref{def:CI} below)
and consists of
the functions $S(\lambda,t_2)$ of two variables $\lambda,t_2$ such
that for every $t_2\in {\mathrm I}$ the function $S(\lambda, t_2)$ is
meromorphic in ${\mathbb C}_+$ and  the kernel
\begin{equation}\label{eq:PosS}
\dfrac{\sigma_1(t_2) - S(w,t_2)^* \sigma_1(t_2)S(\lambda,t_2) }{\lambda+w^*}
\end{equation}
is positive for $\lambda$ and $\omega$ in the domain of analyticity
of $S(\lambda, t_2)$ in $\mathbb C_+$. For positive
$\sigma_1(t_2)$, 
the positivity of the kernel implies that $S$ is analytic in
${\mathbb C}_+$; see \cite{bib:donoghue}, 
\cite{bib:AlpSchur}. For general (invertible) $\sigma_1(t_2)$, the
entries of $S$ are of bounded type and $S$ has at most poles in ${\mathbb
  C}_+$; see \cite{bib:adrs}.
It is also required that $S(\lambda,t_2)$ is analytic at
infinity for each $t_2$, with value $I_p$ there, and that $S(\lambda,t_2)$ maps
solutions of the \textit{input} Linear Differential Equation (LDE)
with the spectral parameter $\lambda$
\[ \lambda \sigma_2(t_2) u(\lambda,t_2) - \sigma_1(t_2) \frac{\partial}{\partial t_2}u(\lambda,t_2) +
\gamma(t_2) u(\lambda,t_2) = 0,
\]
to solutions of the \textit{output} LDE
\[ \lambda \sigma_2(t_2) y(\lambda,t_2) - \sigma_1(t_2) \frac{\partial}{\partial t_2}y(\lambda,t_2) +
\gamma_*(t_2) y(\lambda,t_2) = 0.
\]

It is proved in \cite{bib:MyThesis, bib:MelVinC} that elements of
$\CI$
are the transfer functions
of $t_1$-invariant conservative $2D$ systems; see Section \ref{sec:2Dsys}.\\

The purpose of this paper is to study various questions in
the class $\CI$ such as Nevanlinna-Pick interpolation, property of moments, etc. 
A key result is the following theorem, which we prove in the sequel; see Section
\ref{sec:Schur}.

\noindent\textbf{Theorem \ref{thm:gamma*}} 
\textit{Let us fix the parameters
  $\sigma_1,\sigma_2$, and $\gamma$, and the interval
$\mathrm I$. Then for every
  $t_2^0\in{\mathrm I}$
there is a one-to-one correspondence between pairs
$(\gamma_*, S)$ such that $S\in\CI$ and $\gamma_*$ continuous in a
neighborhood of $t_2^0$, and 
functions $Y(\lambda)$, meromorphic in ${\mathbb C}_+$, and
with the following properties
\begin{enumerate}
\item $Y(\infty)=I_p$,
\item $Y(\lambda)^* \sigma_1(t_2^0) Y(\lambda) \leq \sigma_1(t_2^0)$
for $\lambda\in\mathbb C_+$ where $Y$ is analytic, and
\item   $Y(\lambda)^* \sigma_1(t_2^0) Y(\lambda) = \sigma_1(t_2^0)$
  for almost $\lambda$ satisfying $\Re\lambda=0$ and where $Y(\lambda)$ is the non-tangential limit.
\end{enumerate}
}
As mentioned above the $\sigma_1(t_2)$-contractivity of $Y$ implies that
$Y$ is of bounded type in ${\mathbb C}_+$, and thus the asserted
non-tangential limits exist almost everywhere.

\begin{defn}
The class of functions $Y$ with the properties in Theorem \ref{thm:gamma*} will be
will be denoted by $\SC(t_2^0)$ and the functions will be called $\sigma_1(t_2^0)$-inner.
The subclass of rational functions of $\SC(t_2^0)$ will be denoted
by $\RSC(t_2^0)$.
\end{defn}
For the sequel, it is important to notice that the general
tangential Schur algorithm developped in \cite{bib:ad} can be applied
to functions in $\SC(t_2^0)$, and in particular in $\RSC(t_2^0)$.

\begin{defn}
$\RCI$ will denote the subclass of functions in $\CI$ which are
  rational
in $\lambda$ for every $t_2\in \mathrm I$.
\end{defn}
The paper consists of six sections besides the introduction, and we now describe their content.
In Section \ref{sec:2Dsys} we
review the main results from \cite{bib:MyThesis} (see details in \cite{bib:MelVin1}
and \cite{bib:MelVinC}) on $t_1$-invariant conservative $2D$-systems,
relevant to the present work. In particular the class
$\CI$ mentioned above
consists of the transfer functions of these systems.
In Section \ref{sec:SchurFunc} we present the reproducing kernel
space approach to the tangential Schur algorithm for the class $\SC$, as developed in \cite{bib:ad}.
We obtain in particular new formulas which allow us to find the main
operator in a realization of an element of $\SC$ after one iteration
of the tangential Schur algorithm; see formulas
\eqref{eq:BS},  \eqref{eq:AS},  \eqref{eq:XS} in Theorem
\ref{thm:real}. In Section \ref{sec:Schur}
we develop the tangential Schur algorithm for a function
$S(\lambda, t_2)\in\CI$.
Applying directly the theory of the previous section to
$S(\lambda, t_2)$ leads to a new function which need not belong to
$\CI$. Instead, we apply the Schur algorithm
to the $\sigma_1(t_2^0)$-inner function $S(\lambda, t_2^0)$ for
some preassigned $t_2^0\in{\mathrm I}$, and obtain a simpler (in terms
of McMillan degree) $\sigma_1(t_2^0)$-inner function $S_0(\lambda, t_2^0)$.
We use Theorem \ref{thm:gamma*} to obtain an element in a class
$\CI$ from $S_0(\lambda, t_2^0)$. We call this
procedure the tangential Schur algorithm for the class
$\CI$.
We study in Section \ref{sec:MarkovPart} the
coefficients (called the \textit{Markov moments})
$H_i(t_2)$ of the expansion of $S(\lambda,t_2)$ around $\lambda=\infty$
\begin{equation} \label{eq:SMoments}
 S(\lambda,t_2) = I_p -
\sum\limits_{i=0}^\infty \frac{1}{\lambda^{i+1}} H_i(t_2)
\end{equation}
It turns out that the first Markov moment $H_0(t_2)$ satisfies the
Lyapunov equation
\begin{equation} \label{eq:H0}
\gamma_*(t_2) - \gamma(t_2) =
\sigma_2(t_2) H_0(t_2)  -
\sigma_1(t_2) H_0(t_2) \sigma_1^{-1}(t_2)\sigma_2(t_2),\quad
t_2\in\mathrm I,
\end{equation}
which means that given the functions $\sigma_1, \sigma_2,
\gamma$ and $H_0$ on $\mathrm I$, one
can uniquely reconstruct $\gamma_*$, and, as a result of Theorem
\ref{thm:gamma*}, there will exists a unique function $S(\lambda,t_2)$
($t_2\in
\mathrm I$) with the given Markov parameters.
We prove the following theorem \ref{thm:HiRstrct} on the structure of Markov moments:

\noindent\textbf{Theorem} 
\textit{
Suppose that we are given moments $H_i(t_2)$, defined in a neighborhood of the 
point $t_2^0\in\mathrm I$. Then there exists $n_0 \le p^2$ 
such that each element $H_{i+1}$ is uniquely determined from $H_0,
\ldots, H_i$ using the algebraic
formulas \eqref{eq:Moments}, and $n_0$ LDEs with arbitrary initial
conditions, 
obtained from \eqref{eq:HiHi+1}. Moreover, $n_0$ equations must be satisfied by the elements
of $\gamma_*(t_2)$.}

Starting from the Markov coefficients at infinity, constructed under conditions of Theorem \ref{thm:HiRstrct}
we recover in Section \ref{sec:InvKrein} the transfer function $S(\lambda,t_2)$. 
In the arguments we make use of a Krein-space realization theorem
of Dijksma--Langer--de Snoo for analytic functions at infinity \cite[Theorem 3]{bib:KreinReal}. 
By a counterexample, we show that
it is not always possible to reconstruct the function using Hilbert space realizations.

In Section \ref{sec:NevPick} we study two generalizations of Nevanlinna Pick
interpolation problem in the class $\CI$. Let us recall that in the classical Nevanlinna Pick interpolation problem 
\cite{bib:Nev, bib:Pick} Schur analysis plays a special role; see for instance \cite{bib:fk}, \cite{bib:dym_cbms}, 
\cite{bib:AlpSchur}. Schur analysis gives a parametrization
of all solutions (when they exist) for the given data. The first problem \ref{prob:NPclassic} is similar to 
the non commutative and Riemann surface cases and involves specifying the
exact class of input-output mappings, and
a finite number of inputs that are to be mapped to the corresponding outputs. In our case this yields fixing
vessels parameters and a set of values $w_i, i=1,\ldots,N$ and corresponding solutions
$\xi_i(t_2)$ and $\eta_i(t_2)$ of the input and the output LDEs, respectively:
\begin{problem}[Nevanlinna-Pick interpolation]
\label{prob:NPclassic}
Let $\sigma_1, \sigma_2, \gamma, \gamma_*$ be vessel parameters and let $\mathrm I$ be an interval. Let
$N\in\mathbb N$ and $w_i, i=1,\ldots,N$ be complex numbers. Suppose also that $N$ input functions
$\xi_i(t_2)$ satisfying \eqref{eq:InCC} with corresponding spectral parameters $w_i$'s, and
$N$ output functions $\eta_i(t_2)$ satisfying \eqref{eq:OutCC} with the given $w_i$'s.
\begin{enumerate}
\item Give sufficient and necessary conditions, so that there exists
$S(\lambda,t_2) \in \CI$ such
that $S(w_j,t_2)\xi_i(t_2) = \eta_i(t_2)$, $j=1,\ldots, N$, 
on a sub-interval of $\mathrm I$.
\item Describe the set of all solutions for this problem.
\end{enumerate} 
\end{problem}
The second problem uses the fact that we can also specify the data
for different values of $t_2$ and is more similar to the classical one:
\begin{problem} \label{prob:NPhard}
Given ${\mathbb C}^{p\times p}$-valued functions $\sigma_1,\sigma_2,
\gamma$ defined on $\mathrm I$,
and given
$N$ quadruples $<t_2^j,w_j,\xi_j,\eta_j>$, where
$t_2^j\in\mathbb \mathrm I, w_j\in\mathbb C_+$,
$\xi_j,\eta_j\in\mathbb
C^{1\times p}$ $j=1,\ldots N$, then:
\begin{enumerate}
\item Give sufficient and necessary conditions, so that there
exists $\gamma_*$ and $S(\lambda,t_2) \in \CI$ such
that $S(w_j,t_2^j)\xi_j=\eta_j$, $j=1,\ldots, N$, 
on a sub-interval of $\mathrm I$
containing all the $t_2^j$.
\item Describe the set of all solutions for this problem.
\end{enumerate}
\end{problem}

If all the values $t_2^j=t_2^0$ are equal, 
we have to find a
function $S(\lambda,t_2^0)$ satisfying $S(w_j,t_2^0) \xi_j=\eta_j$.
Thus, the above problem is a generalization of the classical
Nevanlinna-Pick interpolation problem.
We also remark that we do not address the question of describing the
set of all solutions.

\textbf{Remarks:} Some of the results presented here have been announced in \cite{bib:amv}.
\section{
\label{sec:2Dsys}
$t_1$ invariant conservative $2D$ systems.}
The material in this section is taken from \cite{bib:MyThesis, bib:MelVin1, bib:MelVinC}, 
where proofs and more details can be found. 
The origin of this theory can be found in the paper \cite{bib:Vortices}.

\subsection{Definition}
An overdetermined conservative $t_1$-invariant
$2D$ system is a linear input-state-output (i/s/o) system,
which consists of operators depending only on the variable $t_2$ and
is of the following form:
\begin{equation} \label{eq:systemCons}
    I\Sigma: \left\{ \begin{array}{lll}
        \frac{\partial}{\partial t_1}x(t_1,t_2) = A_1(t_2) x(t_1,t_2) + \widetilde B_1(t_2) u(t_1,t_2) \\
        x(t_1, t_2) = F (t_2,t_2^0) x(t_1, t_2^0) + \int\limits_{t_2^0}^{t_2} F(t_2, s)  \widetilde B_2(s) u(t_1, s)ds \\
        y(t_1, t_2) = u(t_1, t_2) -  \widetilde B(t_2)^* x(t_1, t_2),
\end{array} \right. 
\end{equation}
where the variable $t_1$ belongs to $\mathbb R$,  and the
variable $t_2$ belongs to some interval $\mathrm I$. Furthermore,
the input $u(t_1,t_2)$ and the output $y(t_1,t_2)$ take values
in some Hilbert space ${\mathcal E}$ and the state
$x(t_1,t_2)$ takes values in the Hilbert space ${\mathcal H}_{t_2}$.
We assume that $u(t_1,t_2)$ and $y(t_1,t_2)$ are continuous functions of
each variable when the other variable is fixed. The operators of
the system are supposed to satisfy the following:
\begin{assump}
\label{assump:reg}
\mbox{}\\
\begin{enumerate}
\item $A_1(t_2):\mathcal H_{t_2}
 \rightarrow \mathcal H_{t_2}$,
and $ \widetilde B(t_2): \mathcal E \rightarrow \mathcal H_{t_2}$ are bounded
       operators for all $t_2$,
\item The functions $\sigma_1, \sigma_2, \gamma, \gamma_*:
\mathcal E \rightarrow \mathcal E$, are continuous in the operator
norm topology.
\item $\sigma_1(t_2)$ is an invertible operator for every
  $t_2\in\mathrm I$.
\item $F(t,s)$ is an evolution continuous semi-group.
\end{enumerate}
\end{assump}
For continuous inputs $u(t_1,t_2)$, the inner state is continuously
differentiable. Requiring now the invariance of the system transition from
$(t_1^0,t_2^0)$ to $(t_1,t_2)$ via the points $(t_1^0,t_2)$ and $(t_1,t_2^0)$ respectively, is equivalent to 
the equality of second order partial derivatives of $x(t_1,t_2)$:
\begin{equation}
\label{eq:2pe}
 \frac{\partial^2}{\partial t_1\partial t_2} x(t_1,t_2) = 
\frac{\partial^2}{\partial t_2\partial t_1} x(t_1,t_2).
\end{equation}
Substituting in this equality the system equations we obtain that for
the free evolution $u(t_1,t_2)=0$ the so called Lax equation
holds
\begin{equation}\label{eq:LaxCond} 
A_1(t_2) = F(t_2, t_2^0) A_1(t_2^0) F(t_2^0, t_2).
\end{equation}
Inserting \eqref{eq:LaxCond} into \eqref{eq:2pe} we see
that the input $u(t_1,t_2)$ has to satisfy the following PDE
\[ \begin{array}{lll}
\widetilde B(t_2)\sigma_2(t_2) \frac{\partial}{\partial t_1}u(t_1,
t_2) -  
\widetilde B(t_2)\sigma_1(t_2) \frac{\partial}{\partial t_2} u(t_1, t_2) - \\
\big( A_1(t_2) \widetilde B(t_2) \sigma_2(t_2) + F(t_2, t_2^0)\frac{\partial}{\partial t_2} [F (t_2^0, t_2) \widetilde B(t_2)\sigma_1(t_2)] \big) u(t_1, t_2) = 0.
\end{array} \]
Assuming the existence of a function $\gamma(t_2)$ satisfying
\begin{equation} 
\label{eq:OverDetCondIn}
A_1(t_2)\widetilde  B(t_2)\sigma_2(t_2) + 
F(t_2, t_2^0)\frac{\partial}{\partial s} [F (t_2^0, t_2) 
\widetilde B(t_2)\sigma_1(t_2)] = - \widetilde B(t_2) \gamma(t_2)
\end{equation}
we obtain that it is enough that $u(t_1,t_2)$ satisfies the PDE
\begin{equation}
\label{eq:InCCP}
 \sigma_2(t_2) \frac{\partial}{\partial t_1}u(t_1,t_2) -
\sigma_1(t_2) \frac{\partial}{\partial t_2}u(t_1,t_2) +
\gamma(t_2) u(t_1,t_2) = 0.
\end{equation}
The output $y(t_1,t_2)$ should satisfy the 
\textit{output compatibility condition} of the same type as for the input compatibility
condition (\ref{eq:InCCP}), namely:
\begin{equation} \label{eq:OutCCP}
\sigma_{2}(t_2) \frac{\partial}{\partial t_1}y(t_1, t_2) -
  \sigma_{1}(t_2)
\frac{\partial}{\partial t_2}y(t_1,t_2) + \gamma_*(t_2) y(t_1,t_2) = 0.
\end{equation}
Inserting here $y(t_1, t_2) = u(t_1,t_2) - \widetilde B(t_2)^*
x(t_2,t_2)$ we obtain that
\begin{eqnarray}
\label{eq:OverDetCondOut}
0 &=& \sigma_2(t_2) \widetilde B(t_2)^* A_1(t_2) F(t_2,t_2^0)-\\
&\hspace{5mm} &\nonumber
  - \sigma_1(t_2) \frac{\partial}{\partial t_2}[\widetilde
B(t_2)^* F(t_2,t_2^0)] + \gamma_*(t_2) \widetilde B(t_2)^* F(t_2,t_2^0) \\
\label{eq:LinkCond}
\hspace{3mm}\gamma(t_2) &=& \sigma_1(t_2)\widetilde B(t_2)^* \widetilde B(t_2)
\sigma_2(t_2) -\\
&\hspace{5mm}&
\nonumber
 -\sigma_2(t_2) \widetilde B(t_2)^* \widetilde B(t_2) \sigma_1(t_2) + \gamma_*(t_2).
\end{eqnarray}
The fact that the system is lossless comes from the requirement of the
so called \textit{energy balance} equations:
\[
\begin{split}
\dfrac{\partial}{\partial t_i} \langle x(t_1,t_2),x(t_1,t_2)
\rangle_{\mathcal H_{t_2}} + \langle\sigma_i(t_2)
y(t_1,t_2),y(t_1,t_2)
\rangle_\mathcal{E}=\\
&\hspace{-5cm} = \langle\sigma_i(t_2)
u(t_1,t_2),u(t_1,t_2)\rangle_\mathcal{E},~~~~~ i=1,2,
\end{split}
\]
which means that the energy of the output is distributed
between the energy of the input and the change of the energy of the
state of the system. Immediate consequences of this requirement are
\begin{eqnarray} \label{eq:ColCond}
0&=A_1(t_2) + A_1^*(t_2) + \widetilde B(t_2) \sigma_1(t_2)
\widetilde B(t_2)^* , \\
\label{eq:ColCond2} \hspace{1cm}
\dfrac{d}{dt_2} [F^*(t_2,t_2^0) F(t_2,t_2^0)] &= F^*(t_2,t_2^0)
\widetilde B(t_2)^* \sigma_2(t_2) \widetilde B(t_2)F(t_2,t_2^0).
\end{eqnarray}
In this manner we obtain the notion of \textit{conservative vessel in
  the integral form}, which is a collection of operators and spaces
\[ \mathfrak{V} =
(A_1(t_2), F(t_2,t_2^0),
\widetilde B(t_2); \sigma_1(t_2), \sigma_2(t_2), \gamma(t_2), \gamma_*(t_2);
\mathcal H_{t_2}, \mathcal{E})
\]
where the operators
satisfy the regularity assumptions \ref{assump:reg}, 
and the following vessel conditions:
\[ \begin{array}{llllllll}
   0 = A_1(t_2) + A_1^*(t_2) + \widetilde B(t_2)^*
\sigma_1(t_2) \widetilde B(t_2)  & (\text{\ref{eq:ColCond}}) \vspace{2mm} \\
    \| F(t_2,t_2^0) x(t_1,t_2^0)\|^2 - \| x(t_1,t_2^0)\|^2 = \\
\hspace{2cm} = \int_{t_2^0}^{t_2} \langle \sigma_2(s)
\widetilde B(s) x(t_1, s), \widetilde B(s) x(t_1, s)
\rangle ds & (\text{\ref{eq:ColCond2}}) \vspace{2mm}  \\
F(t_2,t_2^0) A_1(t_2^0) = A_1(t_2) F(t_2,t_2^0)  & (\text{\ref{eq:LaxCond}}) \vspace{2mm}  \\
0 =  \frac{d}{dt_2} (F (t_2^0, t_2) \widetilde B(t_2) \sigma_1(t_2)) + \\
\hspace{2cm} + F(t_2^0, t_2) A_1(t_2) \widetilde B(t_2) \sigma_2(t_2) 
+ F(t_2^0, t_2) \widetilde B(t_2) \gamma(t_2)
     & (\text{\ref{eq:OverDetCondIn}}) \vspace{2mm}  \\
0 = \sigma_1(t_2) \frac{\partial}{\partial t_2}[\widetilde B(t_2)^* F(t_2,t_2^0)] - \\
\hspace{2cm} - \sigma_2(t_2) \widetilde B(t_2)^* A_1(t_2) F(t_2,t_2^0) - \gamma_*(t_2) \widetilde B(t_2)^* F(t_2,t_2^0)
             & (\text{\ref{eq:OverDetCondOut}}) \vspace{2mm}  \\
\gamma(t_2) = -\sigma_2(t_2) \widetilde B(t_2)^* \widetilde B(t_2) \sigma_1(t_2) + \\
\hspace{2cm} + \sigma_1(t_2) \widetilde B(t_2)^* \widetilde B(t_2)
\sigma_2(t_2) + \gamma_*(t_2) & (\text{\ref{eq:LinkCond}})
\end{array} \]
In order to simplify some notations we introduce the following definition.
\begin{defn} Let
\[ \bU = \bU(\lambda_1,\lambda_2;t_2) = 
\sigma_2(t_2) \lambda_1 - \sigma_1(t_2) \lambda_2 + \gamma(t_2),
\]
and similarly
\[ \bU_* = \bU_*(\lambda_1,\lambda_2;t_2) = 
\sigma_2(t_2) \lambda_1 - \sigma_1(t_2) \lambda_2 + \gamma_*(t_2).
\]
\end{defn}
The vessel $\mathfrak{V}$ is naturally associated to the system (\ref{eq:systemCons})
\begin{equation*}
     \Sigma: \left\{ \begin{array}{lll}
    \frac{\partial}{\partial t_1}x(t_1,t_2) = A_1(t_2) ~x(t_1,t_2) + \widetilde B(t_2) \sigma_1(t_2) ~u(t_1,t_2) \\[5pt]
    x(t_1, t_2) = F (t_2,t_2^0) x(t_1, t_2^0) + \int\limits_{t_2^0}^{t_2} F(t_2, s) \widetilde B(s) \sigma_2(s) u(t_1, s)ds \\[5pt]
    y(t_1,t_2) = u(t_1,t_2) - \widetilde B(t_2)^*~ x(t_1,t_2).
    \end{array} \right.
  \end{equation*}
with inputs and outputs satisfying the
compatibility conditions (\ref{eq:InCCP}) and (\ref{eq:OutCCP}), i.e. satisfy:
\[ \bU(\dfrac{\partial}{\partial t_1}, \dfrac{\partial}{\partial t_2};t_2) u(t_1,t_2) = 0,
\quad \bU_*(\dfrac{\partial}{\partial t_1}, \dfrac{\partial}{\partial t_2};t_2) y(t_1,t_2) = 0 \]
The theory of such vessels, developed in \cite{bib:MyThesis,bib:MelVin1}
enables to find a more convenient form of the vessel. Denoting
$\mathcal H=
\mathcal H_{t_2^0}$,
$A_1 = A_1(t_2^0)$, $F^*(t_2,t_2^0) F(t_2,t_2^0) = \mathbb
X^{-1}(t_2)$ and
$B(t_2) = F(t_2^0,t_2)\widetilde B(t_2) $, we shall obtain the
following notion, first introduced in \cite{bib:SLVessels}.
\begin{defn} 
A (differential) conservative vessel associated to the vessel parameters
is a collection of operators and spaces
\begin{equation} \label{eq:DefV}
\mathfrak{V} = (A_1, B(t_2), \mathbb X(t_2); \sigma_1(t_2), 
\sigma_2(t_2), \gamma(t_2), \gamma_*(t_2);
\mathcal{H},\mathcal{E}), 
\end{equation}
where the operators satisfy the following vessel conditions:
\begin{align}
\label{eq:DB} 0  =  \frac{d}{dt_2} (B(t_2)\sigma_1(t_2)) + A_1 B(t_2) \sigma_2(t_2) + B(t_2) \gamma(t_2), \\
\label{eq:XLyapunov} A_1 \mathbb X(t_2) + \mathbb X(t_2) A_1^*  =   B(t_2) \sigma_1(t_2) B(t_2)^*, \\
\label{eq:DX} \frac{d}{dt_2} \mathbb X(t_2)  =  B(t_2) \sigma_2(t_2) B(t_2)^*, \\
\label{eq:Link}
\gamma_*(t_2)  =  \gamma(t_2) + \sigma_2(t_2) B(t_2)^* \mathbb X^{-1}(t_2) B(t_2) \sigma_1(t_2) - \\
\nonumber - \sigma_1(t_2) B(t_2)^* \mathbb X^{-1}(t_2) B(t_2) \sigma_2(t_2)
\end{align}
\end{defn}
It turns out that Lyapunov equation \eqref{eq:XLyapunov} is partially redundant.
\begin{lemma} \label{lemma:LyapREdund}
Suppose that $B(t_2)$ satisfies \eqref{eq:DB} and $\mathbb X(t_2)$ satisfies \eqref{eq:DX},
then if the Lyapunov equation \eqref{eq:XLyapunov}
\[ A_1 \mathbb X(t_2) + \mathbb X(t_2) A_1^*  + B(t_2) \sigma_1(t_2) B(t_2)^* = 0\]
holds for a fixed $t_2^0$, then it holds for all $t_2$. If $\mathbb X(t_2^0) = \mathbb X^*(t_2^0)$, then
$\mathbb X(t_2) = \mathbb X(t_2)^*$ for all $t_2$.
\end{lemma}
\textbf{Proof:} By differentiating the left hand side of \eqref{eq:XLyapunov}, we will obtain that 
it is zero. The derivative of $\mathbb X(t_2)$ is selfadjoint. \qed

This representation of a vessel is the most convenient
when one focuses on the notion of transfer function,
as we do in the next subsection.

\subsection{\label{sec:SepVar}Transfer function}
Performing a separation of variables as follows
\[ \begin{array}{lll}
u(t_1,t_2) = u_\lambda(t_2) e^{\lambda t_1}, \\
x(t_1,t_2) = x_\lambda(t_2) e^{\lambda t_1}, \\
y(t_1,t_2) = y_\lambda(t_2) e^{\lambda t_1},
\end{array} \]
we arrive at the notion of a transfer function. Note that
$u(t_1,t_2)$ and $y(t_1,t_2)$ satisfy PDEs, but
$u_\lambda(t_2)$ and $y_\lambda(t_2)$ are solutions of LDEs with spectral
parameter $\lambda$,
\begin{eqnarray}
\label{eq:InCC} \bU(\lambda, \dfrac{\partial}{\partial t_2};t_2) u(t_1,t_2) = 0, \\
\label{eq:OutCC}
\bU(\lambda, \dfrac{\partial}{\partial t_2};t_2) y(t_1,t_2) = 0.
\end{eqnarray}
The corresponding i/s/o system becomes
\[ \left\{ \begin{array}{lll}
    x_\lambda(t_2) = (\lambda I - A_1(t_2))^{-1}  B(t_2) \sigma_1(t_2) u_\lambda(t_2) \\
    \frac{\partial}{\partial t_2}x_\lambda(t_2) = F (t_2,t_2^0) x_\lambda(t_2^0) + \int\limits_{t_2^0}^{t_2} F(t_2, s)  \widetilde B_2(s) u_\lambda(s)ds \\
        y_\lambda(t_2) = u_\lambda(t_2) -  \widetilde B(t_2)^* x_\lambda(t_2)
    \end{array} \right.
\]
The output $y_\lambda(t_2) = u_\lambda(t_2) -\widetilde B(t_2)^* x_\lambda(t_2)$ may be found from the first i/s/o equation:
\[ y_\lambda(t_2) = S(\lambda, t_2) u_\lambda(t_2), \]
using the \textit{transfer function}
\[ \begin{array}{lll}
    S(\lambda, t_2) & = I - \widetilde B(t_2)^* (\lambda I - A_1(t_2))^{-1} \widetilde B(t_2) \sigma_1(t_2)\\
 & = I - \widetilde B(t_2)^* 
(\lambda I - F(t_2,t_2^0)A_1(t_2^0) F(t_2^0,t_2))^{-1} 
\widetilde B(t_2) \sigma_1(t_2)\\
    & = I - \widetilde B(t_2)^* F(t_2,t_2^0) 
(\lambda I - A_1(t_2^0))^{-1} F(t_2^0,t_2) \widetilde B(t_2) \sigma_1(t_2) \\
    & = I - \widetilde B(t_2)^* F^*(t_2^0,t_2) F^*(t_2,t_2^0) F(t_2,t_2^0) (\lambda I - A_1)^{-1} B(t_2) \sigma_1(t_2) \\
    & = I - B(t_2)^* \mathbb X^{-1}(t_2) (\lambda I - A_1)^{-1} B(t_2) \sigma_1(t_2)
\end{array} \]
and we obtain that 
\begin{equation} \label{eq:S}
S(\lambda, t_2) = I - B(t_2)^* \mathbb X^{-1}(t_2) (\lambda I - A_1)^{-1} B(t_2) \sigma_1(t_2).
\end{equation}
\begin{prop} 
\label{prop:PropS}
The transfer function $S(\lambda,t_2)$ defined by \eqref{eq:S} has the
following properties:
\begin{enumerate}
\item For all $t_2$, $S(\lambda, t_2)$ is an analytic function of $\lambda$
                in the neighborhood of $\infty$, where it satisfies:
                \[ S(\infty, t_2) = I_p \]
\item For all $\lambda$, $S(\lambda, t_2)$ is a continuous function of $t_2$.
\item For $\lambda$ in the domain of analyticity of $S(\lambda,t_2)$:
\begin{equation}
\label{eq:Scont}
S(\lambda, t_2)^*
\sigma_1(t_2) S(\lambda, t_2) \leq \sigma_1(t_2), \quad \Re{\lambda} >
0, 
\end{equation}
and
\begin{equation}
\label{eq:Scont2}
S(\lambda, t_2)^*
\sigma_1(t_2) S(\lambda, t_2) = \sigma_1(t_2), \quad \Re{\lambda} = 0.
\end{equation}
\item Maps solutions of the input LDE
\eqref{eq:InCC} with spectral parameter $\lambda$ to the output
LDE \eqref{eq:OutCC} with the same spectral parameter.
\end{enumerate}
\end{prop}
\textbf{Proof:}
These properties are easily checked, and follow from the definition
of $S(\lambda,t_2)$:
\[ S(\lambda, t_2) = I - \widetilde B(t_2)^* (\lambda I -
A_1(t_2))^{-1}
\widetilde B(t_2) \sigma_1(t_2). \]
The function $S(\lambda, t_2)$ is analytic
for $\lambda > \| A_1(t_2)\|$ and since all the operators are
bounded, we have $ S(\infty, t_2)
= I_p$. The second property follows from the regularity assumptions
\ref{assump:reg}. The third property follows from straightforward
calculations:
\[ \begin{array}{lllll}
S(\lambda, t_2)^* \sigma_1(t_2) S(\lambda, t_2) - \sigma_1(t_2) = \\
~~~~~~ - 2 \Re (\lambda) \sigma_1(t_2)
\widetilde B(t_2)^* (\bar\lambda I - A_1^*(t_2))^{-1}(\lambda I - A_1(t_2))^{-1} \widetilde B(t_2) \sigma_1(t_2)
\end{array} \]
Here the sign of $-\Re (\lambda)$ determines the sign of
$S(\lambda, t_2)^* \sigma_1(t_2) S(\lambda, t_2) - \sigma_1(t_2)$ and
thus the third property is obtained. The fourth property follows
directly from our construction. \qed\mbox{}\\

\noindent\textbf{Remark:} When $\dim\mathcal H <\infty$, we obtain
that $S(\lambda,t_2)$ is a rational function
of $\lambda$ for every $t_2$.\\

It is an interesting fact that also the converse of Proposition
\ref{prop:PropS} holds; see \cite{bib:MyThesis},
\cite[chapter 5]{bib:MelVinC}.
\begin{thm} For any function of two variables $S(\lambda,t_2)$,
satisfying the conditions of Proposition
\ref{prop:PropS}, there is a conservative $t_1$-invariant vessel whose
transfer function is $S(\lambda,t_2)$.
\end{thm}
We define the class of transfer functions mentioned in the
introduction as follows:
\begin{defn}[\cite{bib:MelVinC}] The class
$\CI=\CI(\bU(\lambda,\dfrac{\partial}{\partial t_2};t_2),\bU_*(\lambda,\dfrac{\partial}{\partial t_2};t_2))$
consists of functions $S(\lambda,t_2)$ of two variables, which
\begin{enumerate}
		\item are analytic in a neighborhood of $\lambda=\infty$ for all 
			$t_2$ and where it holds $S(\infty,t_2)=I_p$,
    \item are continuous as functions of $t_2$ for all $\lambda$,
    \item satisfy \eqref{eq:Scont} and \eqref{eq:Scont2} in the domain
			of analyticity of $S$,
    \item map solutions of the input LDE
(\ref{eq:InCC}) with spectral parameter $\lambda$ to the output
LDE (\ref{eq:OutCC}) with the same spectral parameter
\end{enumerate}
\label{def:CI}
\end{defn}

Recall (see \cite{bib:CoddLev}) that to every LDE
can be associated an
invertible matrix (or operator) function $\Phi(t_2,t_2^0)$, called
the {\sl fundamental solution}, which
takes value $I$ at some preassigned value $t_2^0$ and such that any other solution
$u(t_2)$ of the LDE, with initial condition $u(t_2^0) = u_0$ is of the form
\[ u(t_2) =
\Phi(t_2,t_2^0) u_0.
\]
Let  $\Phi(\lambda,t_2, t_2^0)$ and $\Phi_*(\lambda,t_2,
t_2^0)$ be the fundamental solutions of the input LDE \eqref{eq:InCC} and the
output LDE \eqref{eq:OutCC} respectively, where we have added in the
notation the dependence in $\lambda$. Then,
\begin{equation} \label{eq:SInttw}
S(\lambda, t_2) \Phi(\lambda,t_2,t_2^0) =
\Phi_*(\lambda,t_2,t_2^0) S(\lambda, t_2^0)
\end{equation}
and consequently $S(\lambda,t_2)$ satisfies the following LDE
\begin{equation}
\label{eq:DS}
\begin{array}{lll}
    \frac{\partial}{\partial t_2} S(\lambda,t_2)
    = \sigma_{1}^{-1}(t_2) (\sigma_{2}(t_2)
\lambda + \gamma_*(t_2)) S(\lambda,t_2)-\\
\quad \quad -   S(\lambda,t_2)\sigma_1^{-1}(t_2) (\sigma_2(t_2) \lambda + \gamma(t_2)).
\end{array}
\end{equation}
The symmetry conditions \eqref{eq:Scont} and \eqref{eq:Scont2}, which are a result of Lyapunov equation 
\eqref{eq:Lyap} gives us
\begin{equation} \label{eq:Symmetry}
S(\lambda,t_2) = \sigma_1^{-1}(t_2) S^{-1*}(-\bar\lambda,t_2) \sigma_1(t_2)
\end{equation}
\section{\label{sec:SchurFunc}
Schur analysis for the classical case}
One of the approaches to the tangential Schur algorithm for
a rational matrix function $S(\lambda)\in\mathbb C^{p\times p}$ is based
on the theory of reproducing kernel Hilbert spaces of the kind
introduced by de Branges and Rovnyak; see \cite{bib:dbr1},
\cite{bib:dbr2}, \cite{bib:dym_cbms} for information on these
spaces. The paper \cite{bib:ad}  considers the case of column-valued
functions. In this section we adapt the results of \cite{bib:ad} to
the case of row-valued functions. We also present a new
formula for
a realization of a Schur function after implementation of the
tangential Schur algorithm.

\subsection{Schur functions and Reproducing Kernel Hilbert spaces}
In this section $\sigma_1$ denotes a fixed self-adjoint and invertible
(but not necessarily unitary)
matrix in ${\mathbb C}^{p\times p}$. Let $S(\lambda)$ be a rational function,
$\sigma_1$-inner in the open right half plane, i.e.
\[
S(\lambda)^*\sigma_1 S(\lambda) - \sigma_1 \leq 0
\]
at all points in the domain of analyticity $\Omega(S)$ of $S$
in ${\mathbb C}_+$, and
\[
S(\lambda)^*\sigma_1 S(\lambda) - \sigma_1 = 0
\]
at all points on the imaginary axis where $S$ is defined.
Then, the kernel
\begin{equation}
\label{eq:KS}
K_S(\lambda,w) = \dfrac{\sigma_1 - S(w)^* \sigma_1 S(\lambda)}
{\bar w + \lambda}
\end{equation}
is positive for $\lambda,w\in\Omega(S)$, and the space of rational
${\mathbb C}^{1\times p}$-valued functions
\[ \mathcal H(S) = \{ \sum\limits_{i=1}^n
\alpha_i c_i K_S(\lambda,w_i) \mid
\alpha_1\in\mathbb C, w_i\in \Omega, c_i\in \mathbb C^{1\times p} \}
\]
is finite dimensional. These well-known facts can be proved  using realization
theory; see for instance \cite{bib:ag}.
Furthermore, $\mathcal H(S)$ is
the reproducing kernel Hilbert space, associated to the kernel
$K_S(\lambda,w)$. The inner product is defined by
\[ \langle c K_S(\lambda,\nu), d K_S(\lambda,w) \rangle_{\mathcal H_S} =
 \langle c K_S(w,\nu)\rangle_{C^{1\times p}} =  c K_S(w,\nu) d^*.
\]
For an
arbitrary $f(\lambda)\in\mathcal H(S)$ we have
the reproducing kernel property
\[ \langle f(\lambda), \xi K(\lambda,w) \rangle_{\mathcal H_S} = f(w) \xi^*.
\]
More generally, let now $\mathcal M$ be a finite dimensional Hilbert space of
${\mathbb C}^{1\times p}$-valued functions defined in some set
$\Omega$, and let $\{ f_1(\lambda),\ldots,f_N(\lambda)\}$
be a basis of $\mathcal M$. Let $\mathbb X\in{\mathbb C}^{p\times p}$
denote the Gram matrix with $\ell,j$ entry given by
\begin{equation} \label{eq:DefX}
{\mathbb X}_{\ell,j} =
\langle f_j(\lambda),f_\ell(\lambda)\rangle_{\mathcal M},\quad
\ell,j=1,\ldots, p.
\end{equation}
It is easily seen that the space $\mathcal M$ is a
reproducing kernel Hilbert space
with kernel given by the formula
\begin{equation}
\label{eq:kernel}
 K(\lambda,w) =
\bbmatrix{f_1(w)^* & \cdots &
f_N(w)^*} \mathbb X^{-1} \bbmatrix{f_1(\lambda) \\ \vdots \\ f_N(\lambda)}.
\end{equation}
We set
\begin{equation}
\label{eq:DefF}
F(\lambda) = \bbmatrix{f_1(\lambda)\\\vdots\\f_p(\lambda)}.
\end{equation}
Assume now that $\mathcal M$ consists of rational functions, defined
on a set $\Omega({\mathcal M})$. For $\alpha\in\Omega(\mathcal M)$
the backward-shift operator $R_\alpha$ is defined by
\[ R_\alpha(f(\lambda)) = \frac{f(\lambda)-f(\alpha)}{\lambda-\alpha}.
\]
Suppose that:
\begin{enumerate}
    \item The space $\mathcal M$ is invariant under the action of $R_\alpha$,
    \item The functions $f_i(\lambda)$ has the property that
$f_i(\infty) = 0$, i.e., $F(\infty) = 0$,
\end{enumerate}
then the function $F$ given by \eqref{eq:DefF} can be written as
\[
F(\lambda) = (\lambda I - A)^{-1} B\sigma_1,
\]
for suitably chosen
matrices $A,B$. In this special case, formula \eqref{eq:kernel}
takes the form
\begin{equation}
\label{eq:K}
K(\lambda,w) = \sigma_1 B^* (\bar w I - A^*)^{-1} \mathbb X^{-1}
(\lambda I - A)^{-1} B \sigma_1.
\end{equation}
As mentioned at the beginning of this section we are interested in
kernels of the form \eqref{eq:KS} for some $\sigma_1$-inner rational
function $S$. We now recall the characterization of these spaces, and
first note the following: equation \eqref{eq:KS}
leads to
\begin{equation}
\label{eq:KTheta}
 \frac{\sigma_1 - S(w)^* \sigma_1 S(\lambda) }{\bar w +
   \lambda} = F(w)^*\mathbb X^{-1}F(\lambda).
\end{equation}
If $S$ is analytic at infinity and
satisfies there $S(\infty) = I_p$, and letting $w\rightarrow\infty$
in this equation, we obtain the formula
\begin{equation}
\label{eq:Theta}
S(\lambda) = I_p - B^* \mathbb X^{-1} (\lambda I - A)^{-1} B \sigma_1.
\end{equation}

\begin{thm}
Let $\mathcal M$ be a finite dimensional Hilbert space of
$\mathbb C^{1\times p}$-valued rational functions, which are zero at infinity.
Suppose that $R_\alpha \mathcal M\subset \mathcal M$ for
$\alpha\in\Omega(\mathcal M)$ and let $\mathbb X$ be its Gram matrix
with respect to $F(\lambda)$.
Then $\mathcal M = \mathcal H(S)$ for $S$ defined by \eqref{eq:Theta}
if and only if the Lyapunov equation
\begin{equation} \label{eq:Lyap}
A \mathbb X + \mathbb X A^* + B\sigma_1B^* = 0
\end{equation}
holds.
\end{thm}

When the spectrum of the operator $A$ is in the open left half plane
$\Re \lambda <0$, one has
\[ \mathcal H(S) = \mathbf H_{2,\sigma_1}\ominus\mathbf H_{2,\sigma_1} S,
\]
where $\mathbf H_{2,\sigma_1}$ is the Hardy space $\mathbf H_2^p$ with
the inner product
\[ [f,g]_{\mathbf H_{2,\sigma_1}} = \langle f, g \sigma_1^{-1}
\rangle_{\mathbf H_2^p}.
\]
We set $J=\bbmatrix{-\sigma_1 & 0 \\ 0 & \sigma_1}$. Note that $J$ is
both invertible and self-adjoint, and one can define $J$-inner
rational
functions. Let
$\Theta=\bbmatrix{\Theta_{11}&\Theta_{12}\\\Theta_{21}&\Theta_{22}}$
be a $J$-inner rational function;
we introduce the linear fractional transformation
\[ T_\Theta(W) = (\Theta_{11} + W\Theta_{21})^{-1}
(\Theta_{12} + W \Theta_{22}).
\]

\begin{thm}
\label{thm:S0Exist}
Let $S$ and $\Theta$ be respectively $\sigma_1$-inner $J$-inner
rational functions.
Then there exists a $\sigma_1$-inner rational function
$W$ such that $S = T_{\Theta}(W)$ if and only if the map
\begin{equation}
\label{eq:t}
 F \mapsto F \bbmatrix{-S(\lambda)\\I_p}
\end{equation}
is a contraction from $\mathcal H(\Theta)$ to $\mathcal H(S)$.
\end{thm}

This theorem originates with the work of de Branges and Rovnyak (see
\cite[Theorem 13 p. 305]{bib:dbr1}), where it is proved in a more
general setting, and is one of the key ingredients
to the reproducing kernel approach to the Schur algorithm. The proof
is the same as for column-valued functions, and is omitted.\\

As a special case of the previous theorem we have:
\begin{cor}
\label{cor:TTheta}
Given $S\in\RSC$, $w_1,\ldots,w_n\in \mathbb C_+$, and
row vectors $\xi_1,\ldots,\xi_n\in\mathbb C^{1\times p}$, define
\[ B = \bbmatrix{-\xi_1 S(w_1)^* & \xi_1\\ \vdots & \vdots \\ -
\xi_n S(w_n)^* & \xi_n},~~~~~ A_1 = \operatorname{diag}[-w_1^*,\ldots,-w_n^*].
\]
Let $\mathcal M$ be a Hilbert space of row vectors spanned by the rows
of the matrix-valued function
\[ F(\lambda) = (\lambda I - A_1)^{-1} B J.
\]
Let $\mathbb X$ be the solution of the Lyapunov equation
\begin{equation}
\label{eq:LyapJ}
A \mathbb X + \mathbb X A^* + BJB^* = 0
\end{equation}
and assume $\mathbb X>0$.
Let $\Theta$ be defined by
\begin{equation} \label{eq:ThetaJ}
 \Theta(\lambda) = I_{2p} - B^* \mathbb X^{-1} (\lambda I_n - A_1)^{-1} B J.
\end{equation}
Then there exists $S_0\in \RSC$ such that $S = T_\Theta(S_0)$.
\end{cor}
\textbf{Proof:} Let us denote by $t$ the map \eqref{eq:t}.
For an arbitrary element $f(\lambda) = \eta F(\lambda)$, we shall obtain that
\[ \begin{array}{lllll}
t f(\lambda) & = \eta (\lambda I - A_1)^{-1} B J
\bbmatrix{-S(\lambda)
\\I_p} = \\
& = \eta \operatorname{diag}[\frac{1}{\lambda+\bar w_i}] \bbmatrix{-\xi_1 S(w_1)^* & \xi_1\\ \vdots & \vdots \\ -\xi_n S(w_n)^* & \xi_n}
 \bbmatrix{\sigma_1 S \\ \sigma_1} = \\
& = \sum_i \frac{\eta_i}{\lambda+\bar w_i}\xi_i (\sigma_1 - S(w_i)^*\sigma_1S(\lambda)) = \\
& = \sum_i \xi_i \eta_i K_S(\lambda,w_i)
\end{array} \]
and consequently,
\[\begin{array}{lll}
\langle t f, t f \rangle & = \langle \sum_i \xi_i \eta_iK_S(\lambda,w_i),\sum_j \xi_j \eta_iK_S(\lambda,w_j)\rangle = \\
& = \sum_{ij} \eta_i\xi_i K_S(w_j,w_i) \eta_j^* \xi_j^* = \eta \mathbb X \eta^* = \\
& = \langle f,f \rangle.
\end{array} \]
Thus $t$ is an isometry and 
Theorem \ref{thm:S0Exist} allows to conclude.\qed\mbox{}\\

The following result shows that the assumption $\mathbb X>0$ in the
statement of Corollary \ref{cor:TTheta} can always be achieved for $n=1$.

\begin{lemma}
\label{lemma:Sch1OK}
Given $S\in\RSC$ which is not the function identically equal to $I_p$. Then
there exist a pair $(\xi,w)\in{\mathbb C}^{1\times p}\times {\mathbb C}_+$
such that the corresponding $\mathbb X>0$.
\end{lemma}
\textbf{Proof:} We proceed by contradiction. Assume
that for each $w\in {\mathbb C}_+$
and for each vector $\xi$, it holds that 
\[
\xi\sigma_1\xi^*=\xi S(w)^*\sigma_1 S(w)\xi^*.
\]
or, equivalently,
\[\xi K_S(w,w)\xi^*=0.\]
Then, $\xi K_S(\lambda,w)\equiv 0$, and for each $f\in\mathcal H(S)$
\[ \xi f(w) = < f(\lambda,\xi K_S(\lambda,w)> = 0.
\]
The space $\mathcal H(S)$ is thus trivial, and
its kernel is zero:
\[ \frac{\sigma_1 - S(w)^*\sigma_1 S(\lambda)}{\lambda+\bar w} = 0,
\]
from where we conclude that for each $w,\lambda$
\[ S(\lambda) = \sigma_1^{-1} S^{-*}(w) \sigma_1.
\]
and consequently, $S(\lambda)\equiv I_p$.
\qed
\subsection{Analysis of the tangential Schur algorithm}
For $n=1$, the
matrix function in \eqref{eq:ThetaJ} becomes
\begin{equation}
\label{eq:theta_mmmm}
\begin{array}{lll}
\Theta(\lambda) & = I_{2p} - B^* \mathbb X^{-1} (\lambda I - A)^{-1} B J = \\
& = I_{2p} - \bbmatrix{-\eta^*\\\xi^*} \bbmatrix{-\eta&\xi} \bbmatrix{-\sigma_1&0\\0&\sigma_1} \dfrac{1}{\mathbb X (\lambda + w^*)} = \\
& = \bbmatrix{I_p + \dfrac{\eta^*\eta\sigma_1}{\mathbb X (\lambda + w^*)} & \dfrac{\eta^*\xi\sigma_1}{\mathbb X (\lambda + w^*)} \\
-\dfrac{\xi^*\eta\sigma_1}{\mathbb X (\lambda + w^*)} & I_p
 - \dfrac{\xi^*\xi\sigma_1}{\mathbb X (\lambda + w^*)}} =
\bbmatrix{\Theta_{11}&\Theta_{12}\\\Theta_{21}&\Theta_{22}},
\end{array}
\end{equation}
where using the Lyapunov equation \eqref{eq:LyapJ}
\[ \begin{array}{lll}
A \mathbb X + \mathbb X A^* & = \mathbb X(-w^*-w) = \\
            & = -B J B^* = -\bbmatrix{-\eta & \xi}
\bbmatrix{-\sigma_1&0\\0&\sigma_1} \bbmatrix{-\eta^*\\\xi^*}
\end{array} \]
we find that
\begin{equation}\label{eq:FormX}
\mathbb X = \dfrac{\xi\sigma_1\xi^* - \eta\sigma_1\eta^*}{w+w^*}.
\end{equation}

To the best of our knowledge, Theorem \ref{thm:real} below is new and uses a realization theorem
of M. Liv\v sic \cite{bib:SpectrAnal}.
It will be of much use in the following sections. 
\begin{thm}
\label{thm:real}
Let $S_0\in\SC$ with minimal realization
\begin{equation}
\label{eq:S0!!!}
 S_0(\lambda)
= I_p - B_0^* \mathbb X_0^{-1} (\lambda I - A_0)^{-1} B_0 \sigma_1,
\end{equation}
and let $\Theta$ be given by \eqref{eq:theta_mmmm}. Then
$S=T_\Theta(S_0)\in
\RSC$ and a minimal realization of $S$ is given by
\[ S(\lambda) = I_p - B_S^* \mathbb X_S^{-1}
(\lambda I -A_S)^{-1} B_S \sigma_1,
\]
where
\begin{eqnarray}
\label{eq:BS}
B_S & =& \bbmatrix{ B_0 \\ \eta- \xi } \\
\label{eq:XS}\mathbb X_S & =&
\bbmatrix{\mathbb X_0 & 0 \\ 0 & \mathbb X} =
\bbmatrix{\mathbb X_0 & 0 \\ 0 &
\dfrac{\xi\sigma_1\xi^* - \eta\sigma_1\eta^*}{w+w^*}} \\
\label{eq:AS}A_S & =& \bbmatrix{ A_0 &
\dfrac{B_0 \sigma_1 \xi^*}{\mathbb X} \\ -
\eta\sigma_1 B_0^* \mathbb X_0^{-1}  & - w^* -
\dfrac{\eta\sigma_1(\eta^*-\xi^*)}{\mathbb X }}
\end{eqnarray}
\end{thm}

{\bf Proof:} From the definition
\[ \begin{array}{llll}
S(\lambda) & = T_{\Theta} (S_0(\lambda)) = \\
& = (\Theta_{11} + S_0(\lambda) \Theta_{21})^{-1}(\Theta_{12}
+ S_0(\lambda) \Theta_{22}) = \\
& = (I_p + \dfrac{\eta^*\eta\sigma_1}{\mathbb X (\lambda + w^*)}
 - S_0(\lambda)\dfrac{\xi^*\eta\sigma_1}{\mathbb X (\lambda +
   w^*)})^{-1}
\times \\
&
\hspace{5mm}\times\left(\dfrac{\eta^*\xi\sigma_1}{\mathbb X
(\lambda + w^*)} +
S_0(\lambda)[I_p - \dfrac{\xi^*\xi\sigma_1}{\mathbb X (\lambda + w^*)}]\right).
\end{array} \]
Let us denote here  $S_0(\lambda) \xi^* = \alpha^*$; then the 
preceding expression becomes
\[ \begin{array}{lll}
S(\lambda) & = &
(I_p + \dfrac{\eta^*\eta\sigma_1}{\mathbb X (\lambda + w^*)} - \dfrac{\alpha^*\eta\sigma_1}{\mathbb X (\lambda + w^*)})^{-1} \times \\
& & \times (\dfrac{\eta^*\xi\sigma_1}{\mathbb X (\lambda + w^*)} + S_0(\lambda) - \dfrac{\alpha^*\xi\sigma_1}{\mathbb X (\lambda + w^*)}) = \\
& = &(I_p + \dfrac{(\eta^*-\alpha^*)\eta\sigma_1}{\mathbb X 
(\lambda + w^*)})^{-1}
(S_0(\lambda) + \dfrac{(\eta^*-\alpha^*)\xi\sigma_1}{\mathbb X (\lambda + w^*)}) = \\
\end{array} \]
Using $(I + b^* a)^{-1} = I - b^* a \dfrac{1}{1+ab^*}$ we get
\[  \begin{array}{lll}
 S(\lambda) & = & (I_p -\dfrac{(\eta^*-\alpha^*)\eta\sigma_1}{\mathbb X (\lambda + w^*)}\dfrac{1}{1+\dfrac{\eta\sigma_1(\eta^*-\alpha^*)}{\mathbb X (\lambda + w^*)}})
(S_0(\lambda) + \dfrac{(\eta^*-\alpha^*)\xi\sigma_1}{\mathbb X (\lambda + w^*)})  = \\
& = & S_0(\lambda) - \dfrac{(\eta^*-S_0(\lambda)\xi^*)[\eta\sigma_1 S_0(\lambda) - \xi\sigma_1]}{\mathbb X (\lambda + w^*)+\eta\sigma_1(\eta^*-S_0(\lambda)\xi^*)}.
\end{array} \]
Suppose further that there is a realization of the
form \eqref{eq:S0!!!} for the function $S_0(\lambda)$. Inserting it here we shall obtain
\[  \begin{array}{lll}
S(\lambda) & = & S_0(\lambda) - \dfrac{(\eta^*-S_0(\lambda)\xi^*)[\eta\sigma_1 S_0(\lambda) - \xi\sigma_1]}{\mathbb X (\lambda + w^*)+\eta\sigma_1(\eta^*-S_0(\lambda)\xi^*)} = \\
& = & I_p - B_0^* \mathbb X_0^{-1} (\lambda I - A_0)^{-1} B_0 \sigma_1 - \\
& & - \dfrac{(\eta^*-S_0(\lambda)\xi^*)[\eta\sigma_1 S_0(\lambda) - \xi\sigma_1]}{\mathbb X (\lambda + w^*)+\eta\sigma_1(\eta^*-S_0(\lambda)\xi^*)}.
\end{array} \]
Let us denote
\[  \begin{array}{lll} M & = \mathbb X (\lambda + w^*)+\eta\sigma_1(\eta^*-S_0(\lambda)\xi^*) = \\
& = \mathbb X 
(\lambda + w^*)+\eta\sigma_1(\eta^*-\xi^*) + B_0^* \mathbb X_0^{-1} (\lambda I - A_0)^{-1} B_0 \sigma_1 \xi^*.
\end{array} \]
The preceding formula for $S(\lambda)$ becomes
\[  \begin{array}{lll}
S(\lambda) & = & I_p - B_0^* 
\mathbb X_0^{-1} (\lambda I - A_0)^{-1} B_0 \sigma_1 -
\dfrac{(\eta^*-S_0(\lambda)\xi^*)[\eta\sigma_1 S_0(\lambda) - \xi\sigma_1]}{\mathbb X (\lambda + w^*)+\eta\sigma_1(\eta^*-S_0(\lambda)\xi^*)} = \\
& = & I_p - \dfrac{1}{M} \bbmatrix{B_0^* & 
\eta^*-\xi^* } \bbmatrix{\mathbb X_0^{-1} & 0 \\ 0 & 1} \times\\
& & \times \bbmatrix{ \alpha^{-1} M +
\alpha^{-1} \beta \delta \alpha^{-1} & -\alpha^{-1} \beta \\ - \delta \alpha^{-1} & 1 } \bbmatrix{ B_0 \\ \eta- \xi } \sigma_1,
\end{array} \]
where
\[ \begin{array}{lll}
\alpha & = \lambda I - A_0 \\
\beta  & = -B_0 \sigma_1 \xi^* \\
\delta & = \eta\sigma_1 B_0^* \mathbb X_0^{-1}.
\end{array} \]
By definition of $M$ we have
\[
M = \mathbb X (\lambda + w^*) + \eta\sigma_1(\eta^*-\xi^*) - 
\delta \alpha^{-1} \beta,
\]
and hence 
we obtain from the formula
\[
\bbmatrix{\alpha & \beta \\ \delta & D}^{-1} = \dfrac{1}{D - \delta \alpha^{-1} \beta}\bbmatrix{ \alpha^{-1} (D - \delta \alpha^{-1} \beta) +
\alpha^{-1} \beta \delta \alpha^{-1} & - \alpha^{-1}\beta \\ - \delta \alpha^{-1} & 1}
\]
that the last expression for $S(\lambda)$ is
\[  \begin{array}{lll}
S(\lambda) & =& I_p - \dfrac{1}{M} \bbmatrix{B_0^* & \eta^*-\xi^* } \bbmatrix{\mathbb X_0^{-1} & 0 \\ 0 & 1} \times \\
& & \times \bbmatrix{ \alpha^{-1} M + \alpha^{-1} \beta \delta \alpha^{-1} & -\alpha^{-1} \beta \\ - \delta \alpha^{-1} & 1 }
    \bbmatrix{ B_0 \\ \eta- \xi } \sigma_1 = \\
& = &  I_p - \bbmatrix{B_0^* & \eta^*-\xi^* } \bbmatrix{\mathbb X_0^{-1} & 0 \\ 0 & 1} \times \\
& & \times \bbmatrix{ \alpha & \beta \\ \delta & \mathbb X (\lambda + w^*) + \eta\sigma_1(\eta^*-\xi^*)}^{-1} \bbmatrix{ B_0 \\ \eta- \xi } \sigma_1 = \\
& = & I_p - \bbmatrix{B_0^* & \eta^*-\xi^* } \bbmatrix{\mathbb X_0^{-1} & 0 \\ 0 & \mathbb X^{-1}} \times \\
& & \times  \bbmatrix{ \lambda I -A_0& -\dfrac{B_0 \sigma_1 \xi^*}{\mathbb X} \\ \eta\sigma_1 B_0^* \mathbb X_0^{-1}  & \lambda + w^* + \dfrac{\eta\sigma_1(\eta^*-\xi^*)}{\mathbb X }}^{-1} \bbmatrix{ B_0 \\ \eta- \xi } \sigma_1.
\end{array} \]
Thus we have obtained the realization
\eqref{eq:BS}--\eqref{eq:AS}. Furthermore, the Lyapunov equation
(\ref{eq:LyapJ}) holds since
\[ \begin{array}{lll}
A_S \mathbb X_S + \mathbb X_S A_S^* = \\
= \bbmatrix{ A_0 & \dfrac{B_0 \sigma_1 \xi^*}{\mathbb X} \\ -\eta\sigma_1 B_0^* \mathbb X_0^{-1}  & - w^* - \dfrac{\eta\sigma_1(\eta^*-\xi^*)}{\mathbb X }} \bbmatrix{\mathbb X_0 & 0 \\ 0 & \mathbb X} + \\
\hspace{5mm}+ \bbmatrix{\mathbb X_0 & 0 \\ 0 & \mathbb X}
\bbmatrix{ A_0^* & -\mathbb X_0^{-1}B_0 \sigma_1 \eta^*  \\ \dfrac{\xi \sigma_1 B_0^*}{\mathbb X} &  -w - \dfrac{(\eta^*-\xi^*) \sigma_1\eta^*}{\mathbb X }}\\
= \bbmatrix{A_0 \mathbb X_0 + \mathbb X_0 A_0^* & -B_0 \sigma_1 (\eta^*-\xi^*) \\ -(\eta-\xi) \sigma_1 B_0^* & -\mathbb X(w+w^*) -
\eta\sigma_1(\eta^*-\xi^*) - (\eta^*-\xi^*) \sigma_1\eta^*}  \\
&\vspace{0.5mm}\\
= \bbmatrix{A_0 \mathbb X_0 + \mathbb X_0 A_0^* & -B_0 \sigma_1 (\eta^*-\xi^*) \\ -(\eta-\xi) \sigma_1 B_0^* & \eta\sigma_1\eta^*-\xi\sigma_1\xi^* -
\eta\sigma_1(\eta^*-\xi^*) - (\eta^*-\xi^*) \sigma_1\eta^*} \\
&\vspace{0.3mm}\\
= - B_S \sigma_1 B_S^*.
\end{array} \]
The last equality follows easily from the Lyapunov equation
for the given realization of $S_0$ and the formula (\ref{eq:FormX}).
\mbox{}\qed\mbox{}\\

\section{\label{sec:Schur}The tangential Schur algorithm in the class
$\CI$}
Our strategy to the tangential Schur algorithm in the class $\CI$
relies on the following theorem. This theorem shows in particular that
one cannot use the naive approach of applying the classical tangential
Schur algorithm for each $t_2$ (that is, looking at $t_2$ as a mere
parameter).
\begin{thm}
\label{thm:gamma*}
Let us fix the parameters
  $\sigma_1,\sigma_2$, and $\gamma$, and the interval
$\mathrm I$. Then for every
  $t_2^0\in{\mathrm I}$
there is a one-to-one correspondence between pairs
$(\gamma_*, S)$ such that $S\in\CI$ and $\gamma_*$ continuous in a
neighborhood of $t_2^0$, and functions $Y(\lambda)\in\SC(t_2^0)$.
\end{thm}

\noindent\textbf{Proof:} Let $\phi$ be the map which to a pair
$(\gamma_*, S)$ associates the function $S(\lambda,t_2^0)\in\SC(t_2^0)$.
The converse map 
$\psi: Y(\lambda) \rightarrow Y(\lambda,t_2)$
 was introduced in 
\cite{bib:MyThesis}, \cite[chapter 7]{bib:MelVin1} in a more general
setting, and is defined as follows. Suppose that we have realized the 
transfer function $Y(\lambda)$ in the form
\[ Y(\lambda) = I_p - B_0^* \mathbb X_0^{-1} (\lambda I - A_1)^{-1}
B_0 
\sigma_1(t_2^0).
\]
Then we construct $B(t_2)$ from the differential equation with the
spectral matrix parameter $A_1$:
\begin{equation} \label{eq:DBSol}
\frac{d}{dt_2} [B(t_2) \sigma_1(t_2)] +
A_1 B(t_2) \sigma_2(t_2) + B(t_2) \gamma(t_2) = 0,\quad B(t_2^0) = B_0.
\end{equation}
In fact, the function $B(t_2)$ is given by the formula
\begin{equation}
\label{eq:Bint}
B(t_2)=
\oint(\lambda I - A_1)^{-1} B_0 \sigma_1^{-1}\Phi^{-1}
(\lambda,t_2,t_2^0) d\lambda.
\end{equation}

Next we construct $\mathbb X(t_2)$ on the maximal interval
$I$, where it is invertible (\ref{eq:DX}) via the formula
\begin{equation} \label{eq:DXSol}
\frac{d}{dt_2} \mathbb X(t_2) = B(t_2) \sigma_2(t_2) B(t_2)^*,
\quad \mathbb X(t_2^0) = \mathbb X_0.
\end{equation}
Finally, we define
\[ \begin{array}{llll}
\gamma_*(t_2) & = \gamma(t_2) + \sigma_2(t_2) B(t_2)^* \mathbb X^{-1}(t_2) B(t_2) \sigma_1(t_2) - \\
& \hspace{5mm}-
\sigma_1(t_2) B(t_2)^* \mathbb X^{-1}(t_2) B(t_2) \sigma_2(t_2).
\end{array} \]
Then easy computations show that the function
\begin{equation}
\label{eq:S_min}
 S(\lambda,t_2) = I_p - B(t_2)^*
\mathbb X^{-1}(t_2)(\lambda I - A_1)^{-1} B(t_2) \sigma_1(t_2)
\end{equation}
is in the class $\boldsymbol{\mathcal {CI}}$ corresponding to the parameters
$\sigma_1(t_2)$,$\sigma_2(t_2)$,$\gamma(t_2)$, and $\gamma_*(t_2)$
for $t_2\in\mathrm I$. \textbf{ ?????????????????
Moreover, the Lyapunov equation
\eqref{eq:Lyap} holds for every $t_2\in{\mathrm I}$, and thus
the realization \eqref{eq:S_min} is minimal for every
$t_2\in\mathrm I$}.\\

Note that the composition $\phi\circ\psi=id$, since starting from a
function $Y\in\SC(t_2^0)$, constructing $Y(\lambda,t_2)$ and taking its
value at $t_2^0$, we shall obtain again $Y$ from the initial conditions
of the differential equations \eqref{eq:DBSol} and \eqref{eq:DXSol} defining $B(t_2), \mathbb X(t_2)$.\\

In order to show that $\psi\circ\phi=id$, we 
start from a function $S(\lambda,t_2)\in\CI$ and take its value
$S(\lambda,t_2^0)\in\SC(t_2^0)$. Using the construction above, we shall
obtain a function $Y(\lambda,t_2)$. Note that the two functions
$S(\lambda,t_2)$ and $Y(\lambda,t_2)$
have the same value at $t_2^0$ and maps solutions of the same input
LDE to (possibly different) output LDEs, i.e. :
\[ \begin{array}{ll}
S(\lambda,t_2) =
\Phi_*(\lambda,t_2,t_2^0) S(\lambda,t_2^0) \Phi^{-1}(\lambda,t_2,t_2^0), \\
Y(\lambda,t_2) = \Phi'_*(\lambda,t_2,t_2^0)
S(\lambda,t_2^0) \Phi^{-1}(\lambda,t_2,t_2^0).
\end{array} \]
Then the function $S^{-1}(\lambda,t_2)Y(\lambda,t_2)$ is equal to
$I_p$ at infinity and is entire. By Liouville's theorem it is a
constant function and is equal to $I_p$. Thus
\[
\Phi_*(\lambda,t_2,t_2^0) = \Phi'_*(\lambda,t_2,t_2^0),
\]
from where
we obtain that
\[
\Phi_*^{-1}(\lambda,t_2,t_2^0) \Phi'_*(\lambda,t_2,t_2^0) = I_p.
\]
Differentiating both sides of this last equation we get to
\[ \begin{array}{ll}
0 & = \dfrac{\partial}{\partial t_2} [\Phi_*^{-1}(\lambda,t_2,t_2^0)
\Phi'_*(\lambda,t_2,t_2^0)] = \\
& = \Phi_*^{-1}(\lambda,t_2,t_2^0) \sigma_1^{-1}(t_2) (-\gamma(t_2) +
\gamma'(t_2)) \Phi'_*(\lambda,t_2,t_2^0).
\end{array}  \]
Since the matrices $\Phi_*(\lambda,t_2,t_2^0),
\Phi'_*(\lambda,t_2,t_2^0),\sigma_1(t_2) $ are invertible we obtain
that
$\gamma(t_2) = \gamma'(t_2)$. \qed\mbox{}\\

\textbf{Remark:} Last theorems claims that the correspondence is between
"initial" values $S(\lambda,t_2^0)$ and pairs $(S(\lambda,t_2),\gamma_*(t_2))$.
Notice that it is possible to obtain functions with the same $\gamma_*(t_2)$ with different
initial conditions:
\begin{prop} Suppose that there exists a function $Y\in\SC(t_2^0)$, which commutes
with $\Phi(\lambda,t_2,t_2^0)$ and suppose that a function $S\in\CI$ corresponds to certain
vessel parameters $\sigma_1,\sigma_2,\gamma,\gamma_*$. Then the function 
$S Y$ belongs to the class $\CI$ and corresponds to the same vessel parameters $\sigma_1,\sigma_2,\gamma,\gamma_*$.
\end{prop}
\textbf{Proof:} Using formula \ref{eq:SInttw} we obtain that
\[ S(\lambda, t_2) =
\Phi_*(\lambda,t_2,t_2^0) S(\lambda, t_2^0) \Phi^{-1}(\lambda,t_2,t_2^0) .
\]
Consequently,
\[ S(\lambda, t_2) Y(\lambda) = 
\Phi_*(\lambda,t_2,t_2^0) S(\lambda, t_2^0) Y(\lambda) \Phi^{-1}(\lambda,t_2,t_2^0) 
\]
intertwines solutions of the input \eqref{eq:InCC} and the output \eqref{eq:OutCC}
ODEs with the spectral parameter $\lambda$,
and is identity at infinity, because $S$ and $Y$ and their product are such. Thus by the
definition the function $S Y\in\CI$ and corresponds to the same spectral parameters as
$S$.
\qed

For example (to be studied in subsection \ref{sec:SL}), taking Sturm-Liouville (SL) vessel parameters
\[ \sigma_1 = \bbmatrix{0 & 1 \\ 1 & 0},   \sigma_2 = \bbmatrix{1 & 0 \\ 0 & 0}, \gamma = \bbmatrix{0 & 0 \\ 0 & i}
\]
we shall obtain that
\[ \Phi(\lambda,t_2,t_2^0) = V \bbmatrix{e^{-k (t_2-t_2^0)} & 0 \\ 0 & e^{k (t_2-t_2^0)}} V^{-1},
\]
where
\[ V = \bbmatrix{-\sqrt{\dfrac{i}{\lambda}} & \sqrt{\dfrac{i}{\lambda}} \\ 1 & 1},
\quad k = \sqrt{i\lambda}.
\]
Taking $Y$, which commutes with $\sigma_1^{-1}(\sigma_2 \lambda+ \gamma) = \bbmatrix{0 & i\\ \lambda & 0}$, i.e.
of the form
\begin{equation} \label{eq:InCommPhi}
 Y(\lambda) = I_2 - \bbmatrix{a(\lambda)& \dfrac{i c(\lambda)}{\lambda} \\ c(\lambda) & a(\lambda)}
\end{equation}
for functions $a(\lambda), c(\lambda)$, which are zero at infinity,
 we shall obtain that for any $S\in\CI$, the function $SY\in\CI$ and corresponds to the same
vessel parameters. But in this case (which can be also generalized) it turns to be a necessary condition
as the following theorem states:
\begin{lemma} Given SL vessel parameters $\sigma_1, \sigma_2, \gamma,\gamma_*(x)$, there exists a unique
initial value $S(0,\lambda)$ up to a scalar $t_2$-independent symmetric and identity at
infinity function.
\end{lemma}
\textbf{Proof:} For this it is enough to prove that if a function $Y(\lambda,t_2)\in\CI$ 
intertwines solutions of the 
input LDE \eqref{eq:InCC} with itself, then $Y(\lambda,t_2)=Y(\lambda,t_2^0)$, 
i.e is a constant function, which moreover commutes with $\Phi_I(\lambda,t_2,t_2^0)$.
Indeed, if we are given two functions $S_1(\lambda,t_2)$, $S_2(\lambda,t_2)$, then the function 
$S_1^{-1}(\lambda,t_2)S_2(\lambda,t_2)$ will intertwine solutions of the input LDE with itself and as a result
is constant $Y(\lambda)\in\SC$, which is moreover commutes with $\Phi_I(\lambda,t_2)$. Following the paragraph
preceding the theorem, $Y(\lambda)=I$ and this means that $S_1(\lambda,t_1^0)=S_2(\lambda,t_2^0)$.

Let us show first that if $S(\lambda,t_2)$ intertwines solutions of the input LDE \eqref{eq:InCC} 
with itself and is identity at $\lambda=\infty$, then it is a
constant matrix, which commutes with $\Phi(\lambda,t_2,t_2^0)$. Performing simple calculations, we can find that
\[ \sigma_1^{-1} (\sigma_2 \lambda + \gamma) = \bbmatrix{0 & i \\ \lambda & 0},
\quad \Phi(\lambda,t_2,t_2^0) = V \bbmatrix{e^{-k (t_2-t_2^0)} & 0 \\ 0 & e^{k (t_2-t_2^0)}} V^{-1},
\]
where
\[ V = \bbmatrix{-\sqrt{\dfrac{i}{\lambda}} & \sqrt{\dfrac{i}{\lambda}} \\ 1 & 1},
\quad k = \sqrt{i\lambda}.
\]
Consequently, for the expression $S(\lambda,t_2) = \Phi(\lambda,t_2,t_2^0) S_0(\lambda) \Phi^{-1}(\lambda,t_2,t_2^0)$
to  be identity for $\lambda=\infty$, it is necessary to "cancel" the essential singularity arising from two
entire functions $\Phi(\lambda,t_2,t_2^0)$ and $\Phi^{-1}(\lambda,t_2,t_2^0)$ (or more precisely, let them cancel
each other). 

Using the formula for $\Phi(\lambda,t_2,t_2^0)$
\begin{eqnarray*}
\Phi(\lambda,t_2,t_2^0) S_0(\lambda,t_2^0) \Phi^{-1}(\lambda,t_2,t_2^0) = \\
V \bbmatrix{e^{-k (t_2-t_2^0)} & 0 \\ 0 & e^{k (t_2-t_2^0)}} V^{-1} S(\lambda,t_2^0)
V^{-1} \bbmatrix{e^{k (t_2-t_2^0)} & 0 \\ 0 & e^{-k (t_2-t_2^0)}} V^{-1} = I
\end{eqnarray*}
and considering coefficients of the exponents, we may conclude that in order to cancel the exponents 
it is necessary to demand
\[ V^{-1} S(\lambda,t_2^0) V = \bbmatrix{b(\lambda) & 0 \\ 0 & d(\lambda)}.
\]
for some analytic in $\sqrt{\lambda}$ at infinity functions $b(\lambda), d(\lambda)$ which are $1$ there. From here it follows that
\[ S(\lambda,t_2^0) = \bbmatrix{-\sqrt{\dfrac{i}{\lambda}}(b(\lambda) + d(\lambda)) &  \dfrac{i}{\lambda}(d(\lambda) - b(\lambda)) \\ d(\lambda) - b(\lambda) & -\sqrt{\dfrac{i}{\lambda}}(b(\lambda) + d(\lambda))}.
\]
denoting here $1-a(\lambda) = -\sqrt{\dfrac{i}{\lambda}}(b(\lambda) + d(\lambda))$ and
$-c(\lambda) = d(\lambda) - b(\lambda)$, we shall obtain that $S_0(\lambda)$ is of the form \eqref{eq:InCommPhi},
i.e. commutes with the fundamental matrix $\Phi(\lambda,t_2,t_2^0)$.

Finally, let us consider a function of this form and will require the symmetry condition for it:
\[ S^*(\lambda,t_2^0) \sigma_1 S(\lambda,t_2^0) = \sigma_1.
\]
Plugging here the expression \eqref{eq:InCommPhi} for $S(\lambda,t_2^0)$
\[ S(\lambda,t_2^0) = I_2 - \bbmatrix{a(\lambda)& \dfrac{i c(\lambda)}{\lambda} \\ c(\lambda) & a(\lambda)}
\]
we shall obtain that on the one imaginary axis, where $\lambda=-\bar \lambda$,
the following system of equations must hold
\[ \left\{ \begin{array}{lll}
(1-a^*(-\bar\lambda)) c(\lambda) + c^*(-\bar\lambda) (1-a(\lambda)) = 0, \\
(1-a^*(-\bar\lambda)) (1-a(\lambda)) + \dfrac{1}{\lambda} c^*(-\bar\lambda)c(\lambda) = 1, \\
(1-a^*(-\bar\lambda)) (1-a(\lambda)) + \dfrac{1}{\bar\lambda} c^*(-\bar\lambda)c(\lambda) = 1, \\
\dfrac{i c^*(-\bar\lambda) (1-a(\lambda))}{\lambda} + \dfrac{i c(\lambda) (1-a^*(-\bar\lambda))}{\lambda} = 0.
\end{array} \right.
\]
Subtracting the second and the third equation we obtain that (remember that $\lambda=-\bar\lambda$)
\[ c^*(-\bar\lambda)c(\lambda) (\dfrac{1}{\lambda} - \dfrac{1}{\bar\lambda}) =
c^*(-\bar\lambda)c(\lambda) \dfrac{2}{\lambda} = 0.
\]
from where it follows that on the imaginary axis at the points of analyticity of $c(\lambda)$ it holds that
$c(\lambda)=0$. Since it is analytic there, it will be globally zero.
Notice that the last system of equations for the trivial $c(\lambda)$ gives 
\[ (1-a^*(-\bar\lambda)) (1-a(\lambda)) = 1,
\]
Which means that $1-a(\lambda)$ is of norm one on the imaginary axis. \qed

We now focus on the rational case. We first note
that starting from different realizations at
$t_2^0$ 
we shall obtain the same $\gamma_*(t_2)$. More precisely we have the
following theorem:
\begin{thm}
\label{thm:S2min}
Suppose that there are two minimal realizations of the
function $S(\lambda,t_2^0)\in\RSC$:
\[ S(\lambda,t_2^0) = I_p - B_\ell^* \mathbb X_\ell^{-1} (\lambda I -
A_1)^{-1} B_\ell
\sigma_1(t_2^0),\quad \ell=1,2,
\]
with associated similarity matrix $V$.
Then the functions
\[
S_\ell(\lambda, t_2)=I_p-B_\ell(t_2)^*\mathbb X_\ell^{-1}(t_2)
(\lambda I-A_\ell)^{-1}B_\ell(t_2)\sigma_1(t_2),\quad \ell=1,2,
\]
obtained from these realizations via the construction in Theorem
\ref{thm:gamma*} are the same, and this holds if and only if:
\begin{eqnarray}
B_2(t_2) = V B_1(t_2),
\end{eqnarray}
where $t_2$ varies in a neighborhood of $t_2^0$.
Moreover, in this case the following formula holds
\[ 
\mathbb X_2(t_2) = V \mathbb X_1(t_2) V^*.
\]
\end{thm}

\textbf{Proof:}
Equality of the two realizations means that there exists an invertible
matrix $V$ such that
\[ A_2 = V A_1 V^{-1}, ~~~~~ B_2\sigma_1(t_2^0) = V B_1
\sigma_1(t_2^0)
, ~~~~~ B_1^*\mathbb X_1^{-1}V = B_2^*\mathbb X_2^{-1},
\]
from which follows (see \cite[Lemma 2.1 p. 184]{bib:ag} for instance) that
\[
\mathbb X_2=V\mathbb X_1 V^*.
\]
By the construction described in Theorem \ref{thm:gamma*},
the function $B_1(t_2)$ will satisfy \eqref{eq:DBSol}
\[\frac{d}{dt_2} [B_1(t_2) \sigma_1(t_2)] + A_1 B_1(t_2) \sigma_2(t_2)
+ B_1(t_2) \gamma(t_2) = 0,\quad B_1(t_2^0) = B_1, \]
and the function $B_2(t_2)$ will satisfy the same equation with $A_2$
instead of $A_1$:
\[ \frac{d}{dt_2} [B_2(t_2)
\sigma_1(t_2)] + A_2 B_2(t_2) \sigma_2(t_2) + B_2(t_2) \gamma(t_2) =
0, 
\quad B_2(t_2^0) = B_2. \]
Using the equalities $A_2 = V A_1 V^{-1}$, $B_2 = V B_1$ we obtain that
the function $V^{-1} B_2(t_2)$ satisfies
\[ \begin{array}{lll}
\frac{d}{dt_2} [V^{-1} B_2(t_2) \sigma_1(t_2)] + A_1 V^{-1}B_2(t_2) \sigma_2(t_2) + V^{-1} B_2(t_2) \gamma(t_2) = 0, \\
V^{-1} B_2(t_2^0) = B_1,
\end{array} \]
which is the same differential equation as for $B_1(t_2)$.
Thus $B_2(t_2) = V B_1(t_2)$.
Similarly, considering the differential equations
\[ \frac{d}{dt_2} \mathbb X_i(t_2) =
B_i(t_2) \sigma_1 B_i(t_2)^*, \quad \mathbb X_i(t_2^0) = \mathbb X_i, \quad i=1,2,
\]
we obtain that
\[ \mathbb X_2(t_2) = V \mathbb X_1(t_2) V^*.
\]
Consequently,
\[ 
\begin{array}{lll} B_2^*(t_2) 
\mathbb X_2^{-1}(t_2) B(t_2) & = B_1^*(t_2) V^* V^{-1*} 
\mathbb X_1^{-1}(t_2) V^{-1} V B_1(t_2)  \\
& = B_1^*(t_2) \mathbb X_1^{-1}(t_2) B_1(t_2),
\end{array}  
\]
from where we conclude that the same function $\gamma_*$ is associated
to $S_1(\lambda,t_2)$ and $S_2(\lambda,t_2)$.
Since these functions coincide 
for $t_2=t_2^0$ and map the same input ODE, we obtain that  they are
equal in a neighborhood of $t_2^0$.
\qed\mbox{}\\

The following notion has been introduced and studied in
 \cite{bib:SLVessels}.

\begin{defn}
\label{def:tau}
Let $S\in\RCI$ with a realization \eqref{eq:S_min}.
The function
\[
\tau(t_2) = \det \mathbb X(t_2).
\]
is called the $\tau$-function associated to $S$.
\end{defn}
It follows from Theorem \ref{thm:S2min} that the
$\tau$ function is well defined up to a multiplicative strictly
positive constant. Indeed using the notation of the
theorem,
\[ \det \mathbb X_2(t_2) = \det [V \mathbb X_1(t_2) V^*] = \det \mathbb X_1(t_2) \det (VV^*).
\]

We now introduce the counterpart of the tangential Schur algorithm
in the class $\RCI$. Let $S\in\RCI$ and fix $t_2^0\in\mathrm I$. The
function
$S(\lambda, t_2^0)\in\RSC$.
Consider now a space $\mathcal M$ with $\mathbb X>0$ and corresponding function
$\Theta$ as
in Corollary \ref{cor:TTheta}. This is always possible in view of
Lemma \ref{lemma:Sch1OK}. It follows from Theorem \ref{thm:S0Exist} that
there exists $S_0\in\RSC$ such that
\begin{equation}
\label{eq:STThetaS0}
S(\lambda,t_2^0) = T_{\Theta(\lambda)}(S_0(\lambda)).
\end{equation}
Applying Theorem \ref{thm:gamma*} to $S_0$, we obtain a
uniquely defined function $S_0(\lambda,t_2)\in\RCI$, such that
at $t_2^0$ the relation (\ref{eq:STThetaS0}) holds.

\begin{defn}\label{def:SchurAlg}
The map $T_{\Theta,t_2^0}$
\[ S(\lambda,t_2) \mapsto S_0(\lambda,t_2)
\]
is the time-varying counterpart of the linear fractional
transformation \eqref{eq:STThetaS0}. We will call it a generalized
linear fractional transformation.
\end{defn}

If $S(\lambda,t_2)\in\CI$
corresponds to the
vessel parameters $\sigma_1$, $\sigma_2$, $\gamma$, $\gamma_*$, then 
$S_0(\lambda,t_2^0) \in \CI$
corresponds to the vessel parameters $\sigma_1$, $\sigma_2$,
$\gamma$, $\gamma_{0,*}$ for a uniquely defined
function $\gamma_{0,*}(t_2)$. Moreover,
for $t_2=t_2^0$ we have the usual linear fractional transformation 
\eqref{eq:STThetaS0}.\\

As a consequence of Theorem \ref{thm:gamma*} the following lemma holds.
\begin{lemma}
For a given $J$-inner function $\Theta$ and a given point $t_2^0$, 
the map $T_{\Theta,t_2^0}$ is one-to-one from $\RSC$ into $\RCI$.
\end{lemma}
\textbf{Proof:}
Notice that for a given $t_2^0$ the map $T_\Theta$ is injective.
Furthermore, using Theorem \ref{thm:gamma*}, every $S\in\RCI$ is
uniquely defined by the function $S(\lambda, t_2^0)\in\RSC$.
The result follows. \qed\mbox{}\\

Suppose now that we start from $S_0(\lambda)\equiv I_p$
and apply $n$ linear fractional transformations for a
fixed $t_2^0$, using the data $<w_i,\xi_i,\eta_i>$ ($i=1,\ldots n$) to
construct the corresponding $J$-inner functions.
We obtain a function
\[ S_n(\lambda) = T_\Theta(I) = I - B_n^* X_n^{-1} (\lambda I - A_n)^{-1} B_n \sigma_1\in\RSC.
\]
Using iteratively formulas (\ref{eq:BS}), (\ref{eq:XS}), (\ref{eq:AS})
we obtain that
\[ \begin{array}{llll}
B_n = \bbmatrix{\eta_1-\xi_1\\\vdots\\\eta_n-\xi_n},  \\
\mathbb X_n = \operatorname{diag}[\widetilde{\mathbb X}_1,
\ldots,\widetilde{\mathbb X}_n],
\end{array} \]
and 
{\tiny
\begin{equation}
\label{eq:an}
A_n
=
\bbmatrix{-w_1^* - \dfrac{\eta_1\sigma_1(\eta_1^*-\xi_1^*)}{\widetilde{\mathbb X}_1} & \dfrac{(\eta_1-\xi_1)\sigma_1\xi_2^*}{\widetilde{\mathbb X}_2} & \ldots & \dfrac{(\eta_1-\xi_1)\sigma_1\xi_n^*}{\widetilde{\mathbb X}_n} \\
-\dfrac{\eta_2\sigma_1(\eta_1^*-\xi_1^*)}{\widetilde{\mathbb X}_1} &
-w_2^* - \dfrac{\eta_2\sigma_1(\eta_2^*-\xi_2^*)}{\widetilde{\mathbb X}_2} & \ldots & \dfrac{(\eta_2-\xi_2)\sigma_1\xi_n^*}{\widetilde{\mathbb X}_n} \\
\vdots & \vdots & \ddots & \vdots \\
-\dfrac{\eta_n\sigma_1(\eta_1^*-\xi_1^*)}{\widetilde{\mathbb X}_1} &
-\dfrac{\eta_n\sigma_1(\eta_2^*-\xi_2^*)}{\widetilde{\mathbb X}_2} & \ldots &
-w_n^* - \dfrac{\eta_n\sigma_1(\eta_n^*-\xi_n^*)}{\widetilde{\mathbb X}_n}
},
\end{equation}
}
where $\widetilde{\mathbb X}_i$ is defined by (\ref{eq:FormX}):
\[ \widetilde{\mathbb X}_i = \dfrac{\xi_i\sigma_1\xi_i^* -
\eta_i\sigma_1\eta_i^*}{w_i+w_i^*},\quad i=1,\ldots n.
\]
Assume now that, starting from the identity matrix we apply this procedure
for two different sets of data (with the same $n$) at two different
points $t_2^1$ and $t_2^2$. The following theorem answers the question
as when we obtain the same function, that is, when do we have:
\[
T_{\Theta_1,t_2^1} (I_p)= T_{\Theta_2,t_2^2} (I_p).
\]
\begin{thm}
\label{thm:EqualS} Suppose that there are given two sets of $n$ triples \linebreak
$<w^1_i,\xi^1_i,\eta^1_i>$ and $<w^2_i,\xi^2_i,\eta^2_i>$, with corresponding
$\Theta_\ell, \ell=1,2$. Then necessary
and sufficient conditions for equality of the two functions
\[
S_\ell(\lambda,t_2) = T_{\Theta_\ell,t_2^\ell} (I)
= I_p - B^\ell_n \mathbb (\mathbb X_n^\ell)^{-1}
(\lambda I - A_n^\ell)^{-1} B_n^\ell \sigma_1,\quad \ell=1,2
\]
are:
\begin{enumerate}
\item The corresponding matrices $A_n^1$ and $A_n^2$ defined by
\eqref{eq:an} are similar, 
i.e. there exists an invertible
matrix $V$ such that $A_n^1 = V A_n^2 V^{-1}$,
\item $\oint(\lambda I - A_n^1)^{-1} B_n^1 \sigma_1^{-1}\Phi^{-1}
(\lambda,t_2^2,t_2^1) d\lambda = V B_n^2$.
\end{enumerate}
\end{thm}
\textbf{Proof:} From Theorem \ref{thm:gamma*}, a necessary and
sufficient condition for the functions to be equal is that
\[
S_1(\lambda,t_2^2)=S_2(\lambda, t_2^2).
\]
From \ref{thm:S2min} this holds if and only if
\[ A_n^2 = V A_n^1 V^{-1},
~~~~~B_n^2(t_2^2) = V B_n^1, ~~~~~\mathbb X_n^2(t_2^2) =
V \mathbb X_n^1 V^*
\]
for a uniquely defined invertible matrix $V$.
The result follows using formula \eqref{eq:Bint}. \qed

A more general construction in this setting is obtained if one
supposes that at each step different
values of $t_2$ are chosen. In this case the construction of the
function $S_n(\lambda,t_2)$ is more complicated, and can be
computed recursively.
The formulas are very involved
in this case and we can see no real advantage to develop them
at this point.

\section{\label{sec:MarkovPart}Markov
moments and partial realization problem in the class
$\CI$}
Let $S\in\SC$, and consider the Laurent expansion at infinity
\[ S(\lambda) = I_p - B^* \mathbb X^{-1} (\lambda I - A)^{-1} B
\sigma_1
= I_p - \sum\limits_{i=0}^\infty
\frac{1}{\lambda^{i+1}} B^* \mathbb X^{-1} A^i B \sigma_1
\]
in terms of a given minimal realization. The matrices
$H_i = B^* \mathbb X^{-1} A^i B\sigma_1$ are called the {\sl Markov
moments} of $S$. The
Partial Realization Problem (or moment problem at infinity) 
in $\SC$ is defined as follows:
Given the first $n+1$ Markov moments $H_0, \ldots,H_n$, find all
functions (if any) $S\in\SC$ with these first $n+1$ moments. See
\cite{bib:gkl-scl} for a general study of the partial realization
problem. Similarly, one can define the Markov moments for an element
$S\in\CI$.
We now give necessary conditions which the moments of a function
$S\in\CI$ have to satisfy.

It is important to notice the following: fixing $t_2=t_2^0$ and
solving the corresponding classical moment problem will not lead to a
solution
of the problem in the class $\CI$ because we cannot obtain the
function
$\gamma$ (and hence $\gamma_*$) from this solution.

\subsection{Fundamental properties of $S(\lambda,t_2)$}
Before we consider moments of the function $S(\lambda,t_2)$ we present some of its fundamental
properties which shed more light on the deriving of the moment equations. 
\begin{thm} For fixed $\sigma_1, \sigma_2, \gamma$, a necessary condition on $\gamma_*$
so that there exists a vessel with the vessel parameters 
$\sigma_1, \sigma_2, \gamma, \gamma_*$ is
\[ \det \Phi_*(\lambda,t_2,t_2^0) = \det \Phi(\lambda,t_2,t_2^0)
\]
\end{thm}
\textbf{Proof:}
Let $S\in\CI$ be a function corresponding to the parameters $\sigma_1, \sigma_2, \gamma, \gamma_*$. Then
using a realization \eqref{eq:S_min}
\[ S(\lambda,t_2) = I_p - B(t_2)^*
\mathbb X^{-1}(t_2)(\lambda I - A_1)^{-1} B(t_2) \sigma_1(t_2)
\]
and Lyapunov equation \eqref{eq:Lyap} we shall obtain that
\[ \begin{array}{llll}
\det S(\lambda,t_2) & = \det \big( I_p - B(t_2)^*
\mathbb X^{-1}(t_2)(\lambda I - A_1)^{-1} B(t_2) \sigma_1(t_2) \big) \\
& =  \det \big( I - B(t_2) \sigma_1(t_2) B(t_2)^*
\mathbb X^{-1}(t_2)(\lambda I - A_1)^{-1} \big) \\
& = \det \big( I + (A_1\mathbb X(t_2) + \mathbb X(t_2) A_1^*)
\mathbb X^{-1}(t_2)(\lambda I - A_1)^{-1} \big) \\
& = \det \big( I + A_1(\lambda I - A_1)^{-1} +\mathbb X(t_2) A_1^*
\mathbb X^{-1}(t_2)(\lambda I - A_1)^{-1} \big) \\
& = \det \left( (\lambda I + \mathbb X(t_2) A_1^*
 \mathbb X^{-1}(t_2)) (\lambda I - A_1)^{-1} \right) \\
& = \det (\lambda I + A_1^*) \det (\lambda I - A_1)^{-1} \\
& = \det S(\lambda,t_2^0).
\end{array} \]
Consequently, taking determinant of the formula \eqref{eq:SInttw}
\[ S(\lambda,t_2) \Phi(\lambda,t_2,t_2^0) = \Phi_*(\lambda,t_2,t_2^0) S(\lambda,t_2^0)
\]
we shall obtain that $\det \Phi_*(\lambda,t_2,t_2^0) = \det \Phi(\lambda,t_2,t_2^0)$ for all points
$\lambda$, where $\det S(\lambda,t_2^0)$ exists and is different from zero. Since it happens for all points outside
the spectrum of $A_1$ and the functions $\det \Phi_*(\lambda,t_2,t_2^0)$, $\det \Phi(\lambda,t_2,t_2^0)$ are entire
they are equal for all $\lambda$.
\qed

Next theorem put some light on the connection between equations, which determine the transfer function 
$S(\lambda,t_2)$ and its moments. This theorem is similar to the property of a solution of a Riccati
equation \cite[theorem 2.1]{bib:Zelikin}
\begin{thm}\label{thm:S=S-1*}Suppose that $S(\lambda,t_2)$ is a continuous function of $t_2$ for each $\lambda$,
meromorphic in $\lambda$ for each $t_2$ and satisfies $S(\infty,t_2) = I$. Suppose also that
$S(\lambda,t_2)$ is an intertwining function of LDEs \eqref{eq:InCC} and \eqref{eq:OutCC}.
Then if the symmetry condition \eqref{eq:Symmetry}
\[
S(\lambda,t_2) = \sigma_1^{-1}(t_2) S^{-1*}(-\bar\lambda,t_2) \sigma_1(t_2)
\]
holds for $t_2^0$, then it holds for all $t_2$.
\end{thm}
\textbf{Proof:} Since $S(\lambda,t_2)$ intertwines solutions of \eqref{eq:InCC} and \eqref{eq:OutCC}, then
it satisfies the differential equation \eqref{eq:DS}
\begin{eqnarray*} \frac{\partial}{\partial t_2} S(\lambda,t_2)
    = \sigma_{1}^{-1}(t_2) (\sigma_{2}(t_2)
\lambda + \gamma_*(t_2)) S(\lambda,t_2)-\\
-   S(\lambda,t_2)\sigma_1^{-1}(t_2) (\sigma_2(t_2) \lambda + \gamma(t_2)).
\end{eqnarray*}
Consequently, using properties of $\gamma_*, \gamma$ appearing in definition \ref{def:VesPar}
we obtain that the function 
$\sigma_1^{-1}(t_2) S^{-1*}(-\bar\lambda,t_2) \sigma_1(t_2)$ satisfies the same differential equation.
If these two functions are equal at $t_2^0$, from the uniqueness of solution for a differential
equation with continuous coefficients, they are also equal for all $t_2$.
\qed
\subsection{Restrictions on Markov moments for functions in
$\CI$}
We study the Markov moments of a function $S\in\CI$, 
which maps solutions of the
input ODE (\ref{eq:InCC}) to the output ODE (\ref{eq:OutCC}) using a
minimal realization \eqref{eq:S_min} of $S$.

If $A_1$ is a constant matrix, at some stage the elements
$I,A_1,\ldots,A_1^n$ will be linearly dependent and we obtain:
\begin{lemma}
Given $S\in\RCI$, with Markov moments $H_i(t_2)$, $i=0,1,\ldots$. Then
exists $N$ and constants $\mu_j$, $j=1,\ldots,N$ such that
\begin{equation}
\label{eq:HLin}
\sum\limits_{j=0}^{N+1} \mu_j H_{n-j}(t_2) = 0, \quad n\geq N+1.
\end{equation}
\end{lemma}
Next we present formulas, which hold for the moments of a function in $\CI$. They can be easily derived from
theorem \ref{thm:S=S-1*}.
We notice that the first moment $H_0(t_2)$ satisfies
the linkage condition (\ref{eq:H0}):
\[
\gamma_*(t_2) - \gamma(t_2) = \sigma_2(t_2) H_0(t_2) - \sigma_1(t_2) H_0(t_2) \sigma_1^{-1}(t_2) \sigma_2(t_2).
\]
Let us denote 
\[
H_0(t_2) = C(t_2) B(t_2) \sigma_1(t_2) = (B(t_2))^*
X^{-1}(t_2)B(t_2)\sigma_1(t_2). 
\]
The functions $B(t_2), C(t_2)$ satisfy the
following differential equations 
\begin{eqnarray}
\dfrac{d}{dt_2} [B(t_2) \sigma_1(t_2)] + A_1 B(t_2) \sigma_2(t_2) + B(t_2) \gamma(t_2) = 0 \\
\sigma_1(t_2) \frac{d}{dt_2} C(t_2) - 
\sigma_2(t_2) C(t_2) A_1 - \gamma_*(t_2) C(t_2) = 0,
\end{eqnarray}
see \cite{bib:MyThesis,bib:MelVin1}.
Thus, differentiating $H_0(t_2)$, we obtain
\[ \begin{array}{lll}
 \frac{d}{dt_2} H_0(t_2) & = 
\frac{d}{dt_2} [C(t_2) B(t_2)\sigma_1(t_2)]  \\
  & = \sigma_1^{-1}(t_2) \sigma_2(t_2) C(t_2) A_1 B(t_2) \sigma_1(t_2) - C(t_2) A_1 B(t_2) \sigma_2(t_2)+ \\
  & \hspace{5mm}+ \sigma_1^{-1}(t_2) 
\gamma_*(t_2) H_0(t_2) - H_0(t_2) \sigma_1^{-1}(t_2)\gamma(t_2).
\end{array} \]
In other words the second moment $H_1(t_2) = C(t_2) A_1
B(t_2)\sigma_1(t_2)$ 
satisfies the following differential equation
\[ \begin{array}{lll}
\sigma_1^{-1}(t_2) \sigma_2(t_2) H_1(t_2) - H_1(t_2) 
\sigma_1^{-1}(t_2)\sigma_2(t_2) = \\
= \frac{d}{dt_2} H_0(t_2) - \sigma_1^{-1} \gamma_* H_0(t_2) + H_0(t_2) \sigma_1^{-1}(t_2)\gamma(t_2).
\end{array} \]
In the same manner the moments $H_i(t_2)$ and $H_{i+1}(t_2)$ are connected by the following differential equation
\begin{equation} \label{eq:HiHi+1}
\begin{array}{lll}
\sigma_1^{-1}(t_2) \sigma_2(t_2) H_{i+1} - H_{i+1}(t_2) \sigma_1^{-1}(t_2) \sigma_2(t_2) = \\
= \frac{d}{dt_2} H_i(t_2) - \sigma_1^{-1}(t_2) \gamma_*(t_2) H_i(t_2) + H_i(t_2) \sigma_1^{-1}(t_2)\gamma(t_2).
\end{array} \end{equation}
Notice also that 
\[ \begin{array}{lllll}
H_i \sigma_1^{-1} H_0^* \sigma_1 &  = B^* X^{-1} A_1^i B \sigma_1 B^* X^{-1} B  \sigma_1 = \\
& = B^* X^{-1} A_1^i [-A_1 X - X A_1^*] X^{-1} B  \sigma_1 = \\
& = - H_{i+1} - B^* X^{-1} A_1^i X A_1^* X^{-1} B  \sigma_1 = \\
& = - H_{i+1} - B^* X^{-1} A_1^{i-1} [-B\sigma_1 B^* - XA_1^*] A_1^* X^{-1} B  \sigma_1= \\
& = - H_{i+1} + B^* X^{-1} A_1^{i-1} B\sigma_1 B^* A_1^* X^{-1} B  \sigma_1 + \\
& +  B^* X^{-1} A_1^{i-1} XA_1^* A_1^* X^{-1} B \sigma_1  = \\
& = - H_{i+1} + H_{i-1} \sigma_1^{-1} H^*_1 \sigma_1 + B^* X^{-1} A_1^{i-1} X A_1^{2*} X^{-1} B \sigma_1 = \\
& = \cdots = \\
& = - H_{i+1} + H_{i-1} \sigma_1^{-1} H_1^* \sigma_1 - H_{i-2} \sigma_1^{-1} H_2 \sigma_1 + \\
& + \cdots + (-1)^i \sigma_1^{-1} H_{i+1}^* \sigma_1
\end{array} \]
Consequently, the following formula  holds:
\begin{equation}
\label{eq:Moments}
H_{i+1}\sigma_1^{-1} + (-1)^i \sigma_1^{-1} H_{i+1}^*  =
\sum_{j=0}^{i} (-1)^{j+1} H_{i-j} \sigma_1^{-1} H_j^*.
\end{equation}
\textbf{Remark:} One can obtain the same equations \eqref{eq:HiHi+1} and \eqref{eq:Moments} by taking the expansion
of $S(\lambda,t_2)$ into Taylor series around $\lambda=\infty$ and equating the coefficients of each
$\dfrac{1}{\lambda^i}$. More precisely, the algebraic condition \eqref{eq:HiHi+1} is a consequence of the symmetry
condition \eqref{eq:Symmetry} and the condition \eqref{eq:Moments} is a result of the differential equation
\eqref{eq:DS}.

Finally, we show how the third condition in Proposition \ref{prop:PropS} is reflected in the moments $H_i(t_2)$. 
This condition means that the function $S(\lambda,t_2)$ is $\sigma_1(t_2)$ contractive, and,
for example, the first moment $H_0(t_2)$ satisfies $H_0(t_2) \sigma_1^{-1}(t_2) > 0$.
Using the minimality property, we have
\[ 
\operatorname{span} A_1^n B \mathcal E = \mathcal H.
\]
Thus, in order to have $\mathbb X>\epsilon I, \epsilon>0$ (so that  $\mathbb X^{-1}$ exists) it is enough
to demand
\[ \langle \mathbb X^{-1} \sum_{k=0}^n A^k B e_k, \sum_{k=0}^n A^k B e_k \rangle > \epsilon, 
\text{ for any choice of $e_k \in\mathcal E$}.
\]
Thus we obtain that the following Pick matrix arises
\[ \mathbb P = \bbmatrix{B^*\\B^* A_1^*\\\vdots\\B^* (A_1^*)^n} \mathbb X^{-1} 
\bbmatrix{B&A_1 B&\ldots & A_1^nB} > \epsilon I
\]
which is equal to
\[ \mathbb P_n = [p_{ij}]_{i,j}^{n} = [B^* (A_1^*)^i\mathbb X^{-1} A_1^j B]>\epsilon I, \quad n=0,1,2,\ldots
\]
Notice that the elements of the matrix are indexed from zero to $n$ for convenience.
Using Lyapunov equation \eqref{eq:Lyap} iteratively, we shall obtain that
\[ \begin{array}{llll}
B^* (A_1^*)^m \mathbb X^{-1} A_1^n B & =
B^* (A_1^*)^{m-1} (-\mathbb X^{-1}B\sigma_1 B^* \mathbb X^{-1} - \mathbb X^{-1} A_1 ) A_1^n B \\
& = - \sigma_1^{-1} H^*_{m-1} \sigma_1 H_n \sigma_1^{-1} - B^* (A_1^*)^{m-1} \mathbb X^{-1} A_1^{n+1} B \\
& = \cdots = \\
& = \sum\limits_{k=0}^m (-1)^{k+1} \sigma_1^{-1} H^*_{m-1-k} \sigma_1 H_{n+k} \sigma_1^{-1}, \quad m>0, \\
B^* \mathbb X^{-1} A_1^n B & = H_n.
\end{array} \]
So the Pick matrix is of the following form
\begin{equation} \label{eq:Pick}
\mathbb P = [p_{ij}]_{i,j=0}^n, 
\left\{ \begin{array}{lll}
p_{ij} & = \sum\limits_{k=0}^{i} (-1)^{k+1} \sigma_1^{-1}  H^*_{i-1-k} \sigma_1 H_{j+k} \sigma_1^{-1}, & i>1 \\
p_{0,j} & = H_j.
\end{array} \right.
\end{equation}
Next theorem shows that moments of a transfer function in $\SC$ satisfy recursive equations, involving 
algebraic and differential equations:
The positivity of the matrix $\mathbb P$, in that case does not necessarily hold.

\begin{lemma} \label{lemma:Rok}
Denote the real and the imaginary parts of $H_{2n+1} \sigma_1^{-1} = R_{2n+1} + i M_{2n+1}$. Then if the moments
$H_0,\ldots,H_{2n}$ satisfy the differential equations \eqref{eq:HiHi+1} and \eqref{eq:Moments}, then the real 
part $R_{2n+1}$ satisfies the imaginary part of differential equation \eqref{eq:Moments}:
\begin{equation} \label{eq:R}
\begin{aligned}
2[\sigma_1^{-1} \sigma_2 R_{2n+1} - R_{2n+1} \sigma_2 \sigma_1^{-1}] = \\
= \dfrac{d}{dt_2} [H_{2n}]\sigma_1^{-1} - \sigma_1^{-1}\gamma_* H_{2n}\sigma_1^{-1} + H_{2n} \sigma_1^{-1} \gamma \sigma_1^{-1} - \\
 - \sigma_1^{-1} \dfrac{d}{dt_2} [H_{2n}^*] + \sigma_1^{-1} H^*_{2n} \gamma_*^* \sigma_1^{-1} -
\sigma_1^{-1} \gamma^*  \sigma_1^{-1} H^*_{2n}.
\end{aligned}
\end{equation}
\end{lemma}
\textbf{Proof:} Multiplying equation \eqref{eq:HiHi+1} by $\sigma_1^{-1}$ on the right and using the notation
$H_{2n+1} \sigma_1^{-1} = R_{2n+1} + i M_{2n+1}$ we shall obtain that
\[ \begin{array} {ll}
\sigma_1^{-1} \sigma_2 (R_{2n+1} + i M_{2n+1}) - (R_{2n+1} + i M_{2n+1}) \sigma_2 \sigma_1^{-1} = \\
\quad \quad = \sigma_1^{-1} \sigma_2 R_{2n+1} - R_{2n+1} \sigma_2 \sigma_1^{-1} + i
(\sigma_2 \sigma_1^{-1} M_{2n+1} - M_{2n+1} \sigma_2 \sigma_1^{-1}) =  \\
\quad \quad \quad \quad = 
\dfrac{d}{dt_2} [H_{2n}]\sigma_1^{-1} - \sigma_1^{-1}\gamma_* H_{2n}\sigma_1^{-1} + H_{2n} \sigma_1^{-1} \gamma \sigma_1^{-1}
\end{array} \]
Since $\sigma_1^{-1} \sigma_2 R_{2n+1} - R_{2n+1} \sigma_2 \sigma_1^{-1}$ and
$\sigma_2 \sigma_1^{-1} M_{2n+1} - M_{2n+1} \sigma_2 \sigma_1^{-1}$ are skew adjoint, last equation is equivalent to equality of real and imaginary parts of both sides:
\begin{equation} \label{eq:M} 
\begin{array}{ll}
2i[\sigma_1^{-1} \sigma_2 M_{2n+1} - M_{2n+1} \sigma_2 \sigma_1^{-1}] =
\dfrac{d}{dt_2} [H_{2n}]\sigma_1^{-1} +\sigma_1^{-1} \dfrac{d}{dt_2} [H_{2n}^*] + \\
 + H_{2n} \sigma_1^{-1} \gamma \sigma_1^{-1}  - \sigma_1^{-1}\gamma_* H_{2n}\sigma_1^{-1} - \sigma_1^{-1} H^*_{2n} \gamma_*^* \sigma_1^{-1} + 
\sigma_1^{-1} \gamma^*  \sigma_1^{-1} H^*_{2n},
\end{array} \end{equation}
and \eqref{eq:R}
\begin{eqnarray}
\nonumber 2[\sigma_1^{-1} \sigma_2 R_{2n+1} - R_{2n+1} \sigma_2 \sigma_1^{-1}] =
\dfrac{d}{dt_2} [H_{2n}]\sigma_1^{-1} - \sigma_1^{-1}\gamma_* H_{2n}\sigma_1^{-1} + \\
\nonumber  + H_{2n} \sigma_1^{-1} \gamma \sigma_1^{-1} - \sigma_1^{-1} \dfrac{d}{dt_2} [H_{2n}^*] + \sigma_1^{-1} H^*_{2n} \gamma_*^* \sigma_1^{-1} -
\sigma_1^{-1} \gamma^*  \sigma_1^{-1} H^*_{2n}.
\end{eqnarray}
Let us show that $R_{2n+1}$, defined in equation \eqref{eq:Moments} also satisfies
\eqref{eq:R}. Indeed,
\[ \begin{array} {lllll}
2[\sigma_1^{-1} \sigma_2 R_{2n+1} - R_{2n+1} \sigma_2 \sigma_1^{-1}] = \\
= \sigma_1^{-1} \sigma_2[\sum_{j=0}^{2n} (-1)^{j+1} H_{2n-j} \sigma_1^{-1} H_j^*] - 
[\sum_{j=0}^{2n} (-1)^{j+1} H_{2n-j} \sigma_1^{-1} H_j^*] \sigma_2 \sigma_1^{-1} = \\
= \dfrac{d}{dt_2} [H_{2n}]\sigma_1^{-1} - \sigma_1^{-1}\gamma_* H_{2n}\sigma_1^{-1} + H_{2n} \sigma_1^{-1} \gamma \sigma_1^{-1} - \\
\quad \quad - \sigma_1^{-1} \dfrac{d}{dt_2} [H_{2n}^*] + \sigma_1^{-1} H^*_{2n} \gamma_*^* \sigma_1^{-1} -
\sigma_1^{-1} \gamma^*  \sigma_1^{-1} H^*_{2n}.
\end{array} \]
Performing rearrangement and using \eqref{eq:Moments} for $i=2n-1$, we shall obtain
\[ \begin{array} {lllll}
\sigma_1^{-1} \sigma_2[\sum_{j=0}^{2n} (-1)^{j+1} H_{2n-j} \sigma_1^{-1} H_j^*] - 
[\sum_{j=0}^{2n} (-1)^{j+1} H_{2n-j} \sigma_1^{-1} H_j^*] \sigma_2 \sigma_1^{-1} = \\
= \dfrac{d}{dt_2} [H_{2n} \sigma_1^{-1} - \sigma_1^{-1} H_{2n}^*] - \sigma_1^{-1}\gamma_* H_{2n}\sigma_1^{-1} + H_{2n} \sigma_1^{-1} \gamma \sigma_1^{-1} - H_{2n}\dfrac{d}{dt_2} \sigma_1^{-1} + \\
\quad \quad + \sigma_1^{-1} H^*_{2n} \gamma_*^* \sigma_1^{-1} -
\sigma_1^{-1} \gamma^*  \sigma_1^{-1} H^*_{2n} + \dfrac{d}{dt_2} [\sigma_1^{-1}] H_{2n}^* = \\
= \dfrac{d}{dt_2} [\sum_{j=0}^{2n-1} (-1)^{j+1} H_{2n-1-j} \sigma_1^{-1} H_j^*] - \sigma_1^{-1}\gamma_* H_{2n}\sigma_1^{-1} + H_{2n} \sigma_1^{-1} \gamma \sigma_1^{-1} - \\
\quad \quad  - H_{2n}\dfrac{d}{dt_2} \sigma_1^{-1}
+ \sigma_1^{-1} H^*_{2n} \gamma_*^* \sigma_1^{-1} -
\sigma_1^{-1} \gamma^*  \sigma_1^{-1} H^*_{2n} + \dfrac{d}{dt_2} [\sigma_1^{-1}] H_{2n}^*,
\end{array} \]
Here we can use formula \eqref{eq:HiHi+1} in order to differentiate the expression $H_{2n-1-j} \sigma_1^{-1} H_j^*$.
From \eqref{eq:HiHi+1} it follows that
\begin{equation} \label{eq:DHkS1Hm*}
\begin{array} {lllll}
\dfrac{d}{dt_2} [H_k \sigma_1^{-1} H_m^*] = \\
\quad = [\sigma_1^{-1} \sigma_2 H_{k+1} - 
H_{k+1}\sigma_1^{-1} \sigma_2 + \sigma_1^{-1} \gamma_* H_k - H_k \sigma_1^{-1} \gamma]\sigma_1^{-1} H_m^* + \\
\quad \quad + H_k \dfrac{d}{dt_2} [\sigma_1^{-1}] H_m^* + \\
\quad \quad \quad + H_k \sigma_1^{-1} [H_{m+1}^*\sigma_2 \sigma_1^{-1} -\sigma_2 \sigma_1^{-1}  H^*_{m+1} +
H_m^* \gamma_*^* \sigma_1^{-1} - \gamma^* \sigma_1^{-1} H_m^* ] = \\
\quad =  \sigma_1^{-1} \sigma_2 H_{k+1} \sigma_1^{-1} H_m^* + H_k \sigma_1^{-1} H_{m+1}^*\sigma_2 \sigma_1^{-1} - \\
\quad \quad - H_{k+1} \sigma_1^{-1} \sigma_2\sigma_1^{-1} H_m^*
- H_k \sigma_1^{-1} \sigma_2  \sigma_1^{-1} H^*_{m+1} + \\
\quad \quad \quad +\sigma_1^{-1} \gamma_* H_k \sigma_1^{-1} H_m^* + H_k \sigma_1^{-1}H_m^* \gamma_*^* \sigma_1^{-1}
\end{array} \end{equation}
Now we can insert the formula \eqref{eq:DHkS1Hm*} into the left hand side of the previous expression and 
we shall obtain
\[ \begin{array} {lllll}
\dfrac{d}{dt_2} [\sum_{j=0}^{2n-1} (-1)^{j+1} H_{2n-1-j} \sigma_1^{-1} H_j^*] - \sigma_1^{-1}\gamma_* H_{2n}\sigma_1^{-1} + H_{2n} \sigma_1^{-1} \gamma \sigma_1^{-1} - \\
\quad \quad  - H_{2n}\dfrac{d}{dt_2} \sigma_1^{-1}
+ \sigma_1^{-1} H^*_{2n} \gamma_*^* \sigma_1^{-1} -
\sigma_1^{-1} \gamma^*  \sigma_1^{-1} H^*_{2n} + \dfrac{d}{dt_2} [\sigma_1^{-1}] H_{2n}^* = \\
= \sum_{j=0}^{2n-1} (-1)^{j+1} [\sigma_1^{-1} \sigma_2 H_{2n-j} \sigma_1^{-1} H_j^* + H_{2n-1-j} \sigma_1^{-1} H_{j+1}^*\sigma_2 \sigma_1^{-1} - \\
\quad \quad - H_{2n-j} \sigma_1^{-1} \sigma_2\sigma_1^{-1} H_j^*
- H_{2n-1-j} \sigma_1^{-1} \sigma_2  \sigma_1^{-1} H^*_{j+1} + \\
\quad \quad \quad +\sigma_1^{-1} \gamma_* H_{2n-1-j} \sigma_1^{-1} H_j^* + H_{2n-1-j} \sigma_1^{-1}H_j^* \gamma_*^* \sigma_1^{-1}] - \\
\quad - \sigma_1^{-1}\gamma_* H_{2n}\sigma_1^{-1} + H_{2n} \sigma_1^{-1} \gamma \sigma_1^{-1} - \\
\quad \quad  - H_{2n}\dfrac{d}{dt_2} \sigma_1^{-1}
+ \sigma_1^{-1} H^*_{2n} \gamma_*^* \sigma_1^{-1} -
\sigma_1^{-1} \gamma^*  \sigma_1^{-1} H^*_{2n} + \dfrac{d}{dt_2} [\sigma_1^{-1}] H_{2n}^* =
\end{array} \]
It is easy to recognize telescopic sums in this expression and we obtain that finally it equals to
\[ \begin{array} {lllll}
= \sum_{j=0}^{2n-1} (-1)^{j+1} [\sigma_1^{-1} \sigma_2 H_{2n-j} \sigma_1^{-1} H_j^* +
H_{2n-1-j} \sigma_1^{-1} H_{j+1}^* \sigma_2 \sigma_1^{-1} ] + \\
\quad + H_{2n} \sigma_1^{-1} \sigma_2\sigma_1^{-1} H_0^* - H_0 \sigma_1^{-1} \sigma_2\sigma_1^{-1} H_{2n}^* + \\
\quad \quad + \sum_{j=0}^{2n-1} (-1)^{j+1} [\sigma_1^{-1} \gamma_* H_{2n-1-j} \sigma_1^{-1} H_j^* + H_{2n-1-j} \sigma_1^{-1}H_j^* \gamma_*^* \sigma_1^{-1}] - \\
\quad - \sigma_1^{-1}\gamma_* H_{2n}\sigma_1^{-1} + H_{2n} \sigma_1^{-1} \gamma \sigma_1^{-1} - \\
\quad \quad  - H_{2n}\dfrac{d}{dt_2} \sigma_1^{-1}
+ \sigma_1^{-1} H^*_{2n} \gamma_*^* \sigma_1^{-1} -
\sigma_1^{-1} \gamma^*  \sigma_1^{-1} H^*_{2n} + \dfrac{d}{dt_2} [\sigma_1^{-1}] H_{2n}^* = \\
= \sum_{j=0}^{2n} (-1)^{j+1} [\sigma_1^{-1} \sigma_2 H_{2n-j} \sigma_1^{-1} H_j^* -
H_{2n-j} \sigma_1^{-1} H_j^* \sigma_2 \sigma_1^{-1} ] + \\
\quad + \sum_{j=0}^{2n-1} (-1)^{j+1} [\sigma_1^{-1} \gamma_* H_{2n-1-j} \sigma_1^{-1} H_j^* + H_{2n-1-j} \sigma_1^{-1}H_j^* \gamma_*^* \sigma_1^{-1}] - \\
\quad - \sigma_1^{-1}\gamma_* H_{2n}\sigma_1^{-1} + \sigma_1^{-1} H^*_{2n} \gamma_*^* \sigma_1^{-1} =
\end{array} \]
In other words, we obtain that the compatibility condition \eqref{eq:R} becomes
\[ \begin{array} {lllll}
\sum_{j=0}^{2n-1} (-1)^{j+1} [\sigma_1^{-1} \gamma_* H_{2n-1-j} \sigma_1^{-1} H_j^* + H_{2n-1-j} \sigma_1^{-1}H_j^* \gamma_*^* \sigma_1^{-1}] - \\
\quad - \sigma_1^{-1}\gamma_* H_{2n}\sigma_1^{-1} + \sigma_1^{-1} H^*_{2n} \gamma_*^* \sigma_1^{-1} = 0
\end{array} \]
or after multiplying by $\sigma_1$ from both sides
\[ \gamma_* H_{2n} -  H^*_{2n} \gamma_*^* =
\sum_{j=0}^{2n-1} (-1)^{j+1} [ \gamma_* H_{2n-1-j} \sigma_1^{-1} H_j^* \sigma_1 + \sigma_1 H_{2n-1-j} \sigma_1^{-1}H_j^* \gamma_*^* ],
\]
\qed
\begin{thm}
\label{thm:HiRstrct}
Suppose that we are given moments $H_i(t_2)$, defined in a neighborhood of the 
point $t_2^0\in\mathrm I$. Then there exists $n_0 \le p^2$ 
such that each element $H_{i+1}$ is uniquely determined from $H_0,
\ldots, H_i$ using the algebraic
formulas \eqref{eq:Moments}, and $n_0$ LDEs with arbitrary initial
conditions, 
obtained from \eqref{eq:HiHi+1}. Moreover, $n_0$ equations must be satisfied by the elements
of $\gamma_*(t_2)$.
\end{thm}
\textbf{Remarks: 1.} This theorem is of local nature. The
number "$n_0$" appearing in the theorem
may vary with the point $t_2^0$, but is unchanged in a small
neighborhood of $t_2^0$ by continuity.  \\
\textbf{2.} The proof of the theorem allows to produce an algorithm to determine
explicitly, up to $n_0$ initial conditions, the Markov moments. The
arguments in the proof of the theorem are illustrated on an example in
the following subsection. This example exhibits all the difficulties
present in the general case.\\

\textbf{Proof of \ref{thm:HiRstrct}:}
We have seen in lemma \ref{lemma:Rok}, that the real part of $H_{2n+1}=R+iM$, defined by equation
\eqref{eq:HiHi+1} satisfies \eqref{eq:R}. Let as concentrate on the equation \eqref{eq:M} in order to
solve for $M$. For example, if 
\[ \operatorname{spec}\sigma_1^{-1}\sigma_2 \cap \operatorname{spec}(-\sigma_2\sigma_1^{-1}) = \emptyset
\]
then one can uniquely solve for $M$. But in general, there will be some undefined entries of $M$ in this manner.
In the later case, we can rewrite these equations as linear equations in the entries of $M$.
But before we do that we also rewrite the equation \eqref{eq:M} as following
\begin{eqnarray}
\label{eq:MPrevR} 2i[\sigma_1^{-1} \sigma_2 M_{2n+1} - M_{2n+1} \sigma_2 \sigma_1^{-1}] = \\
\nonumber 
= \dfrac{d}{dt_2} R_{2n} - \sigma_1^{-1}\gamma_* H_{2n}\sigma_1^{-1} - H_{2n} \sigma_1^{-1} \gamma^* \sigma_1^{-1} - \\
\nonumber 
\quad  - \sigma_1^{-1} H^*_{2n} \gamma_*^* \sigma_1^{-1} -
\sigma_1^{-1} \gamma  \sigma_1^{-1} H^*_{2n}
\end{eqnarray}
And further, we are able to rewrite this equation as a system of $p^2$ linear equations in $p^2$
variables $M^{kl}$, $k,l=1,\ldots,p^2$:
\begin{equation}\label{eq:HiHi+1Sys} \left\{ \begin{array}{llll}
\sum \alpha_{kl}^{11} M_{2n+1}^{kl} = \frac{d}{dt_2} R_{2n}^{11} + \sum [ \beta_{kl}^{11} R_{2n}^{kl}] + \sum [ \delta_{kl}^{11} M_{2n}^{kl}]\\
\sum \alpha_{kl}^{12} M_{2n+1}^{kl} = \frac{d}{dt_2} R_{2n}^{12} + \sum [ \beta_{kl}^{12} R_{2n}^{kl}] + \sum [ \delta_{kl}^{12} M_{2n}^{kl}]\\
\hspace{2cm} \vdots \\
\sum \alpha_{kl}^{pp} M_{2n+1}^{kl} = \frac{d}{dt_2} R_{2n}^{pp} + \sum [ \beta_{kl}^{pp} R_{2n}^{kl}] + \sum [ \delta_{kl}^{pp} M_{2n}^{kl}]
\end{array} \right. \end{equation}
Here $\alpha_{kl}^{ij}, \beta_{kl}^{ij}, \delta_{kl}^{ij}$ are functions of $t_2$ derived from the vessel parameters.
Using basic algebra manipulations, it is enough to find a maximally 
independent subset of equations in the left side of
(\ref{eq:HiHi+1Sys}), and to obtain a system of $p^2-n_0$ independent equations
\begin{equation} \label{eq:Hi+1Alg}
\left\{ \begin{array}{ll} 
\sum_{kl} \alpha^1_{kl}(t_2) M_{2n+1}^{kl} = f_1 (\dfrac{d}{dt_2} R_{2n}^{kl}, R_{2n}^{kl}, M_{2n}^{kl}) \\
\sum_{kl} \alpha^2_{kl}(t_2) M_{2n+1}^{kl} = f_2 (\dfrac{d}{dt_2} R_{2n}^{kl}, R_{2n}^{kl}, M_{2n}^{kl}) \\
\hspace{2cm} \vdots \\
\sum_{kl} \alpha^{p^2-n_0}_{kl}(t_2) M_{2n+1}^{kl} = f_{p^2-n_0} (\dfrac{d}{dt_2} R_{2n}^{kl}, R_{2n}^{kl}, M_{2n}^{kl})
\end{array} \right. \end{equation}
for linear functions $f_j$ and a system of linear dependent ones
\[ \left\{ \begin{array}{ll}  
0 = g_1(\dfrac{d}{dt_2} R_{2n}, R_{2n}^{kl}, M_{2n}^{kl}) \\
0 = g_2(\dfrac{d}{dt_2} R_{2n}, R_{2n}^{kl}, M_{2n}^{kl}) \\
\vdots \\
0 = g_{n_0}(\dfrac{d}{dt_2} R_{2n}, R_{2n}^{kl}, M_{2n}^{kl})
\end{array} \right.\]
for linear functions $g_j$. 

On the other hand, using the same considerations for $H_{2n+2}$, we shall obtain
$n_0$ equations of the same type, which are additional restrictions on $M_{2n+1}$
\begin{equation} \label{eq:Hi+1D}
0 = g_j(\dfrac{d}{dt_2} M_{2n+1}, M_{2n+1}, R_{2n+1}), j=1,\ldots,n_0.
\end{equation}
Notice that we obtain algebraic (\ref{eq:Hi+1Alg}) and differential (\ref{eq:Hi+1D}) 
equations which are independent (since all algebraic and differentials are, and they are of different nature). 
So, if \textit{all} elements of $M_{2n+1}$ appear at these equations we shall obtain that one can
uniquely solve them up to $n_0$ initial conditions:
\begin{equation} \label{eq:AllBefore}
\left\{ \begin{array}{lllll}
\sum_{kl} \alpha^i_{kl}(t_2) M_{2n+1}^{kl} = f_i (\dfrac{d}{dt_2} R_{2n}^{kl}, R_{2n}^{kl}, M_{2n}^{kl}) \\
\hspace{5cm}  i=1,\ldots, p^2-n_0, \\
0 = g_i(\dfrac{d}{dt_2} M_{2n+1}, M_{2n+1}, R_{2n+1}), \\
\hspace{5cm} i=1,\ldots, n_0
\end{array} \right. \end{equation}
(See example of NLS at section \ref{sec:NLS}).
If this is not the case, some of the elements, say $p_0\leq n_0$ will
not appear in the equations \eqref{eq:AllBefore}. As a result, using the algebraic expressions \eqref{eq:Hi+1Alg}
and plugging them into differential ones, \eqref{eq:Hi+1D} we shall obtain $p_0\leq n_0$ 
differential equation of second order on the elements $R_{2n}$ and $M_{2n}$:
\[ 0 = h_i(\dfrac{d^2}{dt_2^2} R_{2n}, \dfrac{d}{dt_2} R_{2n}, R_{2n}, \dfrac{d}{dt_2} M_{2n}, M_{2n}),
	 j=1,\ldots,p_0.
\]
By induction the same type of equation will hold for $H_{2n+1}$ (for $j=1,\ldots,p_0$):
\begin{equation} \label{eq:Hi+1D2}
0 = h_i(\dfrac{d^2}{dt_2^2} M_{2n+1}, \dfrac{d}{dt_2} M_{2n+1}, M_{2n+1}, \dfrac{d}{dt_2} R_{2n+1}, R_{2n+1}).
\end{equation}
Finally, we shall obtain the following system of differential equations
\begin{equation} \label{eq:All}
\left\{ \begin{array}{lrlrlr}
\sum_{kl} \alpha^i_{kl}(t_2) M_{2n+1}^{kl} = f_i (\dfrac{d}{dt_2} R_{2n}^{kl}, R_{2n}^{kl}, M_{2n}^{kl}) \\
\hspace{5cm}  i=1,\ldots, p^2-n_0, \\
0 = g_i(\dfrac{d}{dt_2} M_{2n+1}, M_{2n+1}, R_{2n+1}), \\
\hspace{5cm} i=1,\ldots, n_0-p_0 \\
0 = h_i(\dfrac{d^2}{dt_2^2} M_{2n+1}, \dfrac{d}{dt_2} M_{2n+1}, M_{2n+1}, \dfrac{d}{dt_2} R_{2n+1}, R_{2n+1}), \\
\hspace{5cm} i=1,\ldots, p_0
\end{array} \right. \end{equation}
\qed

\noindent\textbf{Remark:} This theorem means that the equations \eqref{eq:HiHi+1} and \eqref{eq:Moments} are compatible,
i.e. one can always solve them, once $n_0$ equations are satisfied for $\gamma_*(t_2)$. But if we consider
the differential equation \eqref{eq:HiHi+1} only and ignore \eqref{eq:Moments}, then one obtain the following 
\begin{cor} \label{cor:Hn+1HnDE}
Suppose that Moments of $H_i$ satisfy the differential equations \eqref{eq:HiHi+1}, then there exist
$K_0(t_2), K_1(t_2), K_2(t_2), K_3(t_2)$ depending on vessel parameters only such that
\begin{equation} \label{eq:Hn+1HnDE}
 \dfrac{d}{dt_2} H_{n+1} = K_0 H_n + K_1 \dfrac{d}{dt_2} H_n + K_2 \dfrac{d^2}{dt_2^2} H_n+ K_3 \dfrac{d^3}{dt_2^3} H_n
\end{equation}
\end{cor}
\textbf{Proof:} All we have to do it to imitate the proof of theorem \ref{thm:HiRstrct} ignoring the additional
constrain \eqref{eq:Moments} on the real/imaginary parts of the moments. Denoting by $H_{n}^{ij}$ the entries
of $H_n$, we will obtain from \eqref{eq:HiHi+1} $p^2$ equations (like \eqref{eq:HiHi+1Sys})
\begin{equation*} \left\{ \begin{array}{llll}
\sum \alpha_{kl}^{11} H_{n+1}^{kl} = \dfrac{d}{dt_2} H_n^{kl} + \sum[\beta_{kl}^{11} H_n^{kl}]\\
\sum \alpha_{kl}^{12} H_{n+1}^{kl} = \dfrac{d}{dt_2} H_n^{kl} + \sum[\beta_{kl}^{12} H_n^{kl}]\\
\hspace{2cm} \vdots \\
\sum \alpha_{kl}^{pp} H_{n+1}^{kl} = \dfrac{d}{dt_2} H_n^{kl} + \sum[\beta_{kl}^{pp} H_n^{kl}]
\end{array} \right. \end{equation*}
and analogously to deriving the formulas \eqref{eq:AllBefore}, we will obtain
\begin{equation*} 
\left\{ \begin{array}{lllll}
\sum_{kl} \alpha^i_{kl}(t_2) H_{n+1}^{kl} = f_i (\dfrac{d}{dt_2} H_n, H_n), &  i=1,\ldots, p^2-n_0, \\
0 = g_i(\dfrac{d}{dt_2} H_{n+1}, H_{n+1}), & i=1,\ldots, n_0
\end{array} \right. \end{equation*}
So, if these equations are independent, differentiating the first one, we shall obtain 
\[ \dfrac{d}{dt_2} H_{n+1} = K_0 H_n + K_1 \dfrac{d}{dt_2} H_n + K_2 \dfrac{d^2}{dt_2^2} H_n.
\]
Otherwise, there will be $p_0$ equations analogous to \eqref{eq:Hi+1D2} of the form
\[ 0 = h_i(\dfrac{d^2}{dt_2^2} H_{n+1}, \dfrac{d}{dt_2} H_{n+1}, H_{n+1}), j=1,\ldots,p_0.
\]
Moreover, the independent entries of $H_{n+1}$ will appear with one differentiation only, because,
otherwise they would appear at the previous stage at the equation $0 = g_i(\dfrac{d}{dt_2} H_{n+1}, H_{n+1})$,
and would be derived from $n_0$ differential equations. So, this last equation means that the remaining $p_0$
elements are found from
\[ \dfrac{d}{dt_2} H_{n+1}^{ij} = L_0 H_{n+1} + L_1 \dfrac{d}{dt_2} H_{n+1} +
L_2 \dfrac{d^2}{dt_2^2} H_{n+1}
\]
so that only the known entries appear at the left hand side. From here it immediately follows the corollary.
\qed

We want to present next a necessary restriction on $\gamma_*(t_2)$, derived from the existence of a finite
dimensional vessel:
\begin{thm}\label{thm:gamma*DE}
Let $\sigma_1, \sigma_2, \gamma, \gamma_*$ be vessel parameters, and $\mathfrak V$ a 
finite dimensional vessel corresponding to the m parameters. Then the entries of the function $\gamma_*$ satisfies
a polynomial differential equation of finite order with coefficients in the differential ring $\boldsymbol{\mathcal R}$, 
generated by (the entries of) $\sigma_1, \sigma_1^{-1}, \sigma_2, \gamma$.
\end{thm}
\textbf{Proof:} Suppose that the transfer function of the vessel $\mathfrak V$, defined in \eqref{eq:DefV} is
\[ S(\lambda,t_2) = I - B^*(t_2) \mathbb X^{-1}(t_2) (\lambda I - A_1)^{-1} B(t_2) \sigma_1(t_2) = 
I - \sum_{i=0}^\infty \dfrac{H_i}{\lambda^{i+1}} ,
\]
and the linkage condition is
\[ \sigma_1^{-1}(t_2) \gamma_*(t_2) = \sigma_1^{-1}(t_2) \gamma(t_2) + [\sigma_1^{-1}(t_2) \sigma_2(t_2), H_0(t_2)].
\]
Notice that if we differentiate this formula and use the equation for the derivative of $H_0(t_2)$ from the equation
\eqref{eq:HiHi+1}, we shall get
\[ \begin{array}{lll}
\dfrac{d}{dt_2}\big(\sigma_1^{-1}(t_2)\gamma_*(t_2)\big) & = \dfrac{d}{dt_2}\big(\sigma_1^{-1}(t_2)\gamma(t_2)\big) + 
[\dfrac{d}{dt_2}(\sigma_1^{-1}(t_2) \sigma_2(t_2)), H_0(t_2)] + \\
& + [\sigma_1^{-1}(t_2) \sigma_2(t_2), \dfrac{d}{dt_2} H_0(t_2)] = \\ 
& = \dfrac{d}{dt_2}\big(\sigma_1^{-1}(t_2)\gamma(t_2)\big) + f_{00}(H_0(t_2)) + f_{01}(H_1(t_2))
\end{array} \]
for linear in moments functions $f_{00}, f_{01}$ with coefficients depending on $\boldsymbol{\mathcal R}$ and
$\gamma_*$.
Similarly, differentiating this expression and using formula \eqref{eq:HiHi+1} for $\dfrac{d}{dt_2} H_0(t_2)$ and
$\dfrac{d}{dt_2} H_1(t_2)$, we shall obtain that the second derivative is
\[ \begin{array}{lll}
\dfrac{d^2}{dt_2^2}\big(\sigma_1^{-1}(t_2)\gamma_*(t_2)\big) = \\
\quad = \dfrac{d^2}{dt_2^2}\big(\sigma_1^{-1}(t_2)\gamma(t_2)\big) +  f_{10}(H_0) + f_{11}(H_1(t_2))+ f_{12}(H_2(t_2)).
\end{array} \]
for linear in the moments functions $f_{10}, f_{11}, f_{12}$
with coefficients depending on $\boldsymbol{\mathcal R}$ and
$\gamma_*, \dfrac{d}{dt_2}\gamma_*$. Continuing this differentiation further at each step $i$ we shall obtain an equation of the form
\[ \dfrac{d^i}{dt_2^i}\big(\sigma_1^{-1}(t_2)\gamma_*(t_2)\big) =
\dfrac{d^i}{dt_2^i}\big(\sigma_1^{-1}(t_2)\gamma(t_2)\big) +
\sum_{j=0}^i f_{ij} (H_j),
\]
where $f_{ij}$ is a linear function of $H_j$ with coefficients, depending on $\boldsymbol{\mathcal R}$
and first $i-1$ derivatives of $\gamma_*$ (this can be immediately seen by the induction). But at some stage the
moments start to repeat themselves due to equation \eqref{eq:HLin}. So, taking $K$ derivatives we shall
obtain equations of the following form
\begin{equation}\label{eq:gamma*SE} \left\{ \begin{array}{lllllll}
\sigma_1^{-1}(t_2) \gamma_*(t_2) = \sigma_1^{-1}(t_2) \gamma(t_2) + f_{00}(H_0(t_2)) \\
\dfrac{d}{dt_2}\big(\sigma_1^{-1}(t_2)\gamma_*(t_2)\big) = \dfrac{d}{dt_2}\big(\sigma_1^{-1}(t_2)\gamma(t_2)\big) + f_{10}(H_0(t_2)) + f_{11}(H_1(t_2)) \\
\dfrac{d^2}{dt_2^2}\big(\sigma_1^{-1}(t_2)\gamma_*(t_2)\big) = \dfrac{d^2}{dt_2^2}\big(\sigma_1^{-1}(t_2)\gamma(t_2)\big) +  \sum\limits_{j=0}^2 f_{2j} (H_j) \\
\hspace{4cm} \vdots \\
\dfrac{d^N}{dt_2^N}\big(\sigma_1^{-1}(t_2)\gamma_*(t_2)\big) =
\dfrac{d^N}{dt_2^N}\big(\sigma_1^{-1}(t_2)\gamma(t_2)\big) +
\sum\limits_{j=0}^N f_{Nj} (H_j) \\
\hspace{4cm} \vdots \\
\dfrac{d^K}{dt_2^K}\big(\sigma_1^{-1}(t_2)\gamma_*(t_2)\big) =
\dfrac{d^K}{dt_2^K}\big(\sigma_1^{-1}(t_2)\gamma(t_2)\big) +
\sum\limits_{j=0}^N f_{Kj} (H_j) \\
\end{array} \right.
\end{equation}
Suppose that that the dimension of the inner state is $n$ (i.e. $\dim\mathcal H= n$),
then we get that each of the matrices $H_j$ has $n^2$ entries and there is the total number of $n^2 (N+1)$ 
entries for the moments $H_0,\ldots,H_N$. So, taking 
"enough" derivatives of $\gamma_*(t_2)$ (i.e., taking $K$ so that $K \dim (\mathcal E)^2> n^2(N+1)$)
and eliminating all the entries of the moments, we shall obtain a finite number of
polynomial differential equation for the entries of $\gamma_*(t_2)$.
\qed

\textbf{Remark:} From this theorem is follows that $\gamma_*$ satisfies an equation of
the form
\[ P(x,x',x'',\ldots,x^{(K)}) = 0, \] 
where $P(x_0,x_1,x_2,\ldots,x_K)$ is a non-commutative polynomial with coefficients in $\boldsymbol{\mathcal R}$.

Using an outer space transformation it is possible to bring the matrix $\sigma_1^{-1}(t_2) \sigma_2(t_2)$ to
a simpler form. Suppose that there is an invertible matrix $V(t_2): \mathcal E \rightarrow \mathcal E$
such that
\[ \sigma_1^{-1}(t_2) \sigma_2(t_2) = V(t_2) L(t_2) V^{-1}(t_2)
\]
for a "simple" matrix $L(t_2)$, where $L(t_2^0)$ is in a Jordan block form.
Define new operators (omitting $t_2$-dependence)
\begin{equation} \label{eq:OuterChange}
\widetilde \sigma_1 = V^* \sigma_1 V, \quad \widetilde \sigma_2 = V^* \sigma_2 V, \quad
\widetilde \gamma = V^* \gamma V - \sigma_1 V', \quad \widetilde \gamma_* = V^* \gamma_* V - \sigma_1 V'
\end{equation}
then simple calculations show that $\widetilde \sigma_1(t_2)$, $\widetilde \sigma_2(t_2)$, $\widetilde \gamma(t_2)$,
$\widetilde \gamma_*(t_2)$ are vessel parameters in the same interval $\mathrm I$ and the collection of operators
and spaces
\[ \mathfrak{\widetilde V} = (A_1, B(t_2), \mathbb X(t_2); \widetilde\sigma_1(t_2), 
\widetilde\sigma_2(t_2), \widetilde\gamma(t_2), \widetilde\gamma_*(t_2);
\mathcal{H},\mathcal{E}), 
\]
is again a Vessel. By multiplying the original vessel conditions of $\mathfrak V$ by $V$ or on $V^*$ we check that
the the vessel conditions are satisfied for $\mathfrak{\widetilde V}$.
But for this new vessel $\mathfrak{\widetilde V}$ we obtain that
\[ \widetilde\sigma_1^{-1}(t_2) \widetilde\sigma_2(t_2) =
V^{-1}(t_2) \sigma_1^{-1}(t_2) \sigma_2(t_2) V(t_2) = L(t_2).
\]
The first step to understand the equations arising in theorem \eqref{thm:HiRstrct} we study the cosntant
case.
\subsection{Moment equations for diagonal constant $\sigma_1^{-1}\sigma_2$}
At the special case $\sigma_1^{-1}\sigma_2 = \operatorname{diag}[s_1,s_2,\ldots,s_E]$ ($E=\dim \mathcal E$)
for some constants $s_i\in\mathbb C$ we can explicitly
find the formulas for the moments using equation \eqref{eq:HiHi+1}.

First we consider the linkage condition
\eqref{eq:H0}, which may be considered as the "$-1$" equation of \eqref{eq:HiHi+1}. Taking diagonal form for
$\sigma_1^{-1}(t_2) \sigma_2(t_2)$ and denoting the $k,j$-th entry of $H_0(t_2)$ by $H_0^{kj}$, 
we shall obtain
\[ \begin{array}{lll}
\sigma_1^{-1}(t_2)[\gamma_*(t_2) - \gamma(t_2)]  = [\sigma_1^{-1}(t_2)\sigma_2(t_2), H_0(t_2)] =
\bbmatrix{H_0^{kj}(t_2) (s_k-s_j)}.
\end{array} \]
Cases for which $s_k-s_j =0$ give restrictions on the entries of \linebreak
$\sigma_1^{-1}(t_2)[\gamma_*(t_2) - \gamma(t_2)]$ and other cases uniquely determine entries of
$H_0(t_2)$:
\begin{eqnarray}
e_k \sigma_1^{-1}(t_2)[\gamma_*(t_2) - \gamma(t_2)] e_j^t = 0, & s_k - s_j = 0, \\
\label{eq:H0kj} H_0^{kj}(t_2) = \dfrac{e_k \sigma_1^{-1}(t_2)[\gamma_*(t_2) - \gamma(t_2)] e_j^t}{s_k-s_j},
& s_k - s_j \neq 0.
\end{eqnarray}
But there more equations arise when we consider \eqref{eq:HiHi+1} for $i=0$:
\[ \begin{array}{lll}
[\sigma_1^{-1}\sigma_2, H_1(t_2)] = \bbmatrix{H_1^{kj}(t_2) (s_k-s_j)} = \\
= \dfrac{d}{dt_2} H_0(t_2) -
\sigma_1^{-1}(t_2) \gamma_*(t_2) H_0(t_2) + H_0(t_2) \sigma_1^{-1}(t_2) \gamma(t_2)
\end{array} \]
Here again if we focus on the cases $s_k-s_j=0$, we shall obtain that the following differential equations
must hold (suppressing notation for $t_2$ dependence)
\begin{eqnarray*}
\dfrac{d}{dt_2} H_0^{kj} - e_k [
\sigma_1^{-1} \gamma_* H_0 + H_0 \sigma_1^{-1} \gamma]e_j^t = 0, & s_k-s_j=0, \\
H_1^{kj}(t_2) = \frac{1}{s_k-s_j} e_k [\dfrac{d}{dt_2} H_0 -
\sigma_1^{-1} \gamma_* H_0 + H_0 \sigma_1^{-1} \gamma]e_j^t, & s_k-s_j\neq 0
\end{eqnarray*}
The first line gives a linear differential equation for the elements $H_0^{kj}$ in terms
of vessel parameters and other known entries of $H_0$:
\[ \dfrac{d}{dt_2} H_0^{kj} = e_k \sigma_1^{-1} \gamma_* H_0 e_j^t -
e_k H_0 \sigma_1^{-1} \gamma e_j^t, \quad 1 \leq k \leq r,
\]
Finally we obtain the following
\begin{lemma} Suppose that 
\[ \sigma_1^{-1}\sigma_2 = \operatorname{diag}[s_1,s_2,\ldots,s_E],
\]
then the entries of $\gamma, \gamma_*$ satisfy the following equations
\begin{equation} \label{eq:Sig1Sig2Diag}
e_k \sigma_1^{-1}(t_2)[\gamma_*(t_2) - \gamma(t_2)] e_k^t = 0, \quad \text{ when }s_k-s_j=0,
\end{equation}
where $e_k$ stay for the standard elementary row vector with 1 at $k$-th place and $0$ everywhere.
\end{lemma}
An interesting question arises, is whether any continuous vessel parameters satisfying the conditions of this
lemma gives rise to a transfer function.

\subsection{Moment equations for Jordan-block constant $\sigma_1^{-1}\sigma_2$}
Suppose now that $\sigma_1^{-1}\sigma_2$ is a Jordan block matrix with arbitrary eigenvalue $s\in\mathbb C$.
\[ \sigma_1^{-1}\sigma_2 = \operatorname{Jordan}(s,r) =
\bbmatrix{s & 1 & \ldots & 0 & 0\\
0 & s & \ldots & 0 & 0 \\
\vdots & \vdots & \ddots & \vdots & \vdots \\
0 & 0 & \ldots & s & 1 \\
0 & 0 & \ldots & 0 & s}
\]
then we can again consider the linkage condition \eqref{eq:H0}:
\[ \begin{array}{lll}
\sigma_1^{-1}(t_2)[\gamma_*(t_2) - \gamma(t_2)]  = [\sigma_1^{-1}(t_2)\sigma_2(t_2), H_0(t_2)] = \\
= [\operatorname{Jordan}(s,r), H_0(t_2)] = [\operatorname{Jordan}(0,r), H_0(t_2)] = \\
= \bbmatrix{H_0^{21} & H_0^{22}-H_0^{11} & H_0^{23}-H_0^{12} & \ldots & H_0^{2r}-H_0^{1,r-1} \\
H_0^{31} & H_0^{32}-H_0^{21} & H_0^{33}-H_0^{22} & \ldots & H_0^{3r}-H_0^{2,r-1} \\
H_0^{41} & H_0^{42}-H_0^{31} & H_0^{43}-H_0^{32} & \ldots & H_0^{4r}-H_0^{3,r-1} \\
\vdots & \vdots & \vdots & \vdots & \vdots \\
H_0^{r-1,1} & H_0^{r-1,2}-H_0^{r-2,1} & H_0^{r-1,3}-H_0^{r-2,2} & \ldots & H_0^{r-1,r}-H_0^{r-2,r-1} \\
H_0^{r1} & H_0^{r2}-H_0^{r-1,1} &  H_0^{r3}-H_0^{r-1,2} & \ldots & H_0^{rr}-H_0^{r-1,r-1} \\
0  & - H_0^{r1} & - H_0^{r2} & \ldots & -H_0^{r,r-1} \\
},
\end{array} \]
where the value $s$ on the main diagonal is canceled. It turns out that sums of elements of each diagonal under the
main one is 0. This may be rewritten also in the following form
\[ \TR [\sigma_1^{-1}(\gamma_*(t_2) - \gamma(t_2)) \operatorname{Jordan}(0,r)^k ] = 0, \quad k=0,1,\ldots,r-1. 
\]
From the other equations we find expressions for the entries of $H_0(t_2)$ using its first row as parameters.
In other words, from the linkage condition \eqref{eq:H0} we find $r$ restrictions on $\sigma_1^{-1}(\gamma_*(t_2) - \gamma(t_2))$ and all the entries of $H_0(t_2)$, while the entries of the first row 
$H_0^{11}, \ldots, H_0^{1r}$ are unknown.

Considering next the equation \eqref{eq:HiHi+1}, we will similarly obtain that
\[ \TR[\big( \dfrac{d}{dt_2} H_0 - \sigma_1^{-1} \gamma_* H_0 + H_0 \sigma_1^{-1} \gamma \big)\operatorname{Jordan}(0,r)^k] = 0, k=0,1,\ldots,r-1.
\]
These equations may be considered as additional equations on the first row of $H_0(t_2)$, which was unknown before.
It is also possible that these differential equations will not determine the first row of $H_0$, in this case we shall
obtain additional differential restrictions on the elements of $\sigma_1^{-1}(\gamma_*(t_2) - \gamma(t_2))$ (see
\eqref{eq:PiBeta} below). In this case, the elements of the first row of $H_0$ will be determined from differential equations on the next level, involving $H_1(t_2)$ (see \eqref{eq:H021} below).

\subsection{\label{sec:SL}Sturm Liouville vessel parameters}
The following example was extensively studied in
\cite{bib:SLVessels}. It deals with the Sturm Liouville differential equation
\[ \frac{d^2}{dt_2^2} y(t_2) - q(t_2) y(t_2) = \lambda y(t_2),
\]
with the spectral parameter $\lambda$. The parameter $q(t_2)$ is
usually called the \textit{potential}.
For $q(t_2) = 0$ this problem is easily solved by exponents and in this case we shall call this equation
\textit{trivial}. In \cite{bib:SLVessels} one connects solutions of the more general problem with
non trivial $q(t_2)$ to the trivial one.
\begin{defn} \label{def:SLpar}the Sturm Liouville vessel parameters are given by
\[ \begin{array}{lll}
\sigma_1 = \bbmatrix{0 & 1 \\ 1 & 0},   \sigma_2 = \bbmatrix{1 & 0 \\ 0 & 0}, \gamma = \bbmatrix{0 & 0 \\ 0 & i},\\
\gamma_*(t_2)= \bbmatrix{-i \pi_{11}(t_2) & -\beta(t_2) \\ \beta(t_2) & i}
\end{array} \]
for \textbf{real} valued continuous functions $\pi_{11}(t_2), \beta(t_2)$.
\end{defn}
Notice that these parameters correspond to the Jordan block form of 
$\sigma_1^{-1}\sigma_2 = \bbmatrix{0 & 0 \\ 1 & 0}$.

The input compatibility differential equation (\ref{eq:InCC}) is equivalent to
\[ \left\{ \begin{array}{lll}
\lambda u_1(\lambda,t_2) - \frac{\partial}{\partial t_2}u_2(\lambda,t_2) = 0 \\
- \frac{\partial}{\partial t_2}u_1(\lambda,t_2) + i u_2(\lambda,t_2) = 0
\end{array}\right.
\]
where we denote $u_\lambda(t_2) = \bbmatrix{u_1(\lambda,t_2)\\ u_2(\lambda,t_2)}$. From the second equation one finds that
$u_2(\lambda,t_2) = -i \frac{\partial}{\partial t_2}u_1(\lambda,t_2)$
and 
plugging it back to the first
equation, we shall obtain the trivial Sturm Liouville differential
equation with the spectral parameter $i\lambda$ for
$u_1(\lambda,t_2)$:
\[ \frac{\partial^2}{\partial t_2^2} u_1(\lambda,t_2) = i\lambda u_1(\lambda,t_2).
\]
For the output $y_\lambda(t_2) = \bbmatrix{y_1(\lambda,t_2)\\y_2(\lambda,t_2)}$, we shall obtain that (\ref{eq:OutCC}) is
equivalent to the system of equations
\[ \left\{ \begin{array}{lll}
(\lambda - i \pi_{11}(t_2))y_1(\lambda,t_2) - (\frac{\partial}{\partial t_2} + \beta(t_2))y_2(\lambda,t_2) = 0 \\
(\beta(t_2)- \frac{\partial}{\partial t_2})y_1(\lambda,t_2) + i y_2(\lambda,t_2) = 0
\end{array}\right.,
\]
from which we immediately obtain that $y_2(\lambda,t_2) = i 
(\beta(t_2)- \frac{\partial}{\partial t_2})y_1(\lambda,t_2)$ and
plugging it into the first equation
\[ \frac{\partial^2}{\partial t_2^2} y_1(\lambda,t_2) - (\pi_{11}(t_2) + \beta'(t_2) + \beta^2(t_2))y_1(\lambda,t_2) = i \lambda y_1(\lambda,t_2), \]
which means that $y_1(\lambda,t_2)$ satisfies the Sturm Liouville 
differential equation with the spectral parameter $i\lambda$ and
the potential $q(t_2) = (\pi_{11}(t_2) + \beta'(t_2) + \beta^2(t_2))$.\\

The first equation (\ref{eq:H0}) considered for
$H_0 = \bbmatrix{H_0^{11} & H_0^{12}\\H_0^{21} & H_0^{22}}$ becomes
\[ \bbmatrix{-i \pi_{11}(t_2) & -\beta(t_2) \\ \beta(t_2) & 0} =
\bbmatrix{H_0^{11}-H_0^{22} & H_0^{12} \\ -H_0^{12} & 0}.
\]
from where we conclude that
\begin{equation} \label{eq:SLH0} H_0^{12} = - \beta(t_2) ,\quad
  H_0^{11}-H_0^{22} 
= -i \pi_{11}(t_2).
\end{equation}
Let us consider the differential equation (\ref{eq:HiHi+1}), where we
use the notation $H_1 = \bbmatrix{H_1^{11} & H_1^{12}\\H_1^{21} &
  H_1^{22}}$. Substituting the expressions for the vessel parameters, we shall obtain
\[ \left\{ \begin{array}{lll}
\frac{d}{dt_2} H_0^{11} -\beta H_0^{11} - i H_0^{21} & = -H_1^{12}, \\
\frac{d}{dt_2} H_0^{12} -\beta H_0^{12} + i (H_0^{11} - H_0^{22}) & = 0, \\
\frac{d}{dt_2} H_0^{21} + i \pi_{11} H_0^{11} +\beta H_0^{21} & =  H_1^{11}-H_1^{22}, \\
\frac{d}{dt_2} H_0^{22} + i \pi_{11} H_0^{12} + \beta H_0^{22} + i H_0^{21} & = H_1^{12}.
\end{array} \right.
\]
and consequently, using the formulas \eqref{eq:SLH0} the second equation results in 
\begin{equation}
\label{eq:PiBeta} \pi_{11} = \frac{d}{dt_2} \beta  - \beta^2.
\end{equation}
Together the first and the fourth equations give
\[ H_1^{12} = -(\frac{d}{dt_2} H_0^{11} -\beta H_0^{11} - i H_0^{21}) =
\frac{d}{dt_2} H_0^{22} + i \pi_{11} H_0^{12} + \beta H_0^{22} + i H_0^{21}
\]
from where we obtain using (\ref{eq:SLH0})
\begin{equation} \label{eq:SLH02}
\frac{d}{dt_2}[H_0^{11} + H_0^{22}] = 0 \Rightarrow H_0^{11} + H_0^{22} = C\in\mathbb C,
\end{equation}
Additionally, $H_0$ has to satisfy $H_0 = \sigma_1^{-1} H_0^* \sigma_1$. Using this relation we shall obtain that
\[ \begin{array}{lll}
\bbmatrix{H_0^{11} & H_0^{12}\\H_0^{21} & H_0^{22}} = \bbmatrix{0 & 1 \\ 1 & 0}\bbmatrix{(H_0^{11})^* & (H_0^{21})^* \\ (H_0^{12})^* & (H_0^{22})^*}  \bbmatrix{0 & 1 \\ 1 & 0} \Rightarrow \\
\Rightarrow \left\{ \begin{array}{llll}
H_0^{11} = (H_0^{22})^* \\
H_0^{12} = (H_0^{12})^* \\
H_0^{21} = (H_0^{21})^* \\
H_0^{22} = (H_0^{11})^* \\
\end{array} \right.
\end{array} \]
from where we conclude that
\begin{equation} \label{eq:SLH0A}
 H_0^{11} = (H_0^{22})^*, \quad H_0^{12} = (H_0^{12})^*, \quad H_0^{21} = (H_0^{21})^*.
\end{equation}
As we can see $H_0^{21}=h_0^{21}$ is a real valued (arbitrary at this stage) function, and $H_0$ is as follows
\begin{equation}\label{eq:H0Mat} 
H_0 = \bbmatrix{ \dfrac{r -i\pi_{11}}{2} & -\beta \\ h_0^{21} &
  \dfrac{r +i\pi_{11}}{2}},
\end{equation}
where $ r\in\mathbb R$.\\

Let us perform the same calculations for $H_1$. From the algebraic equation (\ref{eq:Moments})
\[ H_1 \sigma_1^{-1} + \sigma_1^{-1} H_1^* = -H_0 \sigma_1^{-1} H_0^*
\]
we obtain that
\[ \bbmatrix{H_1^{12} + (H_1^{12})^* & H_1^{11} + (H_1^{22})^* \\H_1^{22} + (H_1^{11})^* & H_1^{21} + (H_1^{21})^*} = -H_0 \sigma_1^{-1} H_0^*
\]
From the equation (\ref{eq:HiHi+1}) with $i=0$ we obtain as before that
\begin{eqnarray}
\label{eq:SLH11} H_1^{11}-H_1^{22} = \frac{d}{dt_2} H_0^{21} + i \pi_{11} H_0^{11} +\beta H_0^{21}, \\
\label{eq:SLH12} H_1^{12} = i H_0^{21} - \frac{d}{dt_2} H_0^{11} + \beta H_0^{11}.
\end{eqnarray}
and the same equation (\ref{eq:HiHi+1}) with $i=1$ produces similarly to the previous case
\begin{eqnarray}
\label{eq:SLH13}\frac{d}{dt_2} H_1^{12} -\beta H_1^{12} + i (H_1^{11} - H_1^{22}) = 0, \\
\label{eq:SLH11+22} \frac{d}{dt_2} (H_1^{11} + H_1^{22}) = -  i \pi_{11} H_1^{12} + \beta (H_1^{11} - H_1^{22}).
\end{eqnarray}
Plugging (\ref{eq:SLH11}) and (\ref{eq:SLH12}) into (\ref{eq:SLH13}), we shall obtain that $H_0^{21}=h_0^{21}$ have to satisfy the following
differential equation of the first order:
\[ \begin{array}{lll}
\dfrac{d}{dt_2} (i H_0^{21} - \frac{d}{dt_2} H_0^{11} + \beta H_0^{11}) - \\
- \beta (i H_0^{21} - \frac{d}{dt_2} H_0^{11} + \beta H_0^{11}) + \\
+ i \frac{d}{dt_2} H_0^{21} - \pi_{11} H_0^{11} +i \beta H_0^{21} = 0.
\end{array} \]
or after cancellations
\begin{equation} \label{eq:H021}
 2i (H_0^{21})' = \frac{d^2}{dt_2^2} H_0^{11} - 2 \beta \frac{d}{dt_2}H_0^{11}.
\end{equation}
Inserting here the expressions for $H_0^{ij}$ appearing in \eqref{eq:H0Mat} we shall obtain that
the real part of the last equality can be derived from the equation \eqref{eq:PiBeta}:
\[ \dfrac{r}{2} [\beta'-\beta^2-\pi_{11}] = 0.
\]
The imaginary part gives the following equation
\[ 4 \frac{d}{dt_2} h_0^{21} = \frac{d}{dt_2} (\pi_{11}\beta) + \beta \pi_{11}'-
\beta^2\pi_{11}-\pi_{11}^2-\pi_{11}''.
\]
Suppose (see \cite{bib:SLVessels})
that there exists a function $\tau$ such that $\beta=-\dfrac{\tau'}{\tau}$. Then using
\eqref{eq:PiBeta} $\pi_{11} = -\dfrac{\tau''}{\tau}$ and inserting these equations into the
formula for $\frac{d}{dt_2} h_0^{21}$ we obtain that
\[ 4 \frac{d}{dt_2} h_0^{21} = \dfrac{\tau^{(4)}}{\tau} - \big(\dfrac{\tau''}{\tau}\big)^2.
\]
Let us write down the formulas for the elements of $H_1$:
\[ 
\left\{ \begin{array}{llll}
H_1^{12} = \frac{d}{dt_2} H_0^{21} -\dfrac{d}{dt_2} H_0^{11} + \beta H_0^{11}, & \eqref{eq:SLH12} \\
H_1^{11}-H_1^{22} = i( \frac{d}{dt_2} H_1^{12} - \beta H_1^{12}) , & \eqref{eq:SLH11} \\
\frac{d}{dt_2} (H_1^{11} + H_1^{22}) = -  i \pi_{11} H_1^{12} + \beta (H_1^{11} - H_1^{22}), & \eqref{eq:SLH11+22}  \\
2i (H_1^{21})' = \frac{d^2}{dt_2^2} H_1^{11} - 2 \beta \frac{d}{dt_2}H_1^{11} & \eqref{eq:H021}'
\end{array} \right. \]
The last equation $\eqref{eq:H021}^\prime$ is obtained from \eqref{eq:H021} by substituting the index $0$ at $H_0^{ij}$
by the index $1$: $H_1^{ij}$.

Finally, we obtain that in the general case $H_{i+1}$ is derived from $H_i$ using a system of
similar equations.
\begin{equation} \label{eq:Hi+1HiFs} 
\left\{ \begin{array}{llll}
H_{i+1}^{12} & = \frac{d}{dt_2} H_i^{21} -\dfrac{d}{dt_2} H_i^{11} + \beta H_i^{11} , \\
H_{i+1}^{11}-H_{i+1}^{22} & = i( \frac{d}{dt_2} H_{i+1}^{12} - \beta H_{i+1}^{12}) , \\
\frac{d}{dt_2} (H_{i+1}^{11} + H_{i+1}^{22}) & = -  i \pi_{11} H_{i+1}^{12} + \beta (H_{i+1}^{11} - H_{i+1}^{22}), \\
2i \frac{d}{dt_2} H_{i+1}^{21} & = \frac{d^2}{dt_2^2} H_{i+1}^{11} - 2 \beta \frac{d}{dt_2}H_{i+1}^{11}. 
\end{array} \right. \end{equation}
from where we see that Markov moments are defined up to initial conditions for $n_0=2$ elements. Notice that
$p_0=1$ in this case.

Let us also demonstrate theorem \ref{thm:gamma*DE} for the Sturm Liouville parameters from definition
\ref{def:SLpar}. We will take the simplest case $n=1$ and as a result the transfer function is
\[ S(\lambda,t_2) = I - \dfrac{1}{\lambda + z} C(t_2) B(t_2) \sigma_1 =
I - \sum\limits_{i=0}^\infty \dfrac{1}{\lambda^{i+1}} (-z)^{i} C(t_2) B(t_2) \sigma_1.
\]
We have already seen (in \eqref{eq:PiBeta}) that $\gamma_*$ is necessarily of the form
\[ \gamma_*(t_2) = \bbmatrix{-i(\beta'-\beta^2) & -\beta \\ \beta & i}
\]
for a real valued function $\beta(t_2)$ on $\mathrm I$. In this case, $N=0$ which means that the first moment is
a multiple of the zero moment: $H_1=-z H_0$. Since the vessel parameters are constant the differential ring
$\boldsymbol{\mathcal R}=\mathbb C$ is trivial. Using the formulas developed in \cite{bib:SLVessels}
we shall obtain that $\tau=\exp{\int\beta}$ satisfies
\[ (\dfrac{d}{dt_2} - k)(\dfrac{d}{dt_2} - \bar k)(\dfrac{d}{dt_2} + k)(\dfrac{d}{dt_2} + \bar k) \tau = 0.
\]
for $k=\sqrt{-iz}$, which may be rewritten as a polynomial differential equation for $\beta$, after
inserting the formula for $\tau=\exp{\int\beta}$ and multiplying by $\tau^{-1}$.

\subsection{\label{sec:NLS}Non-Linear Shr\" odinger equation parameters}
\begin{defn}Non-Linear Shr\" odinger equation parameters are given by
\[ \begin{array}{lll}
 \sigma_1 = \bbmatrix{1&0\\0&1},\quad
\sigma_2 = \dfrac{1}{2} \bbmatrix{1&0\\0&-1},
\quad \gamma(x) =\bbmatrix{0 & 0 \\0 & 0}, \\
\gamma_*(x) = \bbmatrix{0&\beta(x) \\-\beta^*(x)& 0}
\end{array} \]
\end{defn}
Then the output and output compatibility conditions take the form of the classical non linear Schr\" odinger 
equation with the spectral parameter $i\lambda$
\[ \dfrac{\partial}{\partial x} u(x,\lambda) = (i \lambda A + Q(x)) u(x,\lambda),
\]
where 
\[ I=\sigma_1, \quad A=-\dfrac{1}{2}\bbmatrix{i&0\\0&-i} = -i \sigma_2, \quad Q(x) = -\gamma_*(x).
\]
We may perform the same calculations  for the moments as in the Sturm-Liouville case. Using formula
\eqref{eq:SMoments} and denoting n-th moment by $H_n(t_2) = \bbmatrix{H_n^{11} & H_n^{12}\\H_n^{21} & H_n^{22}}$ we shall find that the moments has to satisfy the following relations
\begin{equation} \label{eq:NLSmoments} 
\left\{ \begin{array}{llll}
H_{n}^{21} = - \dfrac{d}{dt_2}\big( H_{n-1}^{21}\big) - \beta^*(t_2) H_{n-1}^{11}, \\
H_{n}^{12} = \dfrac{d}{dt_2}\big( H_{n-1}^{12}\big) - \beta(t_2) H_{n-1}^{22}, \\
\dfrac{d}{dt_2}\big( H_n^{11}\big) = \beta(t_2) H_n^{21}, \\
\dfrac{d}{dt_2}\big( H_n^{22}\big) = -\beta^*(t_2) H_n^{12}.
\end{array} \right. \end{equation}
while the first moment $H_0$ is found from the linkage condition \eqref{eq:Link}
\[ H_0 = \bbmatrix{ H_0^{11} & \beta(t_2) \\ \beta^*(t_2) & H_0^{22}},
\]
and the entries $H_0^{11}, H_0^{22}$ are found using two last equations of \eqref{eq:NLSmoments}:
\[ \begin{array}{ll}
\dfrac{d}{dt_2}\big( H_0^{11}\big) = \beta(t_2) H_0^{21} = \beta(t_2) \beta^*(t_2) , \\
\dfrac{d}{dt_2}\big( H_0^{22}\big) = -\beta^*(t_2) H_0^{12} = -\beta^*(t_2) \beta(t_2) .
\end{array} \]

\section{\label{sec:InvKrein}Inverse problem at infinity for given vessel parameters}
In order to discuss the reconstruction of the transfer function from its moments, we will show first
that for almost arbitrary vessel parameters one can solve the $n_0$ moment equations, appearing in Theorem
\ref{thm:HiRstrct}.
Once a symmetric, equal to identity at infinity, function is constructed, we will use a Krein space
realization for it. Finally, by a counterexample, we will show that not all vessels parameters admit
Hilbert space (with positive $\mathbb X$) realizations.

We have seen that there is a differential equation \eqref{eq:Hn+1HnDE}, which connects the moments, and
means that the next moments is roughly speaking a second derivative of the previous one. If we start to estimate
its growth, we will end with factorial growth for the moments, because consecutive applying of derivative 
produces many terms due to Leibniz rule. Instead, we will assume that there is a "big" function $M(t_2)$ which
majorizes all the moments. Of course such an assumption narrows the class of function, to which the realization 
theorem \ref{thm:hnijBound} is applicable. On the other hand, the algebraic equations \eqref{eq:Moments} connects
the real/ imaginary parts of the moments with all the previous ones. Again, we found assumptions, which will enable
us to handle the growth of such elements.

From theorem \ref{thm:S=S-1*} it follows that we have to demand the symmetry condition \eqref{eq:Moments} only at
$t_2^0$. Examining the equation \eqref{eq:Moments} and supposing that $\|H_n(t_2^0)\| \leq f_n C^{n+1}$,
we will obtain that
\begin{multline*}
\| H_{n+1}\sigma_1^{-1} + (-1)^i \sigma_1^{-1} H_{n+1}^* \| \leq 
\sum\limits_{j=0}^{n} \| H_{n-j}\|\| \sigma_1^{-1}\|\| H_j^* \| \leq \\
\leq \| \sigma_1^{-1}\| \sum\limits_{j=0}^{n} f_{n-j} C^{n-j+1} f_j C^{j+1} =
\| \sigma_1^{-1}\| \sum\limits_{j=0}^{n} f_{n-j} f_j C^{n+2}
\end{multline*}
So, if the series $\{ f_n \}$ has the property
\begin{equation} \label{eq:fn}
 f_{n+1} = \|\sigma_1^{-1}\| \sum\limits_{j=0}^{n} f_{n-j} f_j ,
\end{equation}
we will obtain that
\[ \| H_{n+1}\sigma_1^{-1} + (-1)^i \sigma_1^{-1} H_{n+1}^* \| \leq 
\| \sigma_1^{-1}\| \sum\limits_{j=0}^{n} f_{n-j} f_j C^{n+2} \leq f_{n+1} C^{n+2}.
\]
One can find generating function $F(x)= \sum f_i x^i$. Multiplying \eqref{eq:fn} by $x^{n+1}$ and
summing, we will obtain that
\begin{multline*} x F(x)^2 = F(x) - F(0) \Rightarrow F(x)=\dfrac{1 - \sqrt{1-4 x F(0)}}{2x}
= \sum\limits_{n=0}^\infty f_n x^n
\end{multline*}
from where it follows that for $F(0)=1$
\[ f_n \sim \dfrac{4^n }{\sqrt{\pi n}(2n-1)}
\]
Let $M(t_2)$ be a positive function with positive derivatives on $[t_2^0,t_2]$,
which satisfies the following differential equation
\begin{multline} \label{eq:DM}
\dfrac{C}{3\sqrt{2}} \dfrac{d}{dt_2} M(t_2) = \\
= \| K_0 \| M(t_2) + \| K_1\| \dfrac{d}{dt_2} M(t_2) + \|K_2\| \dfrac{d^2}{dt_2^2} M(t_2) + 
\|K_3\| \dfrac{d^3}{dt_2^3} M(t_2)
\end{multline}
then, if $H^{(k)}_n(t_2) \leq f_n C^{n+1} M^{(k)} (t_2)$ we will obtain that
\begin{multline*} \| \dfrac{d}{dt_2} H_{n+1} \| \leq \\
\leq f_n C^{n+1}(\| K_0 \| \| H_n \| + \| K_1\|\| \dfrac{d}{dt_2} H_n\| + \|K_2\|\| \dfrac{d^2}{dt_2^2} H_n\|+ 
\|K_3\|\| \dfrac{d^3}{dt_2^3} H_n \|) \\
\leq f_n C^{n+1}(\| K_0 \| M(t_2) + \| K_1\| \dfrac{d}{dt_2} M(t_2) + \|K_2\| \dfrac{d^2}{dt_2^2} M(t_2) + 
\|K_3\| \dfrac{d^3}{dt_2^3} M(t_2)) \\
= \dfrac{f_n}{3\sqrt{2}} C^{n+2} M' (t_2) \leq f_{n+1} C^{n+2} M' (t_2) 
\end{multline*}
and by induction on $k$, further differentiating the equation \eqref{eq:Hn+1HnDE} and using triangle inequality and
additional assumption $\| K_i^{(k)}(t_2) \| \leq \dfrac{d^k}{dt_2^k} \| K_i(t_2)\|$, we will obtain that
\[ \| H^{(k)}_{n+1}(t_2) \| \leq \dfrac{f_n}{3\sqrt{2}} C^{n+2} M^{(k)} (t_2)\leq f_{n+1} C^{n+2}M^{(k)} (t_2).
\]

Using these ideas we obtain the following
\begin{thm} \label{thm:hnijBound}
Suppose that the vessels parameters and the norm are chosen so that $K_0, K_1, K_2,K_3$ satisfy for all $k\geq 1$
\begin{equation} \label{eq:AssmDki} \| K_i^{(k)}(t_2) \| \leq \dfrac{d^k}{dt_2^k} \| K_i(t_2) \|, \quad i=0,1,2,3.
\end{equation}
Suppose also that $M(t_2)$ is a solution of the differential equation \eqref{eq:DM} with the initial value $I$,
which majorizes $H_0(t_2)$ and all its derivatives: $\|H_0^{(k)}(t_2)\| \leq f_0 C M(t_2)$ on $[t_2^0,t_2]$. Then if
at each step the imaginary (for even $n$) and real (for odd $n$) part of $H_{n+1}\sigma_1^{-1}$, which is uniquely
determined from the first line in \eqref{eq:All} also satisfies 
\begin{equation} \label{eq:SecPart_t20}
 \| H_{n+1}(t_2^0)\sigma_1^{-1} - (-1)^i \sigma_1^{-1} H_{n+1}^*(t_2^0) \| \leq f_{n+1} C^{n+2},
\end{equation}
then the series 
\[ I - \sum\limits_{n=0}^\infty \frac{1}{\lambda^{i+1}} H_n(t_2) \sigma_1
\]
determines an analytic at $\lambda=\infty$ function $S(\lambda,t_2)$ 
not necessary $\sigma_1(t_2)$-contractive function (i.e. condition 
\eqref{eq:Scont} of proposition \ref{prop:PropS} does not necessary hold), which maps solutions of 
\eqref{eq:InCC} to solutions of \eqref{eq:OutCC} and which is symmetric (i.e. satisfies \eqref{eq:Symmetry}).
Radius of convergence around $\lambda=\infty$ is at least $\dfrac{1}{4C}$.
\end{thm}
\textbf{Proof:} Let us write down all the ideas, preceding the claim of this theorem.
We will prove the theorem, using induction on $n$. The basis of the induction follows from the assumptions, namely.
\[ \| H^{(k)}_0 \| \leq f_0 C M(t_2). \]
Suppose that $\|H_i^{(k)}(t_2)\| \leq f_i C^{i+1} M^{(k)}(t_2)$ for $i=0,1,\ldots,n$, then
from \eqref{eq:Hn+1HnDE} it follows that
\begin{multline*} \| \dfrac{d}{dt_2} H_{n+1} \| \leq \\
\leq f_n C^{n+1}(\| K_0 \| \| H_n \| + \| K_1\|\| \dfrac{d}{dt_2} H_n\| + \|K_2\|\| \dfrac{d^2}{dt_2^2} H_n\|+ 
\|K_3\|\| \dfrac{d^3}{dt_2^3} H_n \|) \\
\leq f_n C^{n+1}(\| K_0 \| M(t_2) + \| K_1\| \dfrac{d}{dt_2} M(t_2) + \|K_2\| \dfrac{d^2}{dt_2^2} M(t_2) + 
\|K_3\| \dfrac{d^3}{dt_2^3} M(t_2)) \\
= \dfrac{f_n}{3\sqrt{2}} C^{n+2} M' (t_2) \leq f_{n+1} C^{n+2} M' (t_2) 
\end{multline*}
and by induction on $k$, we obtain that for $k\geq 1$, 
\begin{multline*}  
\|H_{n+1}^{(k)}(t_2)\| = \\
= \| \sum_{i=0}^k \binom{k}{i} \{ 
K^{(i)}_0 H_n^{(k-i-1)} + K^{(i)}_1 H_n^{(k-i)} + K^{(i)}_2 H_n^{(k-i+1)} + K^{(i)}_3 H_n^{(k-i+2)}\} \\
\leq  \sum_{i=0}^k \binom{k}{i} 
\{ \| K^{(i)}_0 \| \|H_n^{(k-i-1)}\| + \|K^{(i)}_1 \|\|H_n^{(k-i)}\| + \|K^{(i)}_2 \|\|H_n^{(k-i+1)}\| + \|K^{(i)}_3\| \|H_n^{(k-i+2)} \| \} \\
\text{ using assumption \ref{eq:AssmDki} } \\
\leq  \sum_{i=0}^k \binom{k}{i} 
\{ \| K_0 \|^{(i)} \|H_n^{(k-i-1)}\| + \|K_1 \|^{(i)} \|H_n^{(k-i)}\| + \|K_2 \|^{(i)}\|H_n^{(k-i+1)}\| + \|K_3\|^{(i)} \|H_n^{(k-i+2)} \| \} \\
\leq f_n C^{n+1} \sum_{i=0}^k \binom{k}{i} 
\{ \| K_0 \|^{(i)} M^{(k-i-1)} + \|K_1 \|^{(i)} M^{(k-i)} + \|K_2 \|^{(i)} M^{(k-i+1)} + \|K_3\|^{(i)} M^{(k-i+2)} = \\
= f_n C^{n+1} ( \dfrac{C}{3\sqrt{2}} M')^{(k-1)} \leq f_{n+1} C^{n+2} M^{(k)}(t_2).
\end{multline*}
\begin{multline*}
\| H_{n+1}(t_2^0) \sigma_1^{-1} + (-1)^i \sigma_1^{-1} H_{n+1}^*(t_2^0)  \| \leq 
\sum\limits_{j=0}^{n} \| H_{n-j}(t_2^0) \|\| \sigma_1^{-1}\|\| H_j^*(t_2^0)  \| \leq \\
\leq \| \sigma_1^{-1}\| \sum\limits_{j=0}^{n} f_{n-j} C^{n-j+1} f_j C^{j+1} =
\| \sigma_1^{-1}\| \sum\limits_{j=0}^{n} f_{n-j} f_j C^{n+2} \leq f_{n+1} C^{n+2}
\end{multline*}
From \eqref{eq:SecPart_t20} and this inequality it follows that
\begin{multline*}
\| H_{n+1}(t_2^0) \sigma_1^{-1} \| 
\leq \dfrac{1}{2} \| H_{n+1}(t_2^0) \sigma_1^{-1} + (-1)^i \sigma_1^{-1} H_{n+1}^*(t_2^0)  \| + \\
+ \dfrac{1}{2} \| H_{n+1}(t_2^0)\sigma_1^{-1} - (-1)^i \sigma_1^{-1} H_{n+1}^*(t_2^0) \| \leq \\
\leq \dfrac{1}{2} (f_{n+1} C^{n+2} + f_{n+1} C^{n+2}) = f_{n+1} C^{n+2}.
\end{multline*}
Then
\begin{multline*}
\| H_{n+1}(t_2) \sigma_1^{-1} \| = 
\| H_{n+1}(t_2^0) \sigma_1^{-1} + \int_{t_2^0}^{t_2} \| H'_{n+1}(y) \| dy \leq \\
\leq  f_{n+1} C^{n+2} + \int_{t_2^0}^{t_2} f_{n+1} C^{n+2} M'(y) dy = \\
= f_{n+1} C^{n+2} + f_{n+1} C^{n+2} (M(t_2) - I) 
= f_{n+1} C^{n+2} M(t_2).
\end{multline*}
From these estimates we obtain that 
\[ H_n(t_2) \leq f_n C^{n+1} M(t_2)
\]
and as a result its radius of convergence is
\[ \dfrac{1}{R} = \lim\limits_{n\rightarrow\infty} \sqrt[n]{H_n(t)2} \leq
\lim\limits_{n\rightarrow\infty} \sqrt[n]{f_n C^{n+1} M(t_2)} = 4 C.
\]
because $f_n \sim \dfrac{4^n}{\sqrt{\pi n}}$ and the proof is finished.
\qed

Using results in \cite{bib:KreinReal} we can easily show using Cayley transform that 
there exists a Krein space realization for the just constructed function $S(\lambda,t_2)$.
\begin{thm} Let $S(\lambda,t_2)$ be an analytic at infinity function with value $I$ there symmetric
with respect to the imaginary axis (i.e. satisfies \eqref{eq:Symmetry}), then there exists a realization
\[ S(\lambda,t_2) = I - B(t_2)^* \mathbb X^{-1}(t_2) (\lambda I - A_1)^{-1} B(t_2) \sigma_1(t_2)
\]
such that the Lyapunov equation holds
\[ A_1 \mathbb X(t_2) + \mathbb X(t_2) A_1^* + B(t_2)^* \sigma_1(t_2) B(t_2) = 0.
\]
\end{thm}
\textbf{Proof:} The idea behind this theorem is that one can construct a Krein space realization at the given
time $t_2^0$ and then obtain a vessel with an inner Krein space $\widetilde{\mathcal K}$
instead of $\mathcal H$, using the formulas in theorem \ref{thm:gamma*}. Indeed, suppose that
there is constructed a Krein space realization at $t_2^0$ (which will be done in \eqref{eq:St20Krein}):
\[  S(\lambda,t_2^0) = I - B_0^* \mathbb X_0^{-1} (\lambda I - A_1)^{-1} B_0 \sigma_1(t_2^0).
\]
Define next $B(t_2)$ as the solution of \eqref{eq:DBSol}
\[ \frac{d}{dt_2} [B(t_2) \sigma_1(t_2)] +
A_1 B(t_2) \sigma_2(t_2) + B(t_2) \gamma(t_2) = 0,\quad B(t_2^0) = B_0 \]
and $\mathbb X(t_2)$ as a solution of \eqref{eq:DXSol}
\[ \frac{d}{dt_2} \mathbb X(t_2) = B(t_2) \sigma_2(t_2) B(t_2)^*,
\quad \mathbb X(t_2^0) = \mathbb X_0. \]
and define $\gamma^1_*(t_2)$ using the Linkage condition \eqref{eq:Link} for which we shall denote
by $\Phi_*^1(\lambda,t_2,t_2^0)$ the fundamental matrix of \eqref{eq:OutCC}. Then we will obtain that
\[ Y(\lambda,t_2) = \Phi_*^1(\lambda,t_2,t_2^0) S(\lambda,t_2^0) \Phi(\lambda,t_2,t_2^0).
\]
On the other hand, in theorem \ref{thm:hnijBound} we proved that one can construct a 
function $S(\lambda,t_2)$ with suitably chosen initial values, i.e. with $S(\lambda,t_2^0)$ such that 
\[ S(\lambda,t_2) = \Phi_*(\lambda,t_2,t_2^0) S(\lambda,t_2^0) \Phi(\lambda,t_2,t_2^0)
\]
for the given $\gamma_*(t_2)$. It is remained to notice that
\[ Y(\lambda,t_2) S^{-1}(\lambda,t_2) = \Phi_*^1(\lambda,t_2,t_2^0) \Phi_*(\lambda,t_2,t_2^0)
\]
is entire and identity at infinity, i.e. is $I$, which means that
\[ Y(\lambda,t_2) = S(\lambda,t_2)
\]
and from here immediately follows (as at the second part of \ref{thm:gamma*}) that $\gamma_*(t_2) = \gamma_*^1(t_2)$ and
we have obtained a Krein space realization for the function $S(\lambda,t_2)$ constructed in theorem
\ref{thm:hnijBound}.

It remains to show that there always exists a realization at $t_2^0$.
Define a function $Q(\lambda)$ using Cayley transform, which satisfies
\begin{equation} \label{eq:SfromQ}
S(-i\lambda,t_2^0) = ( I + \dfrac{i}{2} Q(\lambda) \sigma_1(t_2^0)) (I - \dfrac{i}{2} Q(\lambda)\sigma_1(t_2^0))^{-1}.
\end{equation}
Actually, this function is given by 
\[ Q(\lambda) = 2 i \sigma_1^{-1}(t_2^0) (I - S(-i\lambda,t_2^0)) (I + S(-i\lambda,t_2^0))^{-1}
\]
and is well-defined at the neighborhood of infinity with value $0$ there. Then from a simple equality, resulting from the
symmetry condition \eqref{eq:Symmetry} considered with $-i\lambda$ instead of $\lambda$
\begin{multline*}
(I-S^*(-i\bar\lambda,t_2^0))\sigma_1^{-1}(t_2^0)(I + S(-i\lambda,t_2^0)) =
\\
= - (I+S^*(-i\bar\lambda,t_2^0))\sigma_1^{-1}(t_2^0)(I - S(-i\lambda,t_2^0))
\end{multline*}
it follows that $Q(\lambda)^* = Q(\bar\lambda)$ and $Q(\lambda,t_2)$ is zero at the neighborhood of
$\lambda=\infty$. Thus \cite[Theorem 3]{bib:KreinReal} $Q(\lambda,t_2)$ admits the following Krein space realization
\[ Q(\lambda) = \Gamma^+ (\widetilde A - \lambda I)^{-1} \Gamma
\]
where for a Krein space $\widetilde{\mathcal K}$, there is a bounded self-adjoint operator
$\widetilde A\in L(\widetilde{\mathcal K})$, and a bounded
$\Gamma\in L(\mathcal E, \widetilde{\mathcal K})$. So, we can insert this realization formula into
\eqref{eq:SfromQ} and simplify:
\[ \begin{array}{lllllll}
S(-i\lambda,t_2^0) & = ( I + \dfrac{i}{2} Q(\lambda) \sigma_1(t_2^0)) (I - \dfrac{i}{2} Q(\lambda)\sigma_1(t_2^0))^{-1} = \\
& = ( 2I - I + \dfrac{i}{2} Q(\lambda) \sigma_1(t_2^0)) (I - \dfrac{i}{2} Q(\lambda)\sigma_1(t_2^0))^{-1}  = \\
& = -I + 2 (I - \dfrac{i}{2} Q(\lambda)\sigma_1(t_2^0))^{-1}  = \\
& = -I + 2 (I - \dfrac{i}{2} \Gamma^+ (\widetilde A - \lambda I)^{-1} \Gamma \sigma_1(t_2^0))^{-1}
\end{array} \]
There is a simple formula \cite{bib:bgr} for evaluating the inverse of a matrix in a realized form:
\[ (I - \dfrac{i}{2} \Gamma^+ (\widetilde A - \lambda I)^{-1} \Gamma \sigma_1(t_2^0))^{-1} =
I + \dfrac{i}{2} \Gamma^+ (\widetilde A^\times - \lambda I)^{-1} \Gamma \sigma_1(t_2^0),
\]
where $\widetilde A^\times = \widetilde A - \dfrac{i}{2} \Gamma \sigma_1(t_2^0) \Gamma^+$. So, the last formula becomes
\begin{equation} \label{eq:SPreKrein} \begin{array}{lllllll}
S(-i\lambda,t_2) & = -I + 2 (I - i \Gamma^+ (\widetilde A - \lambda I)^{-1} \Gamma \sigma_1(t_2^0))^{-1} = \\
& = -I + 2(I+\dfrac{i}{2} \Gamma^+ (\widetilde A^\times - \lambda I)^{-1} \Gamma \sigma_1(t_2^0)) = \\
& = I + i \Gamma^+ (\widetilde A^\times - \lambda I)^{-1} \Gamma \sigma_1(t_2^0) = \\
& = I -  \Gamma^+ (i \widetilde A^\times - i\lambda I)^{-1} \Gamma \sigma_1(t_2^0) 
\end{array} \end{equation}
Let us define $A_1 = - i \widetilde A^\times$ then we obtain that
\begin{equation} \label{eq:LyapPre}
 \begin{array}{lllllll}
A_1 + A_1^{+} & = -i \widetilde A^\times + i (\widetilde A^\times)^+ = \\
& = -i(\widetilde A - \dfrac{i}{2} \Gamma \sigma_1(t_2^0) \Gamma^+) + i (\widetilde A^+ + \dfrac{i}{2} \Gamma \sigma_1(t_2^0) \Gamma^+) = \\
& =  -\Gamma \sigma_1(t_2^0) \Gamma^+,
\end{array} \end{equation}
since $\widetilde A$ is selfadjoint. In the Krein space, the Krein space adjoint may be represented by an
invertible selfadjoint operator $\widetilde{\mathbb X}$ in (the Hilbert space sense),
when one considers the Krein space with its inner 
product $[\cdot,\cdot]$ as a Hilbert space with its initial norm $\langle\cdot,\cdot\rangle$:
\[ [u,v] = \langle \widetilde{\mathbb X} u,v\rangle, \forall u,v \in \widetilde{\mathcal K}.
\]
Then the formulas for Krein space adjoint becomes:
\[ A_1^{+} = \widetilde{\mathbb X}^{-1} A_1^* \widetilde{\mathbb X}, \quad
\Gamma^+ = \widetilde{\mathbb X} \Gamma^*
\]
and the last equation \eqref{eq:LyapPre}
\[ \begin{array}{lllllll}
A_1 +  A_1^{+} = - \Gamma \sigma_1(t_2^0) \Gamma^+ \Leftrightarrow \\
A_1 + \widetilde{\mathbb X}^{-1} A_1^* \widetilde{\mathbb X} =
 \Gamma \sigma_1(t_2^0) \Gamma^* \widetilde{\mathbb X}  \Leftrightarrow \\
A_1 \widetilde{\mathbb X}^{-1} + \widetilde{\mathbb X}^{-1} A_1^* = \Gamma \sigma_1(t_2^0) \Gamma^* 
\end{array}  \]
which is the Lyapunov equation \eqref{eq:Lyap} at $t_2^0$ after defining
$B_0 = \Gamma$ and $\mathbb X_0 = \widetilde{\mathbb X}^{-1} $. Notice that Lyapunov equation
\eqref{eq:Lyap} will also hold for all $t_2$ as a result of lemma \ref{lemma:LyapREdund}.
Moreover, substituting $\lambda$ for $-i\lambda$ in \eqref{eq:SPreKrein}, we also obtain
\begin{equation} \label{eq:St20Krein}
S(\lambda,t_2^0) = I - B_0^* \mathbb X_0^{-1} (\lambda I - A_1)^{-1} B_0 \sigma_1(t_2^0).
\end{equation}
\qed 

It turns out that there are examples for which there is no a realization with strictly positive $\mathbb X$.
For, example considering the Non-Linear Schr\" odinger equation parameters, defined in section \ref{sec:NLS}
for \textit{constant} $\beta$
\[ \gamma_* = \bbmatrix{0 & \beta\\ \beta^* & 0}
\]
we can easily find that the fundamental matrices for the input and output LDEs are (for convenience
$t_2^0=0$)
\begin{eqnarray*}
\Phi(\lambda,t_2) = 
\bbmatrix{e^{\frac{\lambda}{2}t_2} & 0\\0 & e^{\frac{-\lambda}{2}t_2}}, \\
\Phi_*(\lambda,t_2) =
\bbmatrix{ \dfrac{-\beta\beta^*}{2\Gamma} [
\dfrac{E}{\Gamma-\frac{\lambda}{2}} + \dfrac{E^{-1}}{\Gamma+\frac{\lambda}{2}}] &
\dfrac{\beta}{2\Gamma}[E-E^{-1}] \\
\dfrac{-\beta^*}{2\Gamma}[E-E^{-1}] &
\dfrac{(\Gamma-\frac{\lambda}{2})E}{2\Gamma} + \dfrac{(\Gamma+\frac{\lambda}{2})E^{-1}}{2\Gamma}  }
\end{eqnarray*}
where
\[ \Gamma = \sqrt{\frac{\lambda^2}{4}-\beta\beta^*}, \quad
E = \exp(\Gamma t_2).
\]
Let us search for a constant in $t_2$ matrix
 $S(\lambda,0) = \bbmatrix{a(\lambda) & b(\lambda) \\ c(\lambda) & d(\lambda)}$ such
that the constructed $S(\lambda,t_2) = \Phi_*(\lambda,t_2) S(\lambda,0) \Phi^{-1}(\lambda,t_2)$ 
satisfies conditions of the proposition \ref{prop:PropS}. 
Considering the identity at infinity requirement, we come to the conclusion that
\[ c(\lambda) = \dfrac{-\beta^*}{\Gamma+\frac{\lambda}{2}} a(\lambda), \quad
b(\lambda) = \dfrac{-\beta}{\Gamma+\frac{\lambda}{2}} d(\lambda),
\]
where the square root at the Right Half Plane (RHP) is taken so that at the infinity its argument
coincides with the argument of $\dfrac{\lambda}{2}$. In this case the expression $\Gamma+\frac{\lambda}{2}$
is analytic in the whole $\mathbb C$, except for a cut from $-\sqrt{2} |\beta|$ to $\sqrt{2} |\beta|$. In 
other words this can not be defined as an analytic function at the right half plane. On the other hand,
if it were possible to realize this function with $A_1$, satisfying Lyapunov equation \eqref{eq:Lyap}
with positive $\mathbb X$, this would mean that the spectrum of $A_1$ is at the left hand side, since
$\sigma_1=I$ and is a positive matrix:
\[ A_1 \mathbb X+ \mathbb X A_1^* = - B^* B \leq 0.
\]
As a result $S(\lambda,0)$ would have spectrum only at the left hand
side of the complex plane, contradicting the previous argument.

Finally, we present in this section a lemma, which says that one can always define moments so that the corresponding
Pick matrices are positive, but as the example above shows they will not always define an analytic function 
(the radius of convergence may be zero)
\begin{lemma} \label{lemma:Pn+1Moments}
Suppose that equations \eqref{eq:Moments}, \eqref{eq:HiHi+1} hold so that 
the eigenvalues of the Pick matrix $\mathbb P_n$ are strictly positive on a finite interval $I$. 
Suppose that initial values are chosen for $H_{2n+1}$, then there exists an initial values for 
$H_{2n+2}$ such that all eigenvalues of $P_{n+1}$ are positive on arbitrary finite interval around $t_2^0$.
\end{lemma}
\textbf{Proof:} Let us use Sylvester criterion for positive definiteness of a matrix.
Denote the entries of the matrix $(-1)^n H_{2n} \sigma_1^{-1}$ as follows
\[ (-1)^n H_{2n} \sigma_1^{-1} = [h_{ij}] =
\bbmatrix{h_{11} & h_{12} & \ldots & h_{1p} \\
h_{21} & h_{22} & \ldots & h_{2p} \\
\vdots & \vdots & \ddots & \vdots \\
h_{p1} & h_{p2} & \ldots & h_{pp} \\}, \quad
D_n = [d_{ij}]
\]
and suppose that $\mathbb P_n>0$ on a fixed interval $\mathcal I$ including $t_2^0$. We will consider main minors of the
matrix $\mathbb P_{n+1}$. First minors come from $\mathbb P_n$ and as a result are positive, the next one is
\[ \bbmatrix{ \mathbb P_n & C_n^1 \\ (C_n^1)^* & h_{11} + d_{11}},
\]
where $C_n^1$ is the first column of the matric $C_n$. Its determinant, using a formula for block matrix is
\[ \det \bbmatrix{ \mathbb P_n & C_n^1 \\ (C_n^1)^* & h_{11} + d_{11}} =
\det \mathbb P_n \det (h_{11} + d_{11} - (C_n^1)^* \mathbb P_n^{-1} C_n^1).
\]
In order to obtain here a positive result, we have to demand that (since $\det \mathbb P_n > 0$ and since
$h_{11} + d_{11}$ is real valued)
\[ \Re h_{11} > (C_n^1)^* \mathbb P_n^{-1} C_n^1 - \Re d_{11}.
\]
Using the continuity and taking $h_{11}(t_2^0)$ big enough we obtain that the same inequality will hold on the
interval $\mathcal I$. Consider now the next main minor
\[ \bbmatrix{ \mathbb P_n & C_n^1 & C_n^2 \\ (C_n^1)^* & h_{11} + d_{11} & h_{12} + d_{12} \\
(C_n^2)^* & h_{21} + d_{21} & h_{22} + d_{22} }.
\]
Its determinant can be expressed as
\[ \begin{array}{lll}
\det \bbmatrix{ \mathbb P_n & C_n^1 & C_n^2 \\ (C_n^1)^* & h_{11} + d_{11}& h_{12} + d_{12} \\
(C_n^2)^* & h_{21} + d_{21}& h_{22} + d_{22} } = \\
= \det \bbmatrix{ \mathbb P_n & C_n^1 \\ (C_n^1)^* & h_{11}+ d_{11}} \times  \\
 \times \det \big(h_{22} + d_{22} - \bbmatrix{(C_n^2)^* & h_{21}+ d_{21}} 
\bbmatrix{ \mathbb P_n & C_n^1 \\ (C_n^1)^* & h_{11}+ d_{11}}^{-1} \bbmatrix{C_n^2 \\ h_{12}+ d_{12} }\big)
\end{array} \]
from where we obtain that the necessary condition for its positiveness is
\[ \Re h_{22} > \bbmatrix{(C_n^2)^* & h_{21}+ d_{21}} 
\bbmatrix{ \mathbb P_n & C_n^1 \\ (C_n^1)^* & h_{11}+ d_{11}}^{-1} \bbmatrix{C_n^2 \\ h_{12}+ d_{12} } - \Re d_{22}
\]
And so on for all the real parts of diagonal elements $h_{ii}$. It remains to remark that the conditions
\eqref{eq:Moments}, \eqref{eq:HiHi+1} do not involve real parts of $h_{ii}$ and they can be determined from
the positive definiteness of $P_{n+1}$. \qed

\section{\label{sec:NevPick}Generalized Nevanlinna-Pick interpolation problem}
We address the two kinds of Nevanlinna-Pick interpolation problems
\ref{prob:NPclassic} and \ref{prob:NPhard}, defined in the introduction. The first one is solved using
the notion of positive pairs and for the second one we present a criterion in terms of the of the initial data.
\subsection{\label{sec:LFTPairs}Linear Fractional Transformations in terms of intertwining positive pairs}
Suppose that we are given a data of the NP interpolation problem \ref{prob:NPclassic}.
Following the notations of corollary \ref{cor:TTheta} let us write
\[ B(t_2) = \bbmatrix{-\xi_1(t_2)^* S(w_1,t_2)^* & \xi_1(t_2)^* \\ \vdots & \vdots \\ -
\xi_n(t_2)^* S(w_n,t_2)^* & \xi_n(t_2)^*},~~~~~ A_1 = \operatorname{diag}[-w_1^*,\ldots,-w_n^*].
\]
Let us denote by capital Greek letters the following vessel parameters
\begin{eqnarray*} \Sigma_1(t_2) = \bbmatrix{-\sigma_1(t_2) & 0\\ 0 & \sigma_1(t_2)}= J, \\
\Sigma_2(t_2) = \bbmatrix{\sigma_2(t_2) & 0\\ 0 & \sigma_2(t_2)}, \\
\Gamma(t_2) = \bbmatrix{\gamma_*(t_2) & 0\\ 0 & \gamma(t_2)},
\end{eqnarray*}
then simple calculations show that $B(t_2)$ satisfies \eqref{eq:DBSol}
\[\frac{d}{dt_2} [B(t_2) \Sigma_1(t_2)] + A_1 B(t_2) \Sigma_2(t_2)
+ B(t_2) \Gamma(t_2) = 0,\quad B(t_2^0) = B. \]
Suppose that $\mathbb X(t_2)>0$ is a solution of
\begin{eqnarray*} 
A_1 \mathbb X(t_2) + \mathbb X(t_2) A_1^* + B(t_2) \Sigma_1(t_2) B^*(t_2) = 0, \\
\dfrac{d}{dt_2}\mathbb X(t_2) =  B(t_2) \Sigma_2(t_2) B^*(t_2),
\end{eqnarray*}
which is always possible if $\Re w_i\neq 0$ for each $i=1,\ldots,n$. Then the following collection
\[ \mathfrak{V} = \{ A_1, B(t_2), \mathbb X(t_2); \Sigma_1(t_2), \Sigma_2(t_2), \Gamma(t_2), \Gamma_*(t_2);
\mathbb C^{2n}; \mathcal E\oplus\mathcal E\}.
\]
is a vessel for $\Gamma_*(t_2)$ defined from the linkage condition \eqref{eq:Link}
\begin{eqnarray*}
\Gamma_*(t_2)  =  \Gamma(t_2) + \Sigma_1(t_2) B(t_2)^* \mathbb X^{-1}(t_2) B(t_2) \Sigma_2(t_2) - \\
- \Sigma_2(t_2) B(t_2)^* \mathbb X^{-1}(t_2) B(t_2) \Sigma_1(t_2).
\end{eqnarray*}
Transfer function of the vessel $\mathfrak{V}$ is
\[  \Theta(\lambda,t_2) = I_{2p} - B(t_2)^* \mathbb X^{-1}(t_2) (\lambda I_n - A_1)^{-1} B(t_2) \Sigma_1(t_2),
\]
which is in $\CI(\bU_*\oplus\bU,\widetilde\bU)$ for
\[ \widetilde \bU = \lambda \Sigma_2(t_2) - \Sigma_1(t_2) \dfrac{d}{dt_2} + \Gamma_*(t_2).
\]
If we denote further the decomposition of $\Theta(t_2)$ as
\[ \Theta(\lambda,t_2) = \bbmatrix{\Theta_{11}(\lambda,t_2) & \Theta_{12}(\lambda,t_2) \\
\Theta_{21}(\lambda,t_2) & \Theta_{22}(\lambda,t_2)}
\]
then if one defines $S_0(\lambda,t_2)$ so that
\[ (\Theta_{11} + S_0 \Theta_{21})^{-1}(\Theta_{12} + S_0 \Theta_{22}) = S,
\]
then the function $S(\lambda,t_2)$ usually does not intertwine solutions of LDEs with spectral parameter $\lambda$.

Instead, we define
\[ W(\lambda,t_2) = \bbmatrix{W_1(\lambda,t_2) & W_2(\lambda,t_2)}
\]
so that
\[ W_1 \Theta_{11} + W_2 \Theta_{21} = I_p, \quad W_1 \Theta_{12} + W_2 \Theta_{22} = S,
\]
then the following lemma holds
\begin{lemma}\label{lemma:WPos} The pair of functions $W(\lambda,t_2)$ is in $\CI(\widetilde \bU,\bU_*)$ and 
$\dfrac{W(\lambda,t_2)^* \Sigma_1(t_2) W(\mu,t_2)}{\bar\lambda+\mu} \geq 0$ on the domain of analyticity of $W(\lambda,t_2)$.
\end{lemma} 
\textbf{Proof:} Since $\Theta(\lambda,t_2)$ is invertible for all $\lambda$ out of the spectrum of $A_1$,
an element of $\widetilde U$ is of the form $\Theta(\lambda,t_2)\bbmatrix{y_\lambda(t_2)\\u_\lambda(t_2)}$,
where $y_\lambda(t_2)$, $u_\lambda(t_2)$ satisfy \eqref{eq:OutCC} and \eqref{eq:InCC} respectively. Then
\begin{eqnarray*}
W(\lambda,t_2) \Theta(\lambda,t_2)\bbmatrix{y_\lambda(t_2)\\u_\lambda(t_2)} = \\
= \bbmatrix{W_1(\lambda,t_2) & W_2(\lambda,t_2)}  \bbmatrix{\Theta_{11}(\lambda,t_2) & \Theta_{12}(\lambda,t_2) \\
\Theta_{21}(\lambda,t_2) & \Theta_{22}(\lambda,t_2)} \bbmatrix{y_\lambda(t_2)\\u_\lambda(t_2)}= \\
= \bbmatrix{W_1 \Theta_{11} + W_2 \Theta_{21} & W_1 \Theta_{12} + W_2 \Theta_{22}}\bbmatrix{y_\lambda(t_2)\\u_\lambda(t_2)} = \\
= \bbmatrix{I & S(\lambda,t_2)}\bbmatrix{y_\lambda(t_2)\\u_\lambda(t_2)} = y_\lambda(t_2) + S(\lambda,t_2) u_\lambda(t_2) \in \bU_*, \\
\end{eqnarray*}
since $y_\lambda(t_2)$ and $S(\lambda,t_2) u_\lambda(t_2)$ are in $\bU_*$.

From the formula $W \Theta = \bbmatrix{I & S}$, it follows that
\[ W(\lambda,t_2) = \bbmatrix{I & S(\lambda,t_2)} \Theta^{-1}(\lambda,t_2).
\]
Consequently, the expression $\dfrac{W(\lambda,t_2)^* \Sigma_1(t_2) W(\mu,t_2)}{\bar\lambda+\mu}$ 
considered on the domain of analyticity of $W(\cdot,t_2)$ becomes
\begin{eqnarray*}\dfrac{W(\lambda,t_2) \Sigma_1(t_2) W^*(\mu,t_2)}{\lambda + \bar\mu} = \\
= \bbmatrix{I & S(\lambda,t_2)} \dfrac{\Theta^{-1}(\lambda,t_2) \Sigma_1(t_2) \Theta^{-1*}(\mu,t_2)}{\lambda + \bar\mu}
\bbmatrix{I\\S^*(\mu,t_2)}.
\end{eqnarray*}
Since $\Theta(\lambda,t_2)$ is a transfer function of a conservative vessel, its inverse
is a transfer function too and satisfies
\[ \dfrac{\Theta^{-1}(\lambda,t_2) \Sigma_1(t_2) \Theta(\mu,t_2)^{-1*}-\Sigma_1(t_2)}{\lambda + \bar\mu} \geq 0
\]
and consequently,
\[ \begin{array}{lll}
\dfrac{W(\lambda,t_2) \Sigma_1(t_2) W^*(\mu,t_2)}{\lambda + \bar\mu} & 
\geq \dfrac{1}{\lambda+\bar\mu}\bbmatrix{I & S(\lambda,t_2)} \Sigma_1(t_2) \bbmatrix{I\\S^*(\mu,t_2)} \\
& \hspace{-3cm} \geq \dfrac{1}{\lambda+\bar\mu}\bbmatrix{I & S(\lambda,t_2)} \bbmatrix{-\sigma_1(t_2) & 0 \\ 0 & \sigma_1(t_2) } \bbmatrix{I\\S^*(\mu,t_2)}  \\
& \hspace{-3cm} \geq \dfrac{S(\lambda,t_2) \sigma_1(t_2)S^*(\mu,t_2) - \sigma_1(t_2)}{\lambda+\bar\mu} \\
& \hspace{-3cm} \geq 0
\end{array} \]
by the properties of transfer functions for vessels.
\qed

As a consequence of this theorem, we define 
\begin{defn} A pair of functions
\[ W(\lambda,t_2) = \bbmatrix{W_1(\lambda,t_2) & W_2(\lambda,t_2)}
\]
such that the matrix $W(\lambda,t_2)$ has full rank 
is called \textbf{positive} if the conditions of lemma \ref{lemma:WPos} hold:
\[ W(\lambda,t_2)\in \CI(\widetilde \bU,\bU_*), \quad W(\lambda,t_2) J W(\lambda,t_2) \geq 0 \text{ on } \mathbb C_+
\]
\end{defn}
Previous arguments result in the following "partial solution" of the Nevanlinna-Pick interpolation problem
\ref{prob:NPclassic}:
\begin{cor} There exists a one-to-one correspondence between solution of the Nevanlinna-Pick interpolation problem
\ref{prob:NPclassic} and the positive pairs $W(\lambda,t_2)$, which are analytic at the interpolation points.
\end{cor}

\subsection{Nevanlinna-Pick interpolation problem \ref{prob:NPhard}}
Using the previous results, we present a criterion for solvability of the Problem \ref{prob:NPhard}.
\begin{thm} Given $\mathbb C^{p\times p}$-valued
functions $\sigma_1$, $\sigma_2$, $\gamma$, an interval $\mathrm I$ and
$N$ quadruples $<t_2^j,w_j,\xi_j,\eta_j>$ where $t_2^j\in\mathrm I$,
$w_j\in\mathbb C_+$, $\xi_j, \eta_j \in \mathbb C^{1\times p}$
$j=1,\ldots N$, and assume that the corresponding matrices 
$\widetilde{\mathbb X_i}>0$. Then
there exists a solution of the Nevanlinna-Pick problem
\ref{prob:NPhard}, 
i.e. there exists a function $S\in\RCI$ satisfying
$S(w_i,t_2^i)\xi_i = \eta_i$
if and only if there exists $n\in\mathbb N$, matrices $A_0^i,\mathbb
X_0^i,
\in\mathbb C^{(n-1)\times (n-1)}$ with $\mathbb X_0^i>0$,
$B_0^i\in\mathbb C^{p\times (n-1)}, V_{ij}\in\mathbb C^{n\times n}$ such that for
$A_i,B_i,\mathbb X_i$ defined by \begin{eqnarray}
\label{eq:Bi}\hspace{5mm} B_i & = & \bbmatrix{ B_0^i \\ \eta_i - \xi_i} \\
\label{eq:Xi}\hspace{5mm} \mathbb X_i & = & \bbmatrix{\mathbb X_0^i & 0 \\ 0 & \widetilde{\mathbb X}_i  } \\
\label{eq:Ai}\hspace{5mm} A_i & = & \bbmatrix{ A_0^i & \dfrac{B_0^i \sigma_1(t_2^i) \xi_i^*}{\widetilde{\mathbb X}_i} \\ -\eta_i\sigma_1(t_2^i) (B_0^i)^* (\mathbb X_0^i)^{-1}  & - w_i^* - \dfrac{\eta_i\sigma_1(t_2^i)(\eta_i^*-\xi_i^*)}{\widetilde{\mathbb X}_i} }
\end{eqnarray}
it holds that
\begin{enumerate}
    \item $A_i = V_{ij} A_j V_{ij}^{-1}$,
\item $\oint(\lambda I - A_i)^{-1} B_i \sigma_1^{-1}(t_2^j) \Phi^{-1}(\lambda,t_2^j,t_2^i) d\lambda = V_{ij} B_j$.
\end{enumerate}
and the matrix $X(t_2)$,
\[
\mathbb X(t_2)=\mathbb X_i+\int_{t_2^i}^{t_2}B_i(y)\sigma_2(y)(B_i(y))^*dy
\]
is invertible on the interval $\mathrm I$.
\end{thm}
{\bf Proof:}
For each $t_2^i$ all the functions which satisfy $S(w_i,t_2^i)
\xi_i=\eta_i$ are of the form $T_{\Theta_i}(S^i_0(\lambda,t_2))$, provided
$\widetilde{\mathbb X}_i = \dfrac{ \xi_i\sigma_1(t_2^i)\xi_i^*  -
\eta_i\sigma_1(t_2^i)\eta_i^*}{w_i^* + w_i} > 0$. Here
\[ \Theta_i(\lambda) =
\bbmatrix{I + \dfrac{\eta_i^*\eta_i\sigma_1}{\widetilde{\mathbb X}_i (\lambda + w_i^*)} & \dfrac{\eta_i^*\xi_i\sigma_1}{\widetilde{\mathbb X}_i (\lambda + w_i^*)} \\
-\dfrac{\xi_i^*\eta_i\sigma_1}{\widetilde{\mathbb X}_i (\lambda + w_i^*)} & I - \dfrac{\xi^*\xi\sigma_1}{\mathbb X (\lambda + w^*)}}
\]
and $S_0^i\in\RSC$. Given a minimal realization
\[ S_0^i(\lambda) = I -  (B_0^i)^* (\mathbb X_0^i)^{-1}
(\lambda I - A_0^i)^{-1} B_0^i \sigma_1(t_2^i)
\]
of $S_0^i$, we shall obtain from formulas (\ref{eq:BS}),
(\ref{eq:XS}), (\ref{eq:AS}) the formulas (\ref{eq:Bi}),
(\ref{eq:Xi}), (\ref{eq:Ai}) in the theorem.
In view of Theorem \ref{thm:EqualS}, a necessary and sufficient
condition to obtain the same function $S(\lambda,t_2)$ for every $i$
is that there exist invertible
constant matrices $V_{ij}$ such that the operators
$A_i$ are similar. In other words there
must exist $V_{ij}$ such that $A_i = V_{ij} A_j V_{ij}^{-1}$. Moreover, the second part of theorem \ref{thm:EqualS} tells that additionally the
equality
\[ \oint(\lambda I - A_i)^{-1} B_i \sigma_1^{-1}(t_2^j) \Phi^{-1}(\lambda,t_2^j,t_2^i) d\lambda = V_{ij} B_j
\]
must hold.\qed\mbox{}\\


\end{document}